\documentclass[preprint,12pt]{elsarticle}

\usepackage{amssymb}
\usepackage{amsmath}

\usepackage{commath,amsmath,amsfonts}
\usepackage{amsthm}
\usepackage{graphicx}

\usepackage{subcaption}	

\usepackage{csvsimple}
\usepackage{fp} 

\usepackage{pgfplotstable}

\usepackage{booktabs}	

\usepackage{mathtools}

\usepackage{mathrsfs}	

\usepackage{hyperref}

\usepackage{algorithm} 		
\usepackage{algpseudocode} 	

\usepackage{color,soul}		

\usepackage{bm}	

\usepackage{gnuplot-lua-tikz} 

\usepackage{comment}

\usepackage[normalem]{ulem} 

\theoremstyle{plain}
\newtheorem{theorem}{Theorem}

\theoremstyle{definition}
\newtheorem{defn}{Definition}[section]

\theoremstyle{remark}
\newtheorem{remark}{Remark}

\theoremstyle{plain}
\newtheorem{result}{Numerical Result}

\begin{document}

\newcommand{\flow}{\Phi_H}
\newcommand{\flown}{\Phi_N}

\newcommand{\SM}{S}
\newcommand{\SMn}{\SM^{N}}
\newcommand{\NHIM}{\Lambda}

\newcommand{\FRM}{\mathcal{F}}

\newcommand{\innerMap}{\mathcal{F}}

\newcommand{\sm}{\sigma}	
\newcommand{\smn}{\sm^{N}}

\newcommand{\tm}{\tau}	

\newcommand{\gf}{\mathcal{L}}	

\newcommand{\gft}{\widetilde{\gf}}	

\newcommand{\NHIMn}{\Lambda^{N}}	
\newcommand{\NHIMc}{\Lambda_c}	
\newcommand{\NHIMcn}{\Lambda_c^{N}}	
\newcommand{\NHIMsec}{\Lambda_c^\Sigma}

\newcommand{\NHIMsecn}{\Lambda_c^{N,\Sigma}}	

\newcommand{\An}{\mathbb{A}^{N}}

\newcommand*\average[1]{\overline{#1}}

\newcommand{\domain}{\mathcal{A}}	

\newcommand{\posDomain}{\mathcal{A}_+}	
\newcommand{\negDomain}{\mathcal{A}_-}	
\newcommand{\neuDomain}{\mathcal{A}_0}	
\newcommand{\Tsecn}{T^{N,\Sigma}}

\newcommand{\domainloc}{\domain_{\text{loc}}}

\newcommand{\tin}{t_\mathrm{in}}
\newcommand{\tout}{t_\mathrm{out}}

\newcommand{\ee}{\mathrm{e}}	

\newcommand{\Anote}[1]{\marginpar{\textcolor{cyan}{\textbf{A: #1}}}} 

\newcommand{\Mnote}[1]{\relax} 
\newcommand{\AB}[1]{{#1}} 
\newcommand{\AR}[1]{\relax} 
\newcommand{\MB}[1]{{#1}} 
\newcommand{\MR}[1]{\relax} 
\newcommand{\PB}[1]{{#1}} 
\newcommand{\PR}[1]{\relax} 

\begin{frontmatter}



\title{Semi-analytic construction of global transfers between quasi-periodic orbits in the spatial R3BP}


\author[inst1]{Amadeu Delshams} 

\author[inst2]{Marian Gidea} 

\author[inst3]{Pablo Roldan\corref{cor1}} 
\ead{Pablo.Roldan@upc.edu}

\cortext[cor1]{Corresponding author}

\affiliation[inst1]{organization={Laboratory of Geometry and Dynamical Systems and IMTech, UPC},
	country={Barcelona}
    }

\affiliation[inst2]{organization={Department of Mathematical Sciences, Yeshiva University},
	city={New York},
	postcode={10016}, 
	state={NY},
	country={USA}}

\affiliation[inst3]{organization={Departament de Matem\`atiques, Universitat Polit\`ecnica de Catalunya (UPC)},
	city={Barcelona},
	postcode={08028}, 
	country={Spain}}

\begin{abstract}
	Consider the spatial restricted three-body problem, as a model for the motion of a spacecraft relative to the Sun-Earth system. We focus on the dynamics near the equilibrium point $L_1$, located between the Sun and the Earth. We show that we can transfer  the spacecraft  from a quasi-periodic orbit that is nearly planar relative to the ecliptic to a quasi-periodic orbit that has large vertical amplitude, at zero energy cost. (In fact, the final orbit has the maximum vertical amplitude that can be obtained through the particular mechanism that we consider. Moreover, the transfer can be made through any prescribed sequence of quasi-periodic orbits in between).

Our transfer mechanism is based on selecting trajectories homoclinic to a normally hyperbolic invariant manifold (NHIM) near $L_1$,
and then gluing them together. \MB{We present a theoretical result  establishing the existence of such transfer orbits, and we verify numerically its applicability to our model.}
We provide several explicit constructions of such transfers, and also develop an algorithm to design trajectories that achieve the \emph{shortest transfer time} for this particular mechanism.

The change in the vertical amplitude along a homoclinic trajectory can be described via the scattering map. We develop a new tool, the `Standard Scattering Map' (SSM), which is a series representation of the exact scattering map. We use the SSM to obtain a complete description of the dynamics along homoclinic trajectories.
The SSM can be used in many other situations, from Arnold diffusion problems to transport phenomena in applications.

\end{abstract}




\end{frontmatter}



\section{Introduction} \label{Sec:Introduction}
In this paper, we consider the spatial circular restricted three-body problem (RTBP for short), as a model for the motion of spacecraft relative to the Sun-Earth system.
We focus on the dynamics near the equilibrium point $L_1$ located between the Sun and the Earth. We show that we can transfer the spacecraft, at zero energy cost,  from a quasi-periodic orbit that is nearly planar relative to the ecliptic  to a quasi-periodic orbit of large vertical amplitude. \MR{That is, we can achieve a change in the vertical amplitude of the orbit of the spacecraft simply by choosing  suitable initial conditions and letting the gravitational fields of the Sun and the Earth drive the motion.} Moreover, we provide several explicit constructions of such trajectories, and also develop an algorithm to design  trajectories that achieve the shortest transfer time.  Our algorithm is flexible and  can be applied to other systems besides Sun-Earth.

\begin{figure}
	{\tiny \scalebox{0.9}{\input{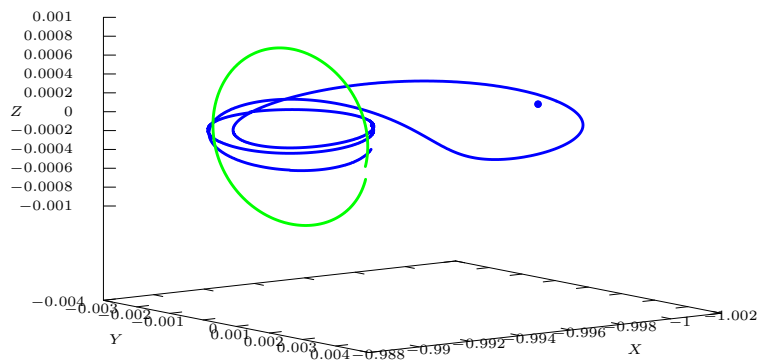}}}
	\caption{Initial segment (blue) and final segment (green) of the fastest transfer trajectory between a quasi-periodic orbit with small vertical amplitude (small $Z$-oscillation), and another with large vertical amplitude ($Z$-oscillation of amplitude $0.001$, roughly 150000 km). The complete trajectory (not shown here) is very complicated, and goes through many homoclinic jumps. This trajectory is computed in Section~\ref{sec:triple_dynamical_system}. (See also  Figures~\ref{fig:diffusion_shortest_path_TM} and \ref{fig:diffusion_orbit}).
	}
	\label{fig:initial_final_segments}
\end{figure}
	
For illustration, Figure~\ref{fig:initial_final_segments} shows the initial and final segments (in blue and green, respectively) of the fastest transfer trajectory between one quasi-periodic orbit with small vertical amplitude, and another with large vertical amplitude. 

\MR{The model that we consider  is a $3$-degrees of freedom Hamiltonian system. We construct trajectories that follow closely geometric structures that organize the dynamics. The main geometric object near $L_1$ is a center manifold  on which the dynamics is nearly integrable. More precisely, in a neighborhood of $L_1$ the Hamiltonian can be approximated by a high-order Birkhoff normal form, which is an integrable Hamiltonian.   In terms of the Birkhoff normal form, the center manifold is  represented by  a $4$-dimensional manifold, which  is foliated by a family of $2$-dimensional invariant tori. The sphere can be parametrized by a system of symplectic coordinates consisting of two action variables $(J_p, J_v)$ and two angle variables $(\phi_p,\phi_v)$, with each torus corresponding to a pair of fixed values of the two actions. The action $J_p$ describes the horizontal amplitude  (relative to the ecliptic) of an orbit lying on a $2$-dimensional torus, and $J_v$ describes the vertical amplitude.
Restricting to an energy level close to that of $L_1$ amounts to fixing the action variable $J_p$. This yields a   $3$-dimensional sphere which is filled with $2$-dimensional tori. Each torus is given by a fixed value of the remaining action variable $J_v$. A change in $J_v$ corresponds to a change in   the vertical amplitude of  the orbit. Since the tori are invariant, by using  only the `inner flow' restricted to the $3$-dimensional sphere, the vertical amplitude of orbits remains constant. In order to change the vertical amplitude, we need to use the `outer dynamics', described below.}

\MB{The model that we consider  is a $3$-degrees of freedom Hamiltonian system. We construct trajectories that follow closely geometric structures that organize the dynamics. The main geometric object near $L_1$ is a 4-dimensional center manifold  on which the dynamics is nearly integrable. More precisely, in a neighborhood of $L_1$ the Hamiltonian can be approximated by a high-order Birkhoff normal form, which is an integrable Hamiltonian.   In terms of the Birkhoff normal form, the center manifold can be parametrized by a system of symplectic coordinates consisting of two action variables $(J_p, J_v)$ and two angle variables $(\phi_p,\phi_v)$.
The action $J_p$ describes the horizontal amplitude  (relative to the ecliptic) of a quasiperiodic orbit, and $J_v$ describes its vertical amplitude.
Restricting to an energy level close to that of $L_1$ amounts to fixing the action variable $J_p$. This yields a   $3$-dimensional sphere which is foliated by $2$-dimensional tori, where each torus is uniquely determined by a fixed value of the remaining action variable $J_v$. 
 A change in $J_v$ corresponds to a change in   the vertical amplitude of  a quasi-periodic orbit. 
 Since the tori are invariant, by using  only the `inner flow' restricted to the $3$-dimensional sphere, the vertical amplitude of orbits remains constant. However, the tori do not separate the energy manifold, so there may exist trajectories that move from one torus to another, thus changing the vertical amplitude.  In order to obtain such trajectories, we need to use the `outer dynamics', described below.
}

The $3$-dimensional sphere is a normally hyperbolic invariant manifold (NHIM), and has stable and unstable manifolds which go around the Earth and intersect transversally along trajectories homoclinic to the NHIM. By \MB{carefully selecting consecutive}  homoclinic trajectories, we show that it is possible to achieve large changes in the vertical amplitude of orbits. \MB{Proving the existence of such trajectories is related to \textit{Arnold diffusion problem}   (see the discussion later).}

\MR{However, if we would pick at random which  homoclinic trajectories to follow, sometimes we will obtain a growth in the vertical amplitude,  and other times we will obtain a decay. In general, one expects that the values of $J_v$ will follow a stochastic process, e.g., a Brownian motion with drift \cite{CapinskiG23}. To consistently achieve a growth in the vertical amplitude, we need to carefully select which homoclinic trajectories to follow.}

The tool that allows us to systematically select suitable homoclinics at each step is the \emph{scattering map}. This is a map defined on the NHIM, which relates the past asymptotic of a homoclinic point to its future asymptotic.
The scattering map was introduced in \cite{DelshamsLS00,DelshamsLS06a,DelshamsLS06b} in the study of Arnold diffusion.
When restricted to a suitable $2$-dimensional Poincar\'e section, the scattering map turns out to be symplectic \cite{DLS08,DelshamsGR12}. An additional advantage that we exploit in our model is that we obtain two scattering maps, which give us more options in the selection of suitable homoclinics. \MR{In general,  Hamiltonian systems similar to the one we consider here exhibit multiple homoclinics and associated scattering maps.}

\MB{We provide a theoretical result (Theorem \ref{th:main})  that gives sufficient conditions on a inner dynamics  and two scattering maps defined on an annulus inside a NHIM (described in action-angle variables), to ensure the existence of true trajectories that start near the lower boundary of the annulus and end near its upper boundary.}

\MB{In our model, the action variable corresponds to the vertical amplitude $J_v$ and the annulus to some range of $J_v\in[J_1,J_2]$. Provided that the assumptions of Theorem \ref{th:main} are verified, the orbits of the Iterated Functions System (IFS) formed by the two scattering maps and the inner map that achieve the prescribed gain in the \MB{vertical amplitude}  yield true trajectories that achieve the same gain (up to a small error).} \AB{To distinguish the orbits of the IFS from true trajectories of the system, we will refer to  the former as  \emph{pseudo-orbits}.}

\AB{
Before going into more details, we now outline the methodology that we will use. The trajectories designed to change the \MB{vertical amplitude}  $J_v$ are guided by homoclinic orbits that depart and land, asymptotically, at points in the NHIM, depending on the scattering (outer) map used. The \MB{vertical amplitude} at landing can be higher or lower than the \MB{vertical amplitude} at departure. The inner dynamics of the NHIM allow us to travel \MR{internally} between points with the same \MB{vertical amplitude} to change, if necessary, the \MB{departure} points. With a simple simulation of these 3 dynamics---the 2 scattering maps and the inner map---we obtain pseudo-orbits with the desired \MB{vertical amplitude} gain.
}

\MR{In perturbative problems, the scattering map can be computed analytically via perturbation theory \cite{DLS08}. However, }

We note that the  problem that we consider is not perturbative. Until now, the only available methods for computing the scattering map in such a setting have been numerical \cite{CanaliasDMR06,DGR16}. However, purely numerical approaches are computationally intensive and offer little insight into the geometric structures determined by the scattering map.

The highlight of this paper is that we provide an analytical approximation of the scattering map for the spatial  circular RTBP. We describe the scattering map via a generating function depending on old and new variables. Then we approximate the generating function using a
Fourier-\AB{Polynomial interpolation}.
As it turns out, the numerical computation of the scattering map at a few points can be used to compute the coefficients of the Fourier-\AB{Polynomial interpolation}, up to some suitable order, and thus to obtain an analytical formula for the generating function, up to some small error. The outcome of this approximation is referred to as the \emph{Standard Scattering Map} (SSM). It is  given explicitly as a perturbation of an integrable twist map. As such,  the phase space of the scattering map is organized by KAM tori, elliptic islands, hyperbolic periodic orbits and their stable and unstable manifolds, and resonant zones. See Figure~\ref{fig:phase_port_SM}. The rich geometric structure unveiled by the analytical approximation of the scattering map was not available through previous approaches. Similar computations of the generating function of the scattering map and of its phase space were obtained purely analytically in \emph{uncoupled} pendulum-rotor systems subject to small perturbations of a special type \cite{DelshamsS17,DS18}.
However, as the unperturbed pendulum-rotor systems considered in those papers are uncoupled, the phase shift phenomenon~\cite{DMR08} does not take place, and the unperturbed scattering maps are just the identity. In particular, they are not twist maps, which makes the dynamics different from those considered in this paper.

The main application of the Standard Scattering Map is that it offers an explicit method to find pseudo-orbits, that is, orbits of the iterated function system consisting of  scattering maps and the inner map (induced by the inner flow on the Poincar\'e section), along which 
the vertical amplitude $J_v$ grows consistently. The method is versatile, in the sense that one can  choose the starting and ending points of such pseudo-orbits.
As mentioned earlier, in this paper we compute two scattering maps, and we compare them in terms of the fastest trajectory to achieve the desired change in vertical amplitude.

It is important to note that the orbits of the scattering maps are not equivalent to true trajectories of the system. Rather, we can approximate  a segment of a homoclinic trajectory by a concatenation of a finite orbit of the inner dynamics, followed by an application of a scattering map, followed by another finite orbit of the inner dynamics.
The map that described  this  concatenation of orbits is referred to as the \emph{transition map}. The transition map is a map on the NHIM, and each application of the transition map approximates a segment of a homoclinic orbit to the NHIM;  see Section~\ref{sec:transition}. Since we have constructed two scattering maps, we have two corresponding transition maps. To obtain approximate trajectories that change the orbital vertical amplitude,  we consider the iterated function system consisting of the two transition maps and the inner map, and we search for optimal trajectories.



\PB{Note the order in which our methodology proceeds. First and foremost, we obtain the Standard Scattering Map, which models the asymptotic homoclinic dynamics in a comprehensive and efficient manner, forming the backbone of all subsequent computations. Second,  we use the SSM  to validate the hypothesis of Theorem~\ref{th:main} and conclude the existence of diffusion. However, to design pseudo-orbits it is necessary to work with finite-time segments of trajectories. Third, from the SSM we derive the transition map, which approximates \emph{finite-time} homoclinic segments. The transition map is used to design pseudo-orbits, which can finally be refined into true orbits.}

To find optimal pseudo-orbits, we leverage the classic Dijkstra algorithm for finding shortest paths in a graph. A surprising finding is that for an optimal orbit of the iterated function system, rather than always selecting a transition map that grows $J_v$,  sometimes we must select a transition map that decreases $J_v$, in order to arrive to a place where the next application of a transition map yields a large increase in $J_v$. Another surprising finding is that an optimal pseudo-orbit involves very few applications of the inner map.

\AB{
More quantitatively, in Section~\ref{Sec:DriftOrbits} we divide the phase space into a grid of $30\times30=900$ small cells, and in Section~\ref{sec:triple_dynamical_system} we obtain a diffusion pseudo-orbit that uses just $17$ iterations of the scattering/transition maps plus intermediate iterations of the inner map. \MR{The points of this pseudo-orbit are those that can later be used to find a real trajectory that shadows them.} We emphasize that only using the image of $7$ tori we can have an approximation of the phase space of the scattering map, good enough to be able to design with several algorithms the appropriate diffusion pseudo-orbits. The use of more invariant tori and of a more refined grid does not present additional technical difficulties, and would not modify in a qualitative way the phase space of the scattering maps represented in Figure~\ref{fig:phase_port_SM}, so the optimal paths would not change substantially.
}


\PB{These optimal pseudo-orbits can be refined into true orbits either through theoretical shadowing results or by applying appropriate numerical methods, such as parallel shooting.}

Our construction described so far is based on approximating the Hamiltonian near $L_1$ by  a Birkhoff normal form.  Since the approximation is quite accurate, the true dynamics associated to the original Hamiltonian follow closely the Birkhoff normal form dynamics. In particular,  the trajectories of the true inner dynamics stay close to invariant tori, and the scattering map for the true dynamics is close to the scattering map derived from the normal form approximation. This implies that there exist diffusing trajectories---that change the orbital vertical amplitude by a significant amount---for the original Hamiltonian system.

Our results are related to the Arnold diffusion problem for Hamiltonian systems, claiming that integrable Hamiltonian systems subjected to small perturbations of generic type have `diffusing orbits' along which the action variable changes by an amount independent of the smallness of the perturbation \cite{Arnold64}. \MR{Arnold illustrated this phenomenon for an uncoupled pendulum-rotor system subject to small perturbations of  special type.} \MB{Arnold proposed a mechanism of diffusion based on} transition chains of tori, which are sequences of invariant tori with consecutive heteroclinic connections between consecutive tori. Arnold conjectured: ``I believe that this mechanism of instability is applicable to the general case (for example, to the problem of three bodies)". This conjecture
has witnessed significant progress in recent years, including \cite{ChierchiaG94,bolotin1999unbounded,DelshamsLS00,Treschev02c,Mather04,Treschev04, DelshamsLS06a,DelshamsLS06b,Piftankin2006,GelfreichT2008,DelshamsHuguet2009,ChengY09,Mather12,KaloshinZ15,bernard2016arnold,
cheng2019variational,
Treschev12,GideaL17,Gelfreich2017,gidea_marco_2017,kaloshin2020arnold}. Some of the progress has been geared towards proving Arnold diffusion in concrete models, under explicit, verifiable conditions on the perturbation. This direction opened up the possibility of implementing Arnold's mechanism of diffusion in applications, such as to Celestial Mechanics. Notably, some papers, including \cite{DelshamsKRS19,guardia2023degenerate,clarke2022inner}, succeeded in providing analytical proofs of Arnold diffusion in some models of the three- and four-body problem.
These papers rely on perturbative methods, and therefore they need to assume that certain parameters (such as ratios of the  masses of the bodies, or ratios of the semi-major axes of the orbits)  are very  small, in fact much smaller than those observed in solar systems. Another line of arguments combine analytical methods with numerical computation (including computer assisted proofs) to show Arnold diffusion in models with realistic parameters, see, e.g. \cite{capinski2016arnold,CapinskiG23,FiguerasHaro24}.

\MB{Our paper follows this latter approach. In the case of the Sun-Earth system we implement  Arnold's mechanism of transition chains of tori via analytical tools and numerical methods, obtaining orbits that follow the transition chain and change their vertical amplitude. The system that we consider is \emph{a priori chaotic}, that is, it contains a hyperbolic basic set (see, e,g., \cite{Gelfreich2017}).}
Since the system \MR{we consider} is non-perturbative, we do not obtain a change that is independent of some small parameter. However, we argue that the change is as large as possible for the geometric mechanism at play.
Specifically, we obtain a change in the action that spans approximately 70\% of the region where primary homoclinic orbits exist. (See Section~\ref{Sec:DomainGlobalSM}).

The construction in this paper can be potentially adapted to astrodynamics applications. Sometimes, a satellite (which typically  carries little fuel) ends up on a wrong orbit,
and  one tries to correct the orbit by exploiting  the gravity of Earth, Sun, Moon as much as possible,  and firing the satellites' thrusters as little as possible; see, e.g.,~\cite{belbruno2004capture}. While our methodology to change the vertical amplitude of a satellite orbit may be too slow from a practical point of view, by combining zero-cost geometric routes with small thrusts, one may be able to design useful trajectories; see, e.g.,~\cite{GomezKLMMR2001}.

\PB{This paper does not focus on trajectory optimization. Nevertheless, in our setting, we can estimate both the duration ($32$ years) and the total required thrust ($271$ m/s) for the pseudo-orbit that exhibits the fastest growth in vertical amplitude. See Section~\ref{sec:triple_dynamical_system} and Remark~\ref{rem:thrust}.}

Moreover, our methodology can be applied to build transfers involving two or more NHIMs connected by heteroclinic orbits, e.g., the NHIM around $L_1$ and that around $L_2$~\cite{DMR08, BarcelonaHaroMondelo24}.

\section{Setup} \label{Sec:Setup}

\subsection{The Spatial Circular RTBP} \label{Sec:RTBP}

\begin{figure}
	\centering
	\includegraphics[width=0.5\linewidth]{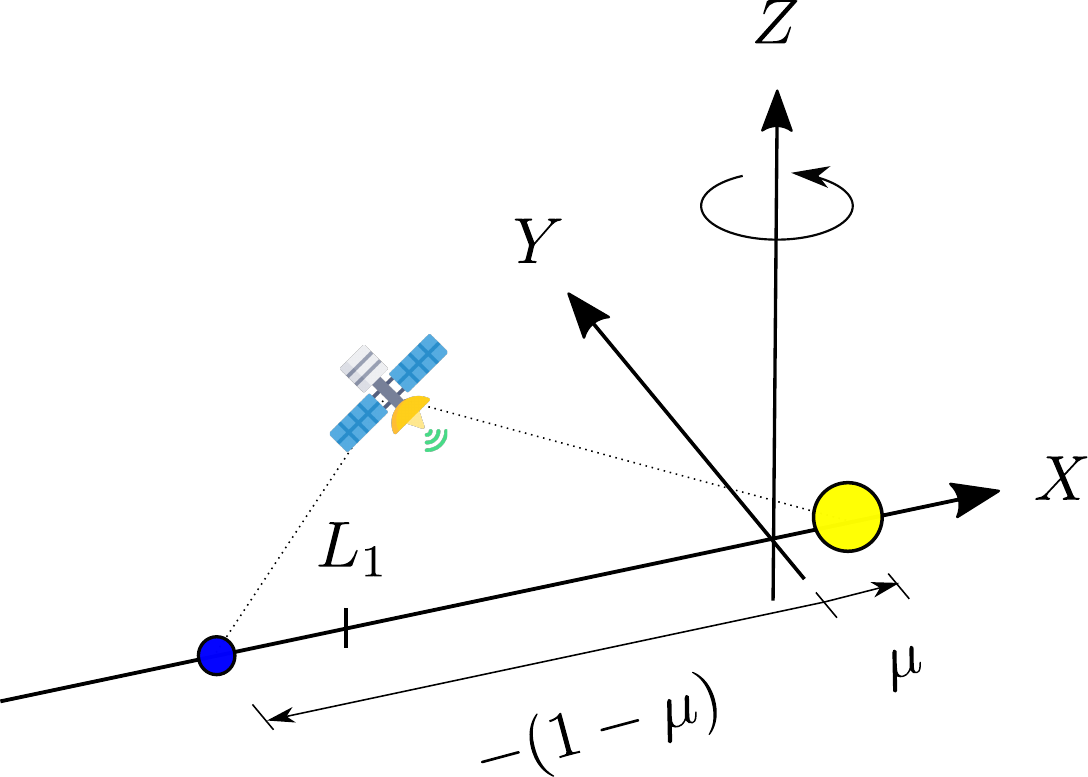}
	\caption{Schematic location of the Sun (in yellow), the Earth (in blue), and the equilibrium point $L_1$.}
	\label{fig:SR3BP}
\end{figure}

\MR{We consider the spatial circular RTBP as a model for the motion of a satellite under the gravitational influence of the Sun and the Earth.
In this model, two heavy bodies (referred to as primaries) move in the same plane along circular orbits about their common center of mass, while a  third, infinitesimal body (referred to as secondary) moves in space under the gravitational influence of the heavy bodies,  without affecting their orbits.
It is convenient to use a co-rotating frame $XYZ$ whose origin $O$ is set at the center of mass of the system, such that the orbits of the primaries lie in the $XY$-plane, and the vertical component of the motion of the secondary is given by the  $Z$-coordinate. The units are normalized so that the masses of the primaries are  $\mu_1$ and $\mu_2$, with $\mu_1+\mu_2=1$,  the  distance between the primaries is   $1$, the period of the motion of the primaries is $2\pi$,   and the gravitational constant is  $G=1$.
Denoting the smaller mass by $\mu=\mu_2$ and the
larger one by  $1-\mu$, the larger body is located on the $X$-axis to the right of $O$ at
$P_1=(\mu,0,0)$, and the smaller body is located on the $X$-axis to the left of $O$ at $P_2=(\mu-1,0,0)$ (see Figure \ref{fig:SR3BP}).}
\MB{We consider the spatial circular RTBP as a model for the motion of an infinitesimal body (e.g. satellite) under the gravitational influence of two massive bodies,
referred to as primaries (e.g. Sun and the Earth).
We describe the problem in normalized units, relative to a  co-rotating frame $XYZ$, such that the larger primary  of mass $1-\mu=\mu_1$ is located at $P_1=(\mu,0,0)$, the smaller primary of mass $\mu$  is located at $P_2=(\mu-1,0,0)$,  and the vertical component of the motion of the infinitesimal body is given by the  $Z$-coordinate (see Figure \ref{fig:SR3BP}).}

In the case of  the Sun-Earth system    $\mu=3.040423398444176 \times 10^{-6}$.

The  motion of the infinitesimal body relative to these coordinates is given by the autonomous system of equations:
\begin{equation}\label{eq:RTBP}
\begin{split}
   \ddot{X}-2\dot{Y}	&=\Omega_X,\\
   \ddot{Y}+2\dot{X}	&=\Omega_Y,\\
   \ddot{Z}		            &=\Omega_Z,
\end{split}
\end{equation}
where the effective potential $\Omega$   is given by
\begin{equation*}
   \Omega=\frac{1}{2}(X^2+Y^2)+\frac{1-\mu}{r_1}+\frac{\mu}{r_2},
\end{equation*}
with $r_1,r_2$ representing  the distances from the secondary to the larger and the
smaller primary, respectively:
\begin{align*}
   r_1&=((X-\mu)^2+Y^2+Z^2)^{1/2},\\
   r_2&=((X-\mu+1)^2+Y^2+Z^2)^{1/2}.
\end{align*}
The phase space is $6$-dimensional. 

The system has an integral of motion (referred to as  the Jacobi
integral) given by:
\[ C=2\Omega-(\dot X^2+\dot Y^2+\dot Z^2). \]

Equivalently, the equations \eqref{eq:RTBP} can be described as a $3$-degree-of-freedom, autonomous Hamiltonian system given by the Hamiltonian function:
\begin{equation}\label{eq:Hamiltonian}
H=\frac{1}{2}(P_X^2+P_Y^2+P_Z^2)+Y P_X-X P_Y-\frac{1-\mu}{r_1}-\frac{\mu}{r_2},
\end{equation}
where  $X$, $Y$, $Z$ are the  generalized coordinates,
$P_X=\dot X-Y$, $P_Y=\dot Y+X$, $P_Z=\dot Z$ are the generalized momenta,
and the symplectic form is:
\[dP_X\wedge dX+dP_Y\wedge dY+dP_Z\wedge dZ.\]
\MB{We denote by $\flow$ the flow of \eqref{eq:Hamiltonian}.}

The Hamiltonian function and the Jacobi integral are equivalent integrals of motion, since
\[H=-{C}/{2}.\]

As the Hamiltonian (and, equivalently, the Jacobi integral) is preserved along the solutions of the system,
each  trajectory lies on a $5$-dimensional  energy manifold $M_c$ corresponding to some energy level $h=-c/2$,
that is, \[M_c=\{H(P_X,P_Y,P_Z,X,Y,Z)=-c/2\}.\]
\MR{or, equivalently, on the level surface of the Jacobi integral}
\begin{equation*}\label{eq:energySurface}
\MR{
M_c=\{C(X,Y,Z,\dot X, \dot Y, \dot Z)=c\}.
}
\end{equation*}

\MR{The system has 5 equilibrium points, denoted $L_1, \dotsc, L_5$.
Here we adopt the convention that $L_1$ is located between the primaries (see Figure~\ref{fig:SR3BP}).
The equilibria $L_1,L_2,L_3$ are of Saddle $\times$ Center $\times$ Center -- linear stability type,
and the  equilibria $L_4,L_5$ are of Center $\times$ Center $\times$ Center -- linear stability type  (provided that $\mu$ is less than Routh's critical value $\mu_{\textrm{cr}}$, which is always the case for the planets in our solar system).
A general reference for the RTBP  is~\cite{Sze}.}

\MB{In this paper we focus on the dynamics around the equilibrium point $L_1$, which is located between the primaries,   and is of  Saddle $\times$ Center $\times$ Center -- linear stability type (see Figure~\ref{fig:SR3BP}).
Let us denote by 
\[\pm\lambda,\pm i\nu_p, \pm i\nu_v, \textrm{  with }\lambda,\nu_p,\nu_v\textrm{ real and positive},\]  the eigenvalues of the linearized system  near $L_1$.
The quantities $\lambda,-\lambda$ represent the exponential expansion and contraction hyperbolic rates,    while  $\nu_p,\nu_v$ represent the frequencies of the planar and vertical components of the motion, respectively.}
\AB{In the case of the Sun-Earth system, $\lambda=2.5326591740529678$, $\nu_p=2.0864535642231075$ and $\nu_v=2.0152106629966386$. }

\MB{By the Center Manifold
Theorem (see e.g., \cite{JorbaMasdemont99}), there exists a $4$-dimen-sional invariant
center manifold $\NHIM=W^{\mathrm{c}}(L_1)$   that is tangent at $L_1$ to the generalized
eigenspace  corresponding to $\pm i\nu_p, \pm i\nu_v$.
The manifold $\NHIM$ is a normally hyperbolic invariant manifold (NHIM) for the flow $\flow$.
If we fix the energy level (or, equivalently the Jacobi constant),
the restriction of the center manifold $W^{\mathrm{c}}(L_1)$ to the energy level $\{H=-{c}/{2}\}$ is a $3$-dimensional sphere 
\[ \NHIMc = \NHIM \cap M_c,\]
which is a NHIM for the flow $\flow$ restricted to the energy manifold $M_c$.} 

\subsection{Local Dynamics around \texorpdfstring{$L_1$}{L1}}

\MB{The Hamiltonian  $H$ can be expanded about $L_1$, via a symplectic coordinate change
\[(P_X,P_Y,P_Z,X,Y,Z)\mapsto ( y_h, J_p,  J_v, x_h, \phi_p,  \phi_v) ,\]
 as a Birkhoff normal form
\[ \mathcal{H}^{(N)} = H_2 + \mathcal{Z}^{(N)} + \mathcal{R}^{(N)}, \]
where $H_2$ denotes the quadratic part of $H$ expressed in terms of the new variables,
$\mathcal{Z}^{(N)}$ is  a polynomial of degree $N$ that Poisson-commutes with $H_2$,
and the remainder $\mathcal{R}^{(N)}$ is small in a neighborhood of $L_1$ (more precisely, of the order of the $(N+1)$-th power of the distance to $L_1$).
The coordinates $x_h,y_h$ represent the hyperbolic directions of motion.
\MR{We refer to $J_p,  J_v$ as the action variables, and to $\phi_p,  \phi_v$ as the angle variables.}
The \emph{action variables} $(J_p, \phi_p)$ correspond to the \textit{planar} component of the motion, while the \emph{angle variables} $(J_v, \phi_v)$ corresponds to the \textit{vertical} component.
A useful reference for the  derivation of such a Birkhoff normal form is  \cite{Jorba99}.
}

\MR{The above power series expansion is not convergent in general, but only asymptotically convergent.} \MB{
In~\cite{DGR16}, we performed the computation of this normal form for the spatial circular RTBP  up to order $N=16$.
We will use this computation in this paper.
}

\MB{The truncated Birkhoff normal form depends only on   $x_hy_h$, $J_p$, $J_v$,
\begin{equation}\label{eq:HN}
\mathcal{H}^{(N)}_{\textrm{trunc}}=H_2 + \mathcal{Z}^{(N)} = \lambda x_h y_h + \nu_p J_p + \nu_v J_v + \mathcal{Z}^{(N)}(x_hy_h, J_p, J_v),
\end{equation}
which are integrals of motion for $\mathcal{H}^{(N)}_{\textrm{trunc}}$.
We denote by $\flown$ the flow associated to \eqref{eq:HN}.
}

\MB{The truncated Birkhoff normal form $\mathcal{H}^{(N)}_{\textrm{trunc}}$ represents  an approximation of the original Hamiltonian $H$, and therefore the invariant objects
for $\mathcal{H}^{(N)}_{\textrm{trunc}}$ give approximations of the corresponding invariant objects of $H$. 
Truncating the normal form at order $N=16$ ensures that the error between the true dynamics in the NHIM and the dynamics induced by  the Birkhoff normal form is within machine precision. See Section~\ref{subsec:NFprecision}.
In particular, we consider the center manifold $\NHIMn$ for the truncated Birkhoff normal $H^{(N)}_{\textrm{trunc}}$, which represents an approximation
of $\NHIM$.  In the action-angle coordinates, we have that
\begin{equation}
\label{eq:CenterManifold} \NHIMn=\{( y_h, J_p,  J_v, x_h, \phi_p,  \phi_v)\,|\,y_h=0,x_h=0, J_p,  J_v\geq0,\ \phi_p,  \phi_v \in\mathbb{T}\},
\end{equation}
}
\AB{is parameterized by the symplectic action-angle coordinates $(J_p,J_v,\phi_p,\phi_v)$.
}%
\MB{
This has a symplectic structure given by the form
\[dJ_p\wedge d\phi_p+dJ_v\wedge d\phi_v.\]
}%
\MR{The center manifold $\NHIMn$ is foliated by a family of $2$-dimensional invariant tori obtained by fixing  the values of the actions $J_p=\bar{J}_p>0$,  $J_v=\bar{J}_v>0$. Moreover, there}\MB{The planar Lyapunov periodic orbits correspond to $J_v=0$ and $J_p>0$, the vertical Lyapunov periodic orbits to $J_p=0$ and $J_v>0$, and the libration point $L_1$ to $J_p=J_v=0$, so these coordinates are degenerate when some action $J_{p,v}=0$, as it is usual with polar coordinates.
}

\MR{By comparison, the center manifold  $\NHIM$ corresponding to the full Hamiltonian $H$ contains a large family of  $2$-dimensional KAM invariant tori, with the gaps between tori of exponentially small size with respect to the actions $J_p,J_v$ (see~\cite{DG96} for quantitative estimates).  Therefore, the trajectories in the center manifold for $H$ closely follow those for the truncated Birkhoff normal form approximation~\eqref{eq:CenterManifold}, for the times considered in this paper.}

\MB{For a fixed energy $\mathcal{H}^{(N)}_{\textrm{trunc}}=h=-{c}/{2}$,  $J_p$ is given  as an implicit function of $J_v, \phi_p, \phi_v, c$.
The corresponding NHIM for the truncated system is given by
\[ \NHIMcn = \NHIMn \cap \left\{\mathcal{H}^{(N)}_{\textrm{trunc}}=-{c}/{2} \right\}, \]
}%
\AB{and is therefore parameterized by the 3 variables $(J_v, \phi_p, \phi_v)$.
}

\MB{
For each fixed
value of the vertical action $J_v = \bar{J}_v$, there is a unique invariant torus for the Birkhoff normal form
\[
T_{\bar{J}_v}^{N}  =  \NHIMcn \cap \left\{ J_v = \bar{J}_v\right\}
\]
}%
\AB{which is therefore parameterized by the 2 angles $(\phi_p, \phi_v)$.
}
\MB{Thus, the NHIM $\NHIMcn$ is foliated by a family of $2$-dimensional invariant tori. 
By comparison, the NHIM $\NHIMc$ corresponding to the full Hamiltonian 
$H$  contains a KAM family of $2$-dimensional  invariant tori, with gaps between the tori that are exponentially small in the action 
 (see [41] for quantitative estimates).
The existence of the NHIM for the full Hamiltonian has been rigorously proven for the planar RTBP in~\cite{CapinskiRoldan12}. }

\MB{Each of the  objects $\NHIMn$, $\NHIMcn$, $T^{N} _{\bar{J}_v}$, have associated stable and unstable (or asymptotic) manifolds of one more dimension than the object itself. The truncated Birkhoff normal form provides an accurate approximation of the local asymptotic manifolds in a neighborhood of $L_1$. In particular, the $4$-dimensional local stable and unstable manifolds of $\NHIMcn$ are given by
\[ W^{s,u}_\textrm{loc}(\NHIMcn) =\NHIMcn\oplus \Upsilon^{s,u}, \]
where $\Upsilon^{s,u}$ is a small interval in the stable (unstable) direction. In practice, we take this interval to be of length $\delta=10^{-3}$; see Section~\ref{subsec:NFprecision} for more details.
}

\MB{Once the local stable and unstable manifolds are obtained, the global manifolds are computed by integrating the local ones by the flow $\flow$ of \eqref{eq:Hamiltonian}.}

\subsection{Scattering Map}\label{sec:SM}
One of the main tools that we use in this paper is the scattering map.
This is an effective tool to quantify the effect of homoclinic excursions to a NHIM.
In general, it can be computed either perturbatively or numerically.

We recall the definition of the scattering map following \cite{DLS08}.
We consider a general setting of  flow $\Phi$ on a manifold $M$, and assume that there is a normally hyperbolic invariant manifold $\NHIM$ for the flow.
We will assume that the flow  as well as the geometric objects referred to below are differentiable enough, without formulating specific assumptions on regularity. (In the case of the spatial circular RTBP the flow is real analytic, and  the geometric objects of interest are at least $C^1$-differentiable, but not necessarily analytic.)

As a consequence of normal hyperbolicity, the stable manifold  $W^s(\NHIM)$ and the unstable manifold $W^u(\NHIM)$  are foliated by stable and unstable manifolds of points $W^s(y)$, $W^u(y)$, respectively, for $y\in\NHIM$. This implies that for each $x\in W^u(\NHIM)$ there exists a unique $x_-\in\NHIM$ such that
$x\in W^u(x_-)$, and for each   $x\in W^s(\NHIM)$ there exists
a unique $x_+\in\NHIM$ such that $x\in W^s(x_+)$.
These correspondences are described via the  \emph{wave maps}  $\Omega_+:W^s(\NHIM)\to \NHIM$ given by
$\Omega_+(x)=x_+$, and $\Omega_-:W^u(\NHIM)\to \NHIM$  given by $\Omega_-(x)=x_-$, respectively.

\MR{Assume that  $W^s(\NHIM)$ and  $W^u(\NHIM)$
  intersect transversally along a
homoclinic manifold $\Gamma$, that is,
 for each $x\in\Gamma  \subseteq W^u(\NHIM) \cap W^s(\NHIM)$, we have
\begin{equation*} \begin{split}\label{eq:dynamical channel}
T_xM=T_xW^u(\NHIM)+T_xW^s(\NHIM),\\
T_x\Gamma=T_xW^u(\NHIM)\cap T_xW^s(\NHIM).
\end{split} \end{equation*}
}

\MR{Further, assume that the intersection of the manifolds satisfy a strong transversality condition, that for each $x\in\Gamma$ we
have
\begin{equation*} \label{eq:transverse foliation}
\begin{split}
T_xW^s(\NHIM)=T_xW^s(x_+)\oplus T_x(\Gamma),\\
T_xW^u(\NHIM)=T_xW^u(x_-)\oplus T_x(\Gamma).
\end{split}
\end{equation*}
}

\MB{Assume the following transversality conditions:
\begin{itemize}
\item  $W^s(\NHIM)$ and  $W^u(\NHIM)$ intersect transversally along a  homoclinic manifold~$\Gamma$;
\item The homoclinic manifold $\Gamma$ is transverse to the stable (unstable) foliation $\{W^{s,u}_x\}_{x\in\NHIM}$ relative to $W^{s,u}(\NHIM)$.
\end{itemize}
}
Then, the restrictions $\Omega_+^\Gamma$, $\Omega_-^\Gamma$ of
$\Omega_+$, $\Omega_-$, respectively, to $\Gamma$ are  local
diffeomorphisms.
We can always choose $\Gamma$  so  that
$\Omega_+^\Gamma$, $\Omega_-^\Gamma$ are  diffeomorphisms onto their images.
A homoclinic manifold $\Gamma$ for which the corresponding
restrictions of the wave maps to $\Gamma$ are diffeomorphisms is
 referred to as a \emph{homoclinic channel}.

\begin{defn} \label{defn:scattering_map_flow}
Given a homoclinic channel $\Gamma$, the scattering map associated to
$\Gamma$ is the  diffeomorphism
\[\SM^\Gamma=\Omega^\Gamma_+\circ (\Omega^\Gamma_-)^{-1}\] from
the   $ \textrm{Dom}(\SM^\Gamma):=\Omega^\Gamma_-(\Gamma)\subseteq\NHIM$ to the
$ \textrm{Im}(\SM^\Gamma):= \Omega^\Gamma_+(\Gamma)\subseteq\NHIM$.
\end{defn}

That is, if $x\in \Gamma$ is a homoclinic point and $\Omega^\Gamma_\pm(x)=x_\pm$, then $\SM^\Gamma(x_-)=x_+$. In general, the domain and range of the scattering map are proper subsets of $\NHIM$. There are examples where the local domain
of the scattering map cannot be extended to a  global  one, that is, on the whole $\NHIM$, as moving along a loop in $\NHIM$ leads to
lack of monodromy (see \cite{DelshamsLS00}).

The scattering map depends on the choice of the homoclinic channel $\Gamma$. When we flow the homoclinic channel $\Gamma$ to $\Phi^t(\Gamma)$, the corresponding scattering maps are conjugate by the flow (see~\cite[Section 2.3]{DLS08}):
\begin{equation}\label{eq:conjugacy}
\SM^{\Phi^t(\Gamma)}=\Phi^{t}\circ\SM^\Gamma\circ \Phi^{-t}.
\end{equation}
Of course, when $\Gamma$ and $\Phi^t(\Gamma)$ overlap, for $x\in  \Gamma\cap\Phi^t(\Gamma)$, we have $\SM^{\Phi^t(\Gamma)}(x_-)=\SM^\Gamma(x_-)$.
This means that $\SM^\Gamma$ can be continued to $\SM^{\Phi^t(\Gamma)}$ for some interval of times  $t$, for as long as the corresponding homoclinic channels are well defined and can be continued into one another. We will regard the result of such continuation to the maximal domain as a single scattering map.

When the choice of the homoclinic channel $\Gamma$ is evident from the context, we drop the superscript  from the notation $\SM^\Gamma$.

In the case of a  discrete-time dynamical system, the scattering map can be defined in a similar fashion.

A remarkable property of the  scattering map is that it is exact symplectic, provided that the manifold and the flow are exact symplectic.
We refer to \cite{DLS08} for details.

\begin{remark} As we shall see in Section~\ref{Sec:DomainGlobalSM}, in our model we can construct two scattering maps that are defined on a whole annulus inside the NHIM. In this sense, the scattering maps are globally defined on the annulus.
Each of this scattering maps is obtained by a continuation of a locally defined scattering map to a maximal domain.
In our model, the scattering map extended to its maximal domain satisfies the monodromy condition, as moving around on a non-trivial loop inside the annulus does not yield a multi-valued map.\end{remark}

\subsection{Transition map}
\label{sec:transition}
Assume the general setting from Section \ref{sec:SM}.

From the definition of the scattering map associated to $\Phi^t$, it follows that,
if $x\in \Gamma$ is a homoclinic point and
\[\SM(x_-)=x_+,\textrm{ where } x_\pm\in\NHIM,  \]
then
\begin{equation*}\label{eqn:asymptotic_SM}
\begin{split}
  d(\Phi ^{-t_-}(x),\Phi^{-t_-}(x_-))& \to 0, \textrm { as } t_-\to\infty,\\
  d(\Phi ^{t_+}(x),\Phi^{t_+}(x_+))& \to 0, \textrm { as } t_+\to\infty.
\end{split}
\end{equation*}

This correspondence represents the so-called \emph{transition map} (see \cite{DelshamsGR12})).
\MB{
\begin{defn} \label{defn:transition}
Given $\delta>0$, let $T_-=T_-(\delta)$, $T_+=T_+(\delta)$ be the infimum of the positive reals $t_-$, $t_+$, respectively, for which
\begin{equation}
\label{eqn:asymptotic_SM delta}
d(\Phi^{-t_-}(x),\Phi^{-t_-}(x_-)) \le \delta,\quad  d(\Phi^{t_+}(x),\Phi^{t_+}(x_+))\le \delta.
\end{equation}
Then the \emph{transition map}
$\mathcal{T}_\delta: \Phi ^{-T_-}(\textrm{Dom} (\SM ))\to \Phi^{T_+}(\PB{\textrm{Im}} (\SM))$
is defined as
\[\mathcal{T}_\delta=\Phi^{T_+}\circ \SM \circ\Phi ^{T_-} .\]
The transition map depends on the choice of the homoclinic channel $\Gamma$ and on the distance $\delta$ to the NHIM.
\end{defn}
}

Thus, for $T_-,T_+$ large enough, the homoclinic orbit segment $\Phi^{[-T_-,T_+]}(x)=\{\Phi^{t}(x)\,|\,t\in[T_-,T_+]\}$ in $M$ is  an approximation of the pseudo-orbit $\Phi^{T_+}\circ \SM\circ\Phi^{T_-}(x_-)$. The former is an orbit segment in the manifold $M$, while the latter is given by an orbit segment of the inner dynamics ${\Phi^{T_-}}_{\mid \NHIM}$, followed by one application of the scattering map $\SM$, followed by another orbit segment of the inner dynamics ${\Phi^{T_+}}_{\mid\NHIM}$.

\MB{
\begin{remark} \label{rem:scattering_vs_transition} 
Note that, for the computation of both the scattering map and the transition map, one needs to know the inner dynamics, given by the restriction to the NHIM.
Compared to the scattering map, the transition map depends on an extra parameter (the threshold distance $\delta$ to the NHIM), and so,
different choices of this parameter yield different transition maps.
\end{remark}
}

\MB{
\begin{remark} \label{rem:flow_vs_flown} In the setting of Section~\ref{Sec:RTBP}, we can consider the dynamics of the original flow $\flow$ of the RTBP, and the dynamics of the flow $\flown$ associated to the Birkhoff normal form, as well as the corresponding NHIMs $\NHIMc$ and $\NHIMcn$, respectively.
In each case, we can define a scattering map and a transition  map. We denote by $\SM$ the scattering map corresponding to $\flow$, and by $\SMn$ the scattering map corresponding to $\flown$. The map $\SM$ is defined on some domain in $\NHIMc$, and $\SMn$ is defined on some domain in $\NHIMcn$. 
Similarly, we denote by $\mathcal{T}_\delta$ the transition map associated to $\flow$, and by $\mathcal{T}_\delta^N$ the transition map associated to $\flown$. 
\end{remark}
}

\subsection{Reduction of the Scattering Map and the Transition Map to a Poincar\'e Section}
\label{sec:reduction}

We now consider the setting from  Section~\ref{Sec:RTBP}.

In our previous paper~\cite{DGR16}, we showed that, in the case of the  spatial circular RTBP, the unstable and stable manifolds of $\NHIMcn$ intersect, giving rise to homoclinic orbits to $\NHIMcn$. We can select a homoclinic channel $\Gamma_c$ and consider  
the scattering map associated to the truncated Birkhoff normal form:
\[ \SMn:\textrm{Dom}(\SMn)\subseteq\NHIMcn\to \textrm{Im}(\SMn)\subseteq\NHIMcn. \]
\MR{We note in this case the NHIM $\NHIM_c^{(N)}$ is $3$-dimensional, hence not symplectic,
so the corresponding scattering map is not symplectic.
As described in Section \ref{sec:SM}, the scattering map assigns to a point $x_-\in\NHIMcn$ another  point $x_+\in\NHIMc$ whenever there is a homoclinic point $x\in \Gamma_c$ such that the orbit of $x$ tends to the orbit of $x_-$ in the past and to the orbit of $x_+$ in the future.}

\MB{To reduce the dimensionality of the scattering map, we consider the Poincar\'e section
$\Sigma = \{ \phi_p=0 \}$, with associated first return map $\FRM \colon \Sigma\to\Sigma$.
Let $\NHIMsecn$ denote the intersection of the NHIM with the Poincar\'e section:
\[ \NHIMsecn= \NHIMcn\cap\Sigma.\]
We have shown in ~\cite[Section 3.2]{DGR16} that $\NHIMsecn$ is a normally hyperbolic invariant manifold for $\FRM$, which we call the \textit{reduced NHIM}.
This is a $2$-dimensional manifold, that can be parameterized in terms of $(J_v,\phi_v)$.
It has a symplectic structure given by $dJ_v\wedge d\phi_v$.
The scattering map induces a \textit{reduced scattering map}:
\[\smn\colon \textrm{Dom}(\smn)\subseteq\NHIMsecn\to \textrm{Im}(\smn)\subseteq \NHIMsecn.\]
The reduced scattering map $\smn$ is exact symplectic. (See~\cite{DLS08}).}

\MB{Inside $\NHIMsecn$ we consider an annulus given by
\[ \domain = \NHIMsecn \cap\{\ J_v\in[J_1,J_2],\ \phi_v\in[0, 2\pi)\} \]
where $J_1<J_2$ will be specified in Section~\ref{Sec:DomainGlobalSM}.}

For each fixed
value of $J_v = \bar{J}_v$ in the range, there is a unique invariant curve for $\FRM$
\MB{
\[
\Tsecn_{\bar{J}} = \domain\cap\left\{J_v = \bar{J}_v\right\}.
\]
}

One can derive the reduced scattering map $\sm$ from the full scattering map $\SMn$ as follows.
Let $x_\pm$ be two points related by the scattering map: $x_+ = \SMn(x_-)$.
Then we flow $x_-$ backwards to the Poincar\'e section $\Sigma$, obtaining a new point $(J_-,\phi_-)\in\NHIMsecn$. Similarly, we flow $x_+$ forward to $\Sigma$, obtaining  $(J_+,\phi_+)\in\NHIMsecn$.
The reduced scattering map $\smn$ takes $(J_-,\phi_-)$ to $(J_+,\phi_+)$.

When it is clear from the context, we will  abbreviate `reduced scattering map' to just `scattering map'. Also, for ease of notation, we will simply write $(J,\phi)$ for $(J_v, \phi_v)$.

\MR{\begin{remark}
	A scattering map is not unique: it depends on the chosen homoclinic channel. In this paper, we will show that there exist two different channels, and thus two scattering maps defined on a common domain $\domain$ inside the NHIM $\NHIMsecn$. The domain will be made explicit in Equation~\eqref{eq:domainJ}.
\end{remark}}

The transition map can be reduced to a surface of section in a similar fashion. Let  $x\in \Gamma^\Sigma_c=\Gamma_c\cap \Sigma$ be  a homoclinic point and
\[\smn(x_-)= x_+,\textrm{ where }  x_\pm\in\Lambda^\Sigma_c. \]
In terms of  the $(J,\phi)$ coordinates we can write $ x_-=(J_-,\phi_-)$ and $ x_+=(J_+,\phi_+)$. 
For a given $\delta>0$ \AB{ where the Birkhoff normal form applies (see Section~\ref{subsec:NFprecision})}, let $K_-,K_+$ be the smallest positive
integers $k_-,k_+$, respectively, such that
\begin{equation}
\label{eqn:asymptotic_sm delta}
d(\FRM^{-k_-}( x),\FRM^{-k_-}( x_-))<\delta,\quad  d(\FRM^{k_+}(  x),\FRM^{k_+}(\ x_+))<\delta.
\end{equation}

Then the transition map is given by 
\[\tau^N_{\delta}  =\FRM^{K_+}\circ \smn\circ\FRM^{K_-}. \]
 
\begin{remark} Since the dynamics along the hyperbolic manifolds $W^{u,s}(\Lambda^\Sigma_c)$ is much faster then the inner dynamics $\FRM_{\mid\Lambda^\Sigma_c}$, it may happen \MR{(as this is the case for the $\delta>0$ where the Birkhoff normal form applies,  see Section~\ref{subsec:NFprecision}),}   \AR{is possible that for some suitable $\delta$,}   that condition  \eqref{eqn:asymptotic_sm delta} is satisfied for  $K_-=K_+=1$. 
In that case, the transition map is given by
\[\tau_\delta^N=\FRM\circ\smn\circ \FRM.\]
Indeed, the transition map used in this paper will be of this type.

\end{remark}

\MB{
\section{Main results}}

In this section we first provide an abstract result for an iterated function system (IFS) consisting of an inner map (given by the restriction of the dynamics to a NHIM) and finitely many scattering maps, saying that, if the IFS satisfies some verifiable, quantitative conditions, then  the system has true `diffusing' orbits. 
Then, we give the main numerical result of the paper, which amounts to the numerical verification of the assumptions of the abstract result.  


\begin{theorem}\label{th:main}
Assume that $\FRM: \Sigma\to\Sigma$ is a diffeomorphism on a manifold $\Sigma$,
$\Lambda$ is a $2$-dimensional NHIM for $\FRM$, $(I,\phi)$ a coordinate system on $\Lambda$, $\omega=dI\wedge d\phi$   a symplectic form on $\Lambda$, and
\[\domain= [I_1, I_2]\times \mathbb{T}^1\subset \Lambda.\]
 
Assume that $W^s(\domain)$ and $W^u(\domain)$ intersect along a homoclinic channel $\Gamma_i$, $i\in\{1,\ldots,L\}$, where $L\ge 1$.
Let $f=\FRM_{\mid\domain}$ be the restriction of $\FRM$ to $\Lambda$ and $\sigma_i$ be the scattering map associated to $\Gamma_i$, $i\in\{1,\ldots,L\}$.

\begin{itemize}
\item[I.] We assume that the \emph {inner map}  $f$    satisfies the following conditions:
\begin{itemize}
\item[(I.i)]  The map $f$ is an exact symplectic\footnote{A map $f$ on $\domain$ is exact symplectic if $f^*(I d\phi)-I d\phi=ds$ for some function $s$ on $\domain$.} twist map, i.e.
\begin{equation}\label{eq:tw}\frac{\partial \pi_\phi f}{\partial I}(I,\phi)>0.\end{equation}
 \item[(I.ii)] There exists $\rho_1>0$ such that for every   $\bar{I}\in [I_1, I_2]$, the level set   $\{I=\bar{I}\}$ is invariant under  $f$  up to an error of $\rho_1$,
 that is, for all $n\ge 0$ and $\bar{\phi}\in[0,2\pi)$, we have
\[ |\pi_I(f^n(\bar{I},\bar{\phi}))- \bar{I}|<\rho_1.\]
\end{itemize}
\item[II.] We assume that each \emph {scattering map}   $\sigma_i$, $i\in\{1,\ldots,L\}$, satisfies the following conditions:
\begin{itemize}
\item [(II.i)] $\sigma_i$ is globally defined on $\domain$, i.e., $\textrm{dom}(\sigma)\supset \domain$;
\item [(II.ii)] There exist $C^1$-functions $\gft_i=\gft_i(I,\phi')$ and $\omega_i=\omega_i(I)$ such that
\begin{equation}\begin{split}\label{eqn:scattering_repr}
I'=&I+\frac{\partial \gft_i}{\partial \phi'}(I,\phi'),\\
\phi=&\phi'+\omega_i(I)+\frac{\partial \gft_i}{\partial \phi'}(I,\phi'),
\end{split}
\end{equation}
where we denote $\sigma_i(I,\phi)=(I',\phi')$;
\item [(II.iii)] There exists $\rho_2>\rho_1$ such that for every     $\bar{I}\in [I_1,I_2]$ there exists $i\in\{1,\ldots,L\}$ and $\bar{\phi}\in[0,2\pi)$ such that
\begin{equation}\label{eqn:der_gen} \frac{\partial \gft_i}{\partial \bar{\phi}'}(\bar{I},\bar{\phi}')>\rho_2.\end{equation}
\end{itemize}
\end{itemize}
Then, for any neighborhood $U_1$ of $\{I_1\}\times\mathbb{T}^1$ and $U_2$ of $\{I_2\}\times\mathbb{T}^1$ in $\Sigma$  there exists an orbit $\{z_n\}_{n=0,\ldots,N}$ of $\FRM$ in $\Sigma$ such that
\[ z_0 \in U_1 \textrm{ and } z_N\in U_2.\]
\end{theorem}

\begin{proof}

The proof follows immediately from the following result in \cite{Moeckel02}, and its extensions in \cite{LeCalvez07} and \cite{gidea_marco_2017}.

\begin{theorem}[\cite{Moeckel02}]\label{thm:Moeckel}
Let $\domain=[0,1]\times \mathbb{T}^1$ and $f,g:\domain\to\domain$ be two $C^1$-diffeomorphisms of the annulus that preserve the boundary circles.

Assume the following
\begin{itemize}
\item [(i)] $f$ is a  twist map;
\item [(ii)] $g$ is an area preserving map;
\item [(iii)] Every essential\footnote{ An essential  circle  is a  simple closed $C^0$-curve in $\domain$ that is nonhomotopic to zero.} $f$-invariant circle $\Gamma$  is not invariant under $g$.
\end{itemize}

Then, for every pair of open neighborhoods $U_0$ of $ \{0\}\times\mathbb{T}$ and $U_1$ of $\{1\}\times\mathbb{T}$ in $\domain$, there exists an orbit  $(x_n)_{n=0,\ldots,N}$ of the IFS generated by $\{f,g\}$ such that $x_0\in U_0$ and $x_N\in U_1$ 
(the orbit of the IFS is defined at every step by a choice $x_{n+1}=f(x_n)$ or $x_{n+1}=g(x_n)$).
\end{theorem}

In Theorem \ref{thm:Moeckel} it is essential that the map $g$ is globally defined on the annulus $\domain$. 
On the other hand, the condition that $f$ preserves the boundary components of the annulus can be replaced by the condition  that $f$ is exact symplectic.
Also, the conditions  that $g$ preserves the boundary components of the annulus and that it does not preserve any  essential $f$-invariant circle can be replaced by the  condition that for any essential $f$-invariant circle, there are points below the circle that are mapped by $g$ to points above the circle. 

The theorem extends immediately when instead of a single map $g$ we have finitely many maps $g_1,\ldots,g_L$ satisfying (ii) such that every essential $f$-invariant circle $\Gamma$ is not invariant under some $g_i$, $i\in\{1,\ldots,L\}$. 

\MR{In the setting of Section~\ref{Sec:RTBP}, $f$ is the inner map defined as $\FRM_{\mid\domain}$, and we have two maps $g_1, g_2$ represented by the scattering maps $\sigma_1$ and $\sigma_2$. }

Condition (I.i) of Theorem \ref{th:main} in $f$ represents condition (i) of Theorem \ref{thm:Moeckel}. Condition (II.i) on $\sigma_i$ from Theorem \ref{th:main}  ensure that  these maps are globally defined on $\domain$.

Condition (I.ii) of Theorem \ref{th:main} implies that each essential invariant circle for $f$ is contained within  $\rho_1$ of some $I$-level set. 
Conditions (II.ii) and (II.iii) of Theorem \ref{th:main}  imply that for each $I$-level set, there is a scattering map $\sigma_i$ and  a point on the $I$-level set whose  $I$-coordinate is increased by  $\sigma_i$ by more than $\rho_2>\rho_1$.  Therefore, any essential 
$f$-invariant circle fails to remain invariant under at least one of the
$\sigma_i$.
Theorem \ref{thm:Moeckel} implies that there is a pseudo-orbit $(x_n)_{n=0,\ldots,N}$ of the IFS generated by $\{f,\sigma_1,\ldots,\sigma_L\}$ such that $x_0\in U_0$ and $x_N\in U_1$. 
To obtain a true orbit, we use the following result that follows directly from \cite[Lemma 3.2 and Theorem 3.7]{GideaLlaveSeara20-CPAM}.
\begin{theorem}\label{thm:GLS2020}
Assume that $(\Sigma, \omega)$ is a symplectic manifold of finite symplectic volume, 
$\FRM:\Sigma\to\Sigma$ is a symplectic diffeomorphism, $\Lambda\subseteq M$ is a NHIM for $\FRM$, $f=\FRM_{\mid \Lambda}$ and $\sigma_i$, $i=1,\ldots,L$, is a family of scattering maps on $\Lambda$.

Then every finite pseudo-orbit   $\{x_i\}_{i=0,\ldots,N}$  of the IFS generated by $\{f,\sigma_1,\allowbreak\ldots,\allowbreak\sigma_L\}$ in $\Lambda$ can be shadowed by a true orbit, 
that is,  for every $\delta>0$, there exists     an orbit $\{z_i\}_{i=0,\ldots,N}$ of $\FRM$ in $\Sigma$, with $z_{i+1}=f^{k_i}(z_i)$ for some $k_i>0$,  such that
$d(z_i,x_i)<\delta$ for all $i=0,\ldots,N$.
\end{theorem}

Applying Theorem \ref{thm:GLS2020} to the $(x_n)_{n=0,\ldots,N}$ of the IFS generated by $\{f,\sigma_1,\ldots,\sigma_L\}$ yields a true orbit of $\FRM$ from $U_0$ to $U_1$. 
\end{proof}

We now provide an outline of the numerical verification of the conditions of Theorem \ref{th:main}. Details are provided  
in Section~\ref{sec:HeuristicArnoldDiffusion}.

\begin{result}\label{thm:main_numerical}
Consider the setting from Section~\ref{Sec:RTBP}.
\begin{itemize}
\item[(a)] In Section~\ref{Sec:NumericalResults} we obtain the approximation of the NHIM $\NHIMc$ by $\NHIMcn$, using the Birkhoff normal form $H^{(N)}_{\textrm{trunc}}$. 
\item[(b)]  In Section \ref{subsec:NFprecision} we show that the error in the numerical approximation of the NHIM is less than $10^{-15}$ (close to machine precision).  
\item[(c)]  In Section~\ref{subsubsec:GlobalizingSM} we construct two scattering maps that are globally defined on an annular $\domain$ inside $\NHIMcn$, corresponding to the range $J\in[J_1,J_2]$, where $J_1=0.001$ and $J_2=0.007$.
See Section~\ref{Sec:DomainGlobalSM} for additional details.
\item[(d)]  In Section~\ref{Sec:SeriesRepresentation} we rescale $I=10^3\cdot J$ and derive series representations of the two scattering maps, which are of the form \eqref{eqn:scattering_repr}.
\item[(e)]  In Section~\ref{sec:ApproximationError} we show that the generating functions $\gft_i$ associated to these scattering maps satisfy condition \eqref{eqn:der_gen} for $\rho_2=0.1$. By (b), the level sets of $I$ are preserved by the inner map up to an error of $\rho_1=10^{-12}$. 
\item[(f)]  In Section~\ref{sec:inner_map} we verify that the inner map satisfies a twist condition.
\end{itemize}
Provided that the conditions (a)--(f) have been verified, then there exists a true, zero-cost orbit of $\flow$ along which $I$ changes from $I_1=1$ to $I_2=7$.
\end{result}

\begin{remark}
A result similar to Theorem \ref{th:main} can be formulated in terms of  scattering maps only (without using the inner map).
More precisely, assume that we have a family of scattering maps $\{\sigma_1,\ldots,\sigma_L\}$ as in
Theorem~\ref{th:main}, and that one of them, say $\sigma_1$, satisfies a twist condition
as in \eqref{eq:tw}. Also assume that there exist $0<\rho_1<\rho_2$ such that
\[ \max_{\bar\phi'}\left\vert\frac{\partial \gft_1}{\partial \bar{\phi}'}(\bar{I},\bar{\phi}')\right\vert<\rho_1,\]
and for each $\bar I$ there exists $i\neq 1$ such that
\[ \max_{\bar\phi'}\frac{\partial \gft_i}{\partial \bar{\phi}'}(\bar{I},\bar{\phi}')>\rho_2.\]
Then there exist an orbit $\{z_n\}_{n=0,\ldots,N}$ of $\FRM$ as in the conclusion of 
Theorem~\ref{th:main}. 

The main idea is that $\sigma_1$ plays the  role of the inner dynamics in Theorem~\ref{th:main}. 
\MR{Since we are not explicitly using this mechanism in this paper, we omit the details. }
\end{remark}

\MB{In  Sections 4--12 below, we compute numerically the flow $\flown$,  the corresponding NHIM  $\NHIMcn$, two scattering maps and two transition maps, as well as their reductions to the surface of section $\Sigma$. 
In Section \ref{subsec:NFprecision} we  give an estimate of the error between the original flow $\flow$ and the flow computed via the  Birkhoff normal form $\flown$, obtaining that $\NHIMcn$ is an accurate enough numerical approximation of $\NHIMc$.
As all subsequent computations are performed via  the Birkhoff normal form, to simplify the notation, from now on we drop the superscript/subscript $N$ from the notation for all objects. 
} 

\section{Numerical Scattering Map on a Grid for \texorpdfstring{$C=3.00088$}{C=3.00088}} \label{Sec:NumericalResults}

Using the methodology presented in our previous paper~\cite{DGR16}, we compute the NHIM $\NHIM_c$ and its stable/unstable manifolds. We show that the asymptotic manifolds intersect transversally along two homoclinic channels, giving rise to two different scattering maps. We compute the scattering maps numerically at a grid of points; they are shown in Figure~\ref{fig:curves}.

For the purpose of this paper we will use the energy value $C=c \coloneqq 3.00088$. This value is chosen after the appearance of the equilibrium point $L_1$ ($c_1 \coloneqq 3.00090$), but before the appearance of \textit{halo} orbits ($c_\text{halo} \coloneqq 3.00082$). The choice of energy is motivated by two reasons.

Firstly, $c$ is close enough to $c_1$ so that the dynamics around $L_1$ is almost integrable. Thus, the (integrable) Birkhoff normal form $H^{(N)}_{\textrm{trunc}}$ provides a good approximation to the local dynamics. (Section~\ref{subsec:NFprecision} quantifies the Birkhoff normal form error).

Secondly, $c$ is close enough to $c_1$ so that we are in the setting of \textit{a priori chaotic Arnold diffusion} (see Section \ref{Sec:Introduction}). 

\MR{Assume that  for the exact RTBP (not the truncated Birkhoff normal form) there exists a NHIM that is close to $\NHIM_c$, and is almost filled with invariant tori. (The existence of the NHIM  was rigorously proven for the planar RTBP in~\cite{CapinskiRoldan12}). Then the tori do not separate the phase space, so there may exist trajectories that escape a given neighborhood of $L_1$. Proving the existence of such orbits is a classical problem that has attracted lots of attention.
In this paper, we provide new tools to study the problem of Arnold diffusion. We hope that this will lead to the explicit construction of diffusion orbits, as well as explicit bounds on their diffusion time.}

Fixing the energy value $c = 3.00088$, the NHIM $\NHIM_c$ consists of a continuous family of 2-dimensional invariant tori around $L_1$, which we parametrize by the vertical action $J$. The vertical action increases along the family from $J=0$ to $J=J_\text{max}\coloneqq 0.052$.
Correspondingly, the planar action $J_p$ decreases from $0.05029$ to $0$.

\begin{figure}\footnotesize

	\scalebox{0.9}{\input{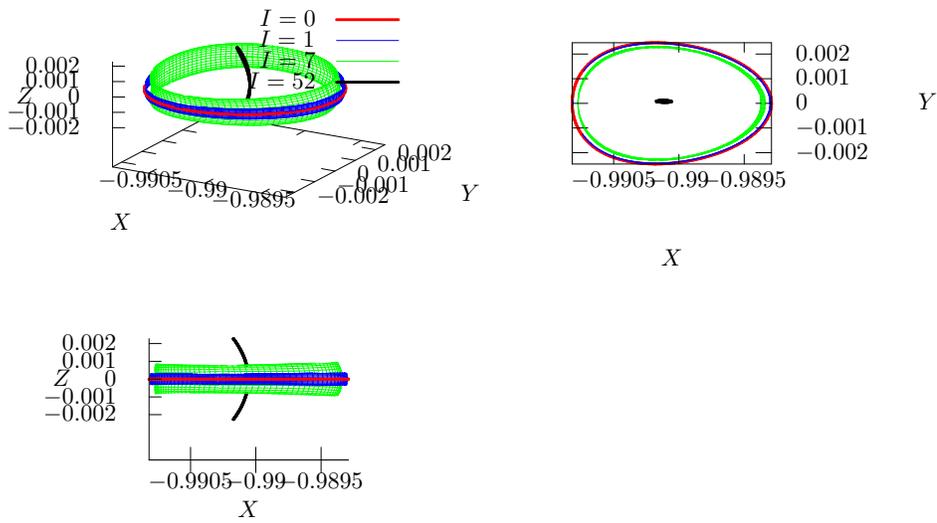}}
	\caption{The NHIM $\NHIM_c$ consists of a continuous family of invariant tori around $L_1$. The endpoints of the family are the planar and vertical Lyapunov orbits (shown in red and black, respectively). In between, there are 2-d tori of increasing vertical amplitude $J$. (For clarity, only two of them are shown). The transfer trajectory shown in Figure~\ref{fig:initial_final_segments} starts at an initial condition very close to the blue torus $I=1$ and ends very close to the green torus $I=7$, where $I=10^3\cdot J$} 
	\label{fig:centermfld}
\end{figure}

Figure~\ref{fig:centermfld} shows some tori in the NHIM $\NHIM_c$.
\begin{itemize}
	\item $J=0$ (i.e. $J_p=0.05029$) corresponds to the unique planar Lyapunov orbit in this energy level.
	
	\item $J=0.052$ (i.e. $J_p=0$) corresponds to the unique vertical Lyapunov orbit in this energy level.
	
	\item Every intermediate action $J\in(0, 0.052)$ corresponds to a 2-dimensional torus located between the planar and the vertical Lyapunov orbits.
\end{itemize}
As seen in Figure~\ref{fig:centermfld}, the NHIM $\NHIM_c$ spans a spherical region of radius $0.002$ Astronomical Units (AU), or roughly 300000 km around the equilibrium point~$L_1$.

\AB{\subsection{Accuracy of the Birkhoff normal form} \label{subsec:NFprecision}}
	The accuracy of the Birkhoff normal form (BNF) expansion has been tested against numerical integration of the RTBP equations, following the same procedure as in~\cite{JorbaMasdemont99}. Compute an initial condition on $\NHIM_c$ by evaluating the BNF up to order $N=16$. This initial condition is integrated for $\pi$ units of adimensional time using two different methods:
	\begin{enumerate}
		\item Using the BNF. No numerical integration is needed. As the Hamiltonian is integrable, and we have it integrated, we simply tabulate the solution at time $\pi$.
		\item Integrate the RTBP equations using a Runge-Kutta-Feldberg numerical integrator of order 7-8, with local error at each step within $10^{-14}$.
	\end{enumerate}
	Then compare the two final conditions.
	
	This test has been performed for several initial conditions on $\NHIM_c$. In all cases, the difference in the Euclidean norm for the final condition is less than $10^{-12}$ adimensional RTBP units. Thus the initial condition was very accurate, in the sense that it is very close to one of the tori computed by the BNF.
	
	In fact, it is known that, due to the hyperbolicity of orbits around the collinear point, errors increase by a factor close to 1500 after $\pi$ units of time. Therefore the error in the initial condition is less than $10^{-12}/1500$ adimensional units, close to machine precision.

\AB{Furthermore, this validity and accuracy of the BNF of order $N$ in a neighborhood of size $\delta$ along the unstable/stable coordinates, say for $\delta=10^{-3}$ and $N=16$, can also be analytically corroborated using theoretical results, for example from~\cite{DG96}. The bound of the error in a neighborhood of the libration point is there given as a function of the bounds of the Taylor expansion of the Hamiltonian, as well as of the small divisors $k_\mathrm{p}\nu_\mathrm{p}+k_\mathrm{v}\nu_\mathrm{v}$ for integers satisfying $\abs{k_\mathrm{p}}+\abs{k_\mathrm{v}}\leq N$.}

\AB{It should be noted that the presence of these small divisors is what affects the convergence of these Birkhoff normal forms. Although they depend on the mass ratio $\mu$, they do so very slightly for libration points; see, for instance,~\cite{Jorba99,JorbaMasdemont99}. Therefore, their influence in detecting strong resonances \MR{, such as halo orbits,} is not noticeable up to normal form orders much higher than $N=16$.}

\subsection{Homoclinic Orbits}
\label{subsec:HomoclinicOrbits}
In~\cite{DGR16}, Section~4.3, we explained in detail how to compute the intersection of the stable and unstable manifolds $W^s(\NHIM_c) \cap W^u(\NHIM_c)$ restricted to the surface of section \[\mathscr{S}=\{(X,Y,Z,\dot X,\dot Y,\dot Z)\,:\, Y=0,\,\dot Y>0\}.\]
We find that the asymptotic manifolds do indeed intersect transversally, giving rise to families of homoclinic orbits to the $\NHIM_c$. These homoclinics will later be encoded by two scattering maps.

Briefly, the numerical procedure to compute each homoclinic consists of finding two initial conditions $y_-$, $y_+$ on the \textit{local} unstable resp. stable manifolds, and a point $x\in\mathscr{S}$ such that:
(1) $\flow^t(y_-) = x$; and
(2) $\flow^{-s}(y_+) = x$.
Then $x$ is a homoclinic point, generating a homoclinic orbit segment from $y_-$ to $y_+$.

\begin{figure}
	\scalebox{0.9}{\input{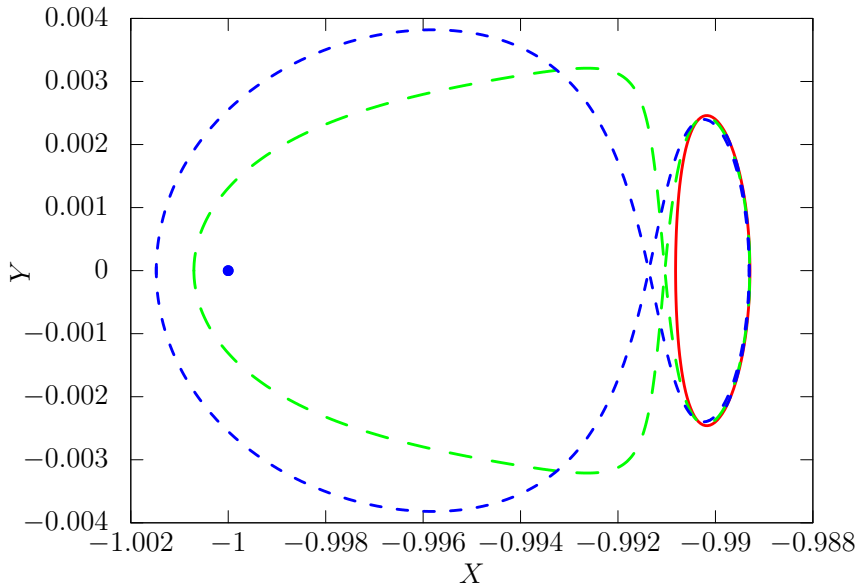}}
	\caption{Planar Lyapunov orbit (solid line), and its two `primary' homoclinics (dashed). Both homoclinics travel around the Earth once (located at $X=-1+\mu$).}
	\label{fig:homocl_orbits_lyap}
\end{figure}

We would like to emphasize some aspects of this computation:

\begin{itemize}
        \item We consider \textbf{only the first intersection} of the stable and unstable manifolds with the surface of section  $\mathscr{S}$. The `primary' homoclinic connections generated in this way travel around the Earth once. 
	\item The initial conditions $y_-, y_+$ are taken at a distance $\delta=10^{-3}$ of the NHIM along the unstable/stable coordinates. This guarantees that they are inside the domain where the Birkhoff normal form is accurate.
	
	\item At the same time, these initial conditions are sufficiently far from the NHIM so that the homoclinic orbit segment does not wind around $L_1$ more than once.
	
	\item All homoclinics take $5.936738 \leq t+s \leq 6.000688$ time units to go from $y_-$ to $y_+$. We will refer to $t+s$ as the `flight time' of the homoclinic. The shortest and longest flight times correspond to the two homoclinics of the Lyapunov orbit. See Figure~\ref{fig:homocl_orbits_lyap}.
\end{itemize}

\MR{\begin{remark}
The section $\mathscr{S}$ corresponds to crossing the $XZ$-plane (with positive $Y$ velocity) in co-rotating SCRTBP coordinates. The choice of the section is somewhat arbitrary; it is only related to the numerical computation. We consider \textbf{only the first intersection} of the stable and unstable manifolds with the surface of section  $\mathscr{S}$. The `primary' homoclinic connections generated in this way travel around the Earth once.
	There exist subsequent (second, third, etc.) intersections with $\mathscr{S}$, but they are not as interesting for applications, since these homoclinic connections are longer (they travel around the Earth twice, thrice, etc.).
	For more details, see~\cite[Section 4.3]{DGR16}.
\end{remark}}

Following the decomposition of $\NHIM_c$ into invariant tori, we first study homoclinic orbits from each $T_J$ to itself. Later we will study homoclinic orbits from $T_J$ to nearby tori.

\begin{itemize}
	\item When $J=0$, the invariant manifolds $W^u(T_0)$ and $W^s(T_0)$ have $2$ transverse intersections in the section $\mathscr{S}$. Hence, the planar Lyapunov orbit $T_0$ has $2$ homoclinic connections with itself. See Figure~\ref{fig:homocl_orbits_lyap}.  
	
	\item For every fixed $J\in (0,0.01)$, the invariant manifolds $W^u(T_J)$ and $W^s(T_J)$ have $8$ transverse intersections. That is, every  torus $T_J$ in this domain has $8$ homoclinic connections with itself. 
	
	\item In contrast, for action values $J$ above $0.01$, the invariant manifolds $W^u(T_J)$ and $W^s(T_J)$ cease to intersect (at their first intersection with the section $\mathscr{S}$).
\end{itemize}

\begin{remark}\label{rem:FlightTimes}
	The planar Lyapunov orbit has two `primary' homoclinics that travel around the Earth (see Figure~\ref{fig:homocl_orbits_lyap}). However, one of them makes a longer excursion than the other: the green homoclinic has flight time $5.936738$, while the blue one has flight time $6.000688$.
\end{remark}

\begin{remark}
	Each of the $2$ intersections in $W^u(T_0) \cap W^s(T_0)$ gives rise to $4$ intersections when we increase the dimension of the manifolds $W^u(T_J), W^s(T_J)$ by one. This is expected by Morse theory~\cite{Milnor63}.
\end{remark}

\PB{\subsection{Computation of Local Transition Map and Scattering Map.}
}

\PB{
Consider a homoclinic point $x\in\mathscr{S}$ as found in Section~\ref{subsec:HomoclinicOrbits}, generating a homoclinic orbit segment from $y_-\in W^u_{loc}(\NHIMc)$ to $y_+\in W^s_{loc}(\NHIMc)$, with flight time $t+s$.
Let 
\[ y_-=(y_h^-=\delta, J_p^-, J_v^-, x_h^-=0, \phi_p^-, \phi_v^-)\]
and 
\[ y_+=(y_h^+=0, J_p^+, J_v^+, x_h^+=\delta, \phi_p^+, \phi_v^+)\]
in Birkhoff normal form coordinates. 
Then, the transition map $\mathcal{T}_\delta \colon\ \NHIMc  \to \NHIMc $ is given by 
\[(J_p^-, J_v^-, \phi_p^-, \phi_v^-) \mapsto (J_p^+, J_v^+, \phi_p^+, \phi_v^+).
\]
As explained in Section~\ref{sec:transition}, the scattering map $S\colon\ \NHIMc\to \NHIMc$ is related to the transition map by the flow. Therefore, the scattering map is given by
\[x_-\coloneq \Phi^{t}(J_p^-, J_v^-, \phi_p^-, \phi_v^-) \mapsto x_+ \coloneq \Phi_H^{-s}(J_p^+, J_v^+, \phi_p^+, \phi_v^+).
\]
}

Every transverse intersection
\[ x \in W^u(T_J) \cap W^s(T_J) \]
implies that there exists a pair of points $x_-, x_+ \in T_J$ such that $W^u(x_-)$ intersects $W^s(x_+)$ at the homoclinic point $x$. 
Since this is an open condition, we can define a local scattering map on some open set containing $x_-$ by
$x_- \mapsto \SM(x_-) \coloneqq x_+$.

In fact, we find that these local scattering maps can be continued to form two \textit{global} scattering maps, which we will denote $\SM_1$ and $\SM_2$. They determine the reduced scattering maps, $\sm_1$ and $\sm_2$ (see Section~\ref{sec:reduction}).

\subsection{Extending the Scattering Map}\label{subsubsec:GlobalizingSM}

Let us explain how $\SM_1$ is numerically constructed. ($\SM_2$ is constructed analogously). The general idea is to compute $\SM_1$ on several tori $T_J$ (for example $J=0.001$, $0.002$, $\dotsc$, $0.007$). \MR{Given that the NHIM $\NHIM_c$ is the union of all tori,} 
\MB{This provides a coarse representation of $\SM_1$.}
\MR{on a whole annulus $\domain$ inside $\NHIM_c$ (see Equation~\eqref{eq:domainJ}).}

\begin{figure}
	\scalebox{0.9}{\input{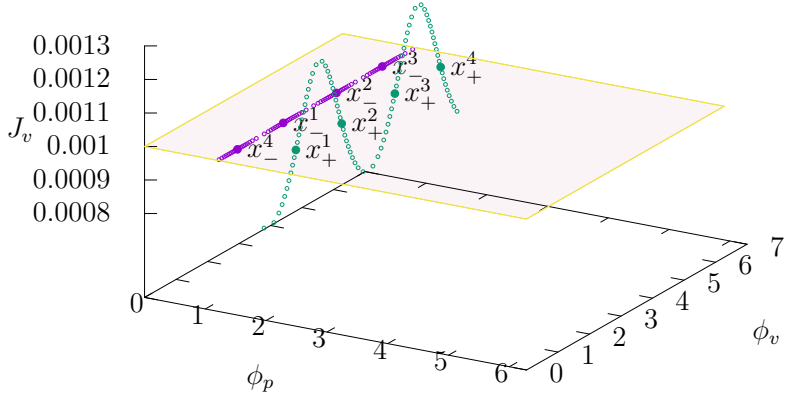}}
	\caption{Action of the scattering map $\SM_1$ on the torus $T_J$ with $J=0.001$. $\SM_1$ maps the purple set to the green set. When the `source' and `destination' tori $T_J$ have the same action ($J=0.001$ in this picture, corresponding to the translucent plane), $\SM_1$ maps the $4$ points $x^i_- \in T_J$ to $x^i_+ \in T_J$. Fixing the source torus $T_J$ and varying the destination torus $T_{J'}$ to action levels $J'$ near $J$, the scattering map $\SM_1$ is continued to the purple and green sets.}
	\label{fig:cont_2_2}
\end{figure}


\begin{figure}
	\scalebox{0.9}{\input{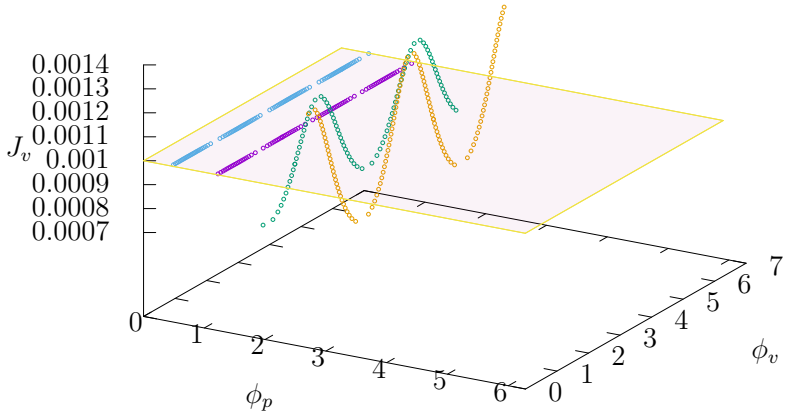}}
	\caption{Action of scattering maps $\SM_1$ and $\SM_2$ on the torus $T_J$ with $J=0.001$. $\SM_1$ maps the purple set to the green set, while $\SM_2$ maps blue to orange.}
	\label{fig:cont_2_2_SM1_SM2}
\end{figure}


%

To compute $\SM_1$ on a given torus $T_J$, fix an action 
$J\in[0.001,0.007]$. As explained above, $W^u(T_J)$ has $8$ transverse intersections with $W^s(T_J)$ in the surface of section $\mathscr{S}$. Four of them, which we will denote $x^i$ for $i=1,2,3,4$, give rise to $4$ pairs of points associated by the local scattering maps:
\[ x^i_+ = \SM(x^i_-) \qquad \text{ for } i=1,\dotsc,4.  \]
Figure~\ref{fig:cont_2_2} shows these $4$ pairs of points: $x^i_-$ in the domain of $\SM_1$ are plotted in purple, while $x^i_+$ in the codomain are plotted in green.

Keeping the action of the source torus fixed to $J$, vary the action of the destination torus to a new value $J'$ close to $J$, and recompute the intersection $W^u(T_J) \cap W^s(T_{J'})$, giving rise to $4$ new homoclinic points and $4$ corresponding pairs of points associated by $\SM_1$. Continue this procedure until the manifolds $W^u(T_J)$ and $W^s(T_{J'})$ cease to intersect, effectively extending the domain of the scattering map $\SM_1$ from $4$ points to the purple set in Figure~\ref{fig:cont_2_2}, and the codomain to the green set.

\MR{The continuation procedure is actually performed in two directions: First increase $J'$ from $J$, producing the portion of the green set located above the translucent plane, until the manifolds cease to intersect ($J'\approx 0.00125$ in the figure). Then decrease $J'$ from $J$, producing the portion below the plane, until the manifolds cease to intersect ($J'\approx 0.0008$ in the figure).}

\begin{remark}
	Starting with the other $4$ homoclinic points, $x^i$ for $i=5,6,7,8$, and applying the same procedure, gives rise to a different scattering map $\SM_2$. See Figure~\ref{fig:cont_2_2_SM1_SM2}.
	$\SM_1$ and $\SM_2$ are fundamentally different, as they are not conjugated by the flow  in the sense of \eqref{eq:conjugacy}.
\end{remark}

\begin{figure}
	\scalebox{0.9}{\input{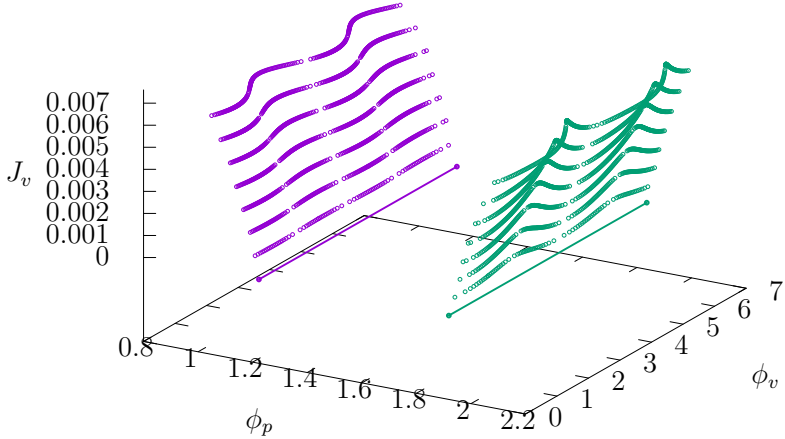}}
	\caption{$\SM_1$ acting on several action levels: $J=0.001, 0.002, \dotsc, 0.007$.}
	\label{fig:SM1}
\end{figure}

We repeat this procedure for the tori $J=0.001, 0.002, \dotsc, 0.007$. See Figure~\ref{fig:SM1}.

\begin{remark}
	For $J\to 0$, the torus $T_J$ degenerates into the horizontal Lyapunov periodic orbit $T_0$. The $4$ homoclinic points $x^i$ converge to a single homoclinic point $x$ for $T_0$. Similarly, the $4$ pairs $x^i_-, x^i_+$ converge to a single pair $x_-, x_+$ of points associated by $\SM_1$. Figure~\ref{fig:SM1} shows $x_-, x_+$ as straight lines (all angles $\phi_v\in[0,2\pi)$ are identified for $J=0$).
\end{remark}

\MR{Notice that one could extend $\SM_1$ from the purple set to the whole NHIM $\NHIM_c$ using the conjugacy property of the scattering map by the flow  \eqref {eq:conjugacy}.
However, we will use the reduced scattering map $\sm_1$ instead.}
\MB{Using the computation of $\SM_1$ we compute the reduced scattering map $\sm_1$.}
Recall from Section~\ref{sec:reduction} that the reduced scattering map can be obtained from the full scattering map simply by flowing the points $x_-$/$x_+$ backwards/forwards to the Poincar\'e section $\Sigma = \{ \phi_p=0 \}$.

Flowing the purple set backwards to $\Sigma$, we obtain a mesh $\left\{(J,\phi)\right\}$ discretizing the reduced NHIM $\NHIMsec$.
Flowing the green set forwards to $\Sigma$, we obtain the image set $\left\{(J',\phi')\right\}$ under $\sm_1$, also on $\NHIMsec$. This way, we have extended the local scattering maps onto a global (reduced) scattering map $\sm_1$ on $\NHIMsec$.

Figure~\ref{fig:curves} (top panel) shows the image set of the global scattering map $\sm_1$. For example, the green set in Figure~\ref{fig:cont_2_2} ($J=0.001$) corresponds to the lowest curve in Figure~\ref{fig:curves} after flowing it forward.


\subsection{Domain of the Global Scattering Map}\label{Sec:DomainGlobalSM}

From above, the global scattering maps $\sm_1, \sm_2$ are well defined on an annulus $\domain$ inside $\NHIMsec$ given by
\begin{equation} \label{eq:domainJ}
\domain = \{ (J, \phi)\colon\ J\in[0.001,0.007] \text{ and } \phi\in[0,2\pi) \}.
\end{equation}

\begin{figure}
	\scalebox{0.9}{\input{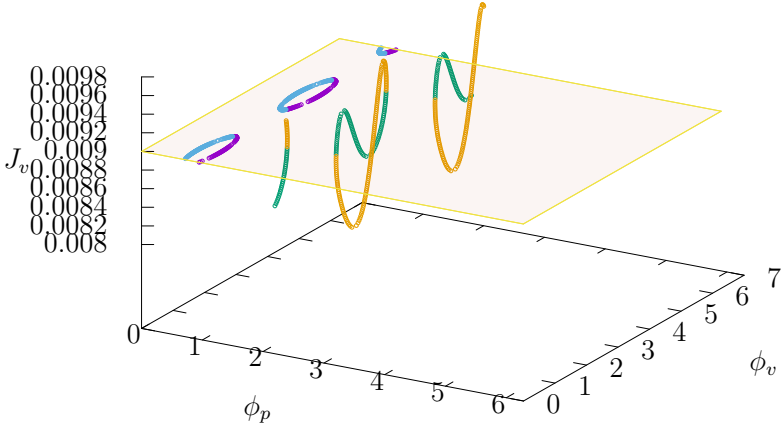}}
	\caption{Action of scattering maps $\SM_1$ and $\SM_2$ on the action level $J=0.009$. $\SM_1$ and $\SM_2$ have `merged', and they can not be extended to the whole torus.}
	\label{fig:cont_10_10}
\end{figure}

For $0.007<J<0.01$, we find that the scattering maps can not be defined on the whole torus $T_J$. For illustration, Figure~\ref{fig:cont_10_10} shows the continuation of all $8$ local scattering maps when $J=0.009$. Notice that the purple and blue sets (which belong to the domain of $\SM_1$ and $\SM_2$ respectively) have become connected, and they form two contractible circles. Moreover, the purple and blue sets do not cover all angles $\phi_v\in[0, 2\pi)$ as before. Thus one can not extend $\SM_i$ to the whole torus $T_J$ by the flow using the conjugacy property~\eqref{eq:conjugacy}.
Consequently, the reduced scattering map $\sm_i$ is not defined on the whole invariant curve $T^\Sigma_{J}$.

For the purpose of this paper, we will restrict the domain of the global scattering maps $\sm_1, \sm_2$ to the annulus $\domain$.


\section{Series Representation of the Scattering Map} \label{Sec:SeriesRepresentation}

The goal of this section is to introduce a series representation of the (global, reduced) scattering maps $\sm_1$ and $\sm_2$. This finite series expansion consists in Equations~\eqref{eq:sm}-\eqref{eq:Taylor} and~\eqref{eq:NewtonOmega} below. It is a more efficient representation than the numerical scattering map computed in the previous section, since it is limited to a small number of terms. Moreover, it allows us to evaluate $\sm_i$, $i=1,2$, at any point of its domain.


Note that $J$ and $\phi$ have different scales: $J$ is of order $10^{-3}$, while $\phi$ is of order $1$. To improve numerical conditioning, we scale  $J$ as follows:
\begin{equation}\label{eq:scaled} I = 1000J. \end{equation}
From now on, we will work with the scaled coordinate $I$ instead of $J$.

The most classical way to represent a symplectic map $(I,\phi)\to(I',\phi')$ is by a generating function depending on old and new variables. 
In~\cite{DMR08}, the scattering map on  Lyapunov periodic orbits was shown to be a phase shift $(I,\phi)\to(I,\phi+\Delta(I))$. Thus, in our setting it is natural to look for a generating function of the form
\[ \gf(I,\phi') = I\phi' + \Omega(I) + \gft(I,\phi'), \]
which will at least be valid for small values of $I$. The generating function $\gf(I,\phi')$ is decomposed into its average part with respect to $\phi'$, denoted $\Omega(I)$, and its oscillatory part $\gft(I,\phi')$, which satisfyies  $\int_{0}^{2\pi} \gft(I,\phi')\textrm{d}\phi'=0$.

%

Hence, the equations for the scattering map $(I',\phi')=\sm(I,\phi)$ are given implicitly by
\begin{subequations} \label{eq:sm}
	\begin{align}
	\phi = \frac{\partial\gf}{\partial I}(I,\phi') &= \phi' + \omega(I) + \frac{\partial\gft}{\partial I}(I,\phi') \label{eq:sm:phi}\\
	I' = \frac{\partial\gf}{\partial\phi'}(I,\phi') &= I \hphantom{+ \omega(I)+}+ \frac{\partial\gft}{\partial\phi'}(I,\phi'), \label{eq:sm:I}
	\end{align}
\end{subequations}
where $\omega(I)=\Omega'(I)$.

We will approximate both $\gft(I,\phi')$ and $\omega(I)$ in Equation~\eqref{eq:sm} using finite series expansions.
Firstly, \AB{since we are going to work with few values of $I$,}  we approximate the function $\gft(I,\phi')$ using a Fourier-\AB{Polynomial interpolation:}
\begin{equation} \label{eq:Fourier}
\gft(I,\phi') =
-\sum_{n=1}^N \frac{B_n(I)}{n}\cos{n\phi'}
+\sum_{n=1}^N \frac{A_n(I)}{n}\sin{n\phi'},
\end{equation}
where
\begin{equation} \label{eq:Taylor}
A_n(I) = \sum_{l=0}^{L} a_l^{(n)} I^l \quad \text{and} \quad
B_n(I) = \sum_{l=0}^{L} b_l^{(n)} I^l.
\end{equation}


\begin{figure}
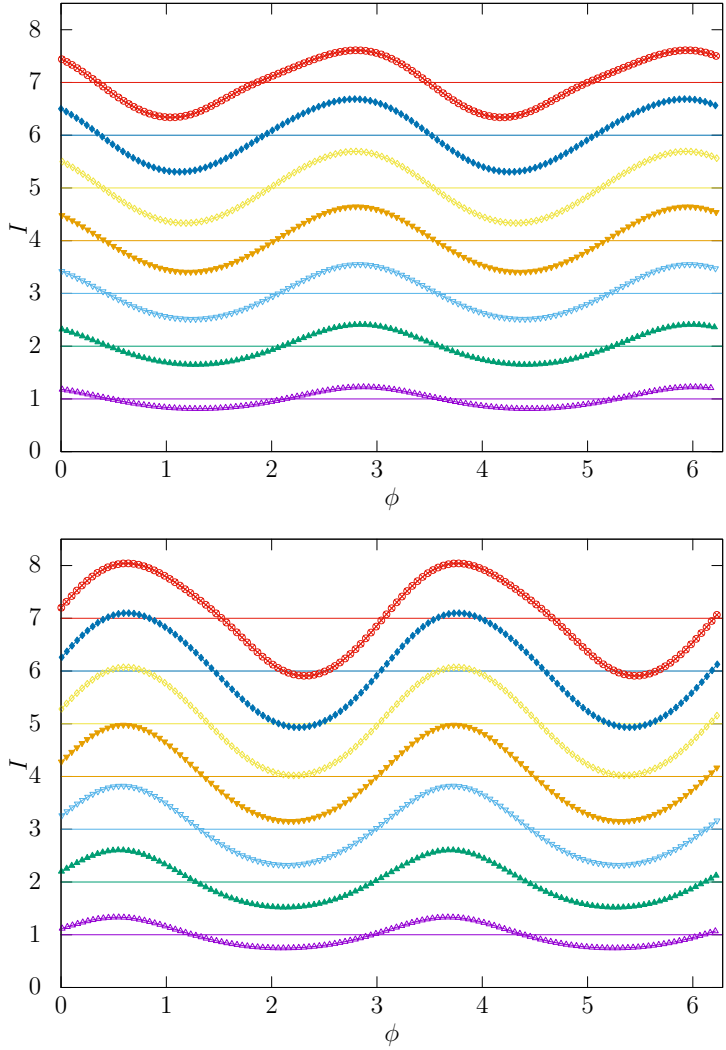

	\scalebox{0.8}{\input{curves}}
	\scalebox{0.8}{\input{curves_SM2}}
	\caption{Image under the scattering map of several tori ($I=1, 2, \dotsc, 7$). Above: action of $\sm_1$, below: action of $\sm_2$. A torus $\{(I,\phi): I=\text{const}, \phi\in[0,2\pi)\}$ and its image are plotted using the same color. Notice that the curves are $\pi$-periodic.}
	\label{fig:curves}
\end{figure}

The goal is to find the coefficients $a_l^{(n)}$ and $b_l^{(n)}$. In Section~\ref{Sec:NumericalResults}, we obtained numerically the scattering map $\sm$ on a grid of equispaced $(I,\phi)$ points. See Figure~\ref{fig:curves}.
It is a simple matter to fit the coefficients to this data, as explained in Section~\ref{Sec:FTApproximation}.

\begin{remark}\label{rem:symmetry}
	The spatial RTBP is invariant with respect to the transformation $(Z,\dot{Z})\to (-Z,-\dot{Z})$. Thus, every trajectory passing through the point $(X,Y,Z)$ has a \textit{symmetric trajectory with respect to the $XY$ plane}, which passes through the point $(X,Y,-Z)$. In particular, every heteroclinic trajectory from torus $I$ to torus $I'$ has a symmetric heteroclinic trajectory from torus $I$ to $I'$.
	In Birkhoff normal form coordinates, the symmetry $Z\to -Z$ corresponds to $\phi_v\to \phi_v+\pi$. This translates to the following fact for the scattering map of the flow. Suppose that $\SM(J_v, \phi_p, \phi_v) = (J_v', \phi_p', \phi_v')$. Then we have $\SM(J_v, \phi_p, \phi_v+\pi) = (J_v', \phi_p', \phi_v'+\pi)$.
	Equivalently for the reduced scattering map: Suppose that $\sm(I, \phi_v) = (I', \phi_v')$. Then we have $\sm(I, \phi_v+\pi) = (I', \phi_v'+\pi)$. This implies that the image under the scattering map of any torus is a $\pi$-periodic curve (see Figure~\ref{fig:curves}). Indeed, we have checked that the curves in Figure~\ref{fig:curves} are $\pi$-periodic up to a tolerance of $10^{-6}$.
	From this point on, \textbf{we will plot all figures involving $\phi_v$ in the domain $[0,\pi)$ only}.
\end{remark}

\subsection{Fourier-\AB{Polynomial interpolation} of the Generating Function}
\label{Sec:FTApproximation}

Let $I$ be fixed, and consider equation~\eqref{eq:sm:I}:
\[ I' = I + \pd{\gft}{\phi'}, \]
where
\begin{equation}
\label{eq:dgft/dphi}
\pd{\gft}{\phi'} =
\sum_{n=1}^N A_n\cos{n\phi'}
+\sum_{n=1}^N B_n\sin{n\phi'}.
\end{equation}
Given a set of $(I',\phi')$ values on a grid (data points composing one curve in Figure~\ref{fig:curves}), we use the discrete Fourier transform to obtain the Fourier coefficients $A_n$, $B_n$.

\begin{remark}
	Since we have 128 $(I',\phi')$ data points for each torus, the maximum possible degree of the Fourier expansion~\eqref{eq:Fourier} is $N=64$.
\end{remark}

\begin{figure}
	\includegraphics[width=.5\linewidth]{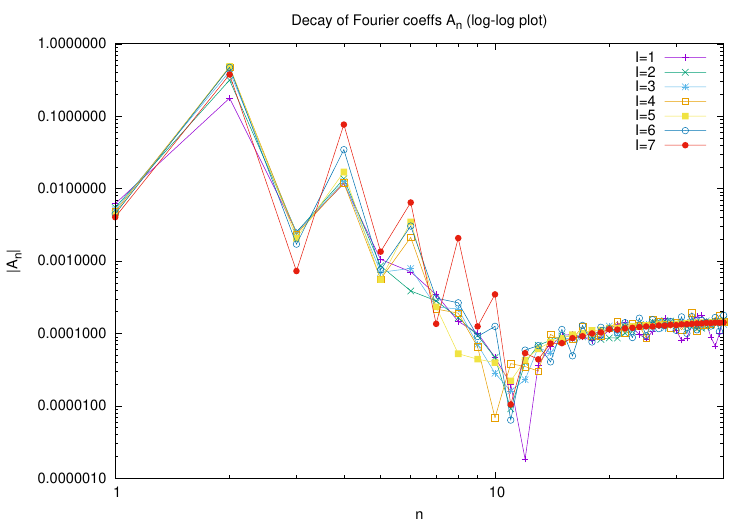}
	\includegraphics[width=.5\linewidth]{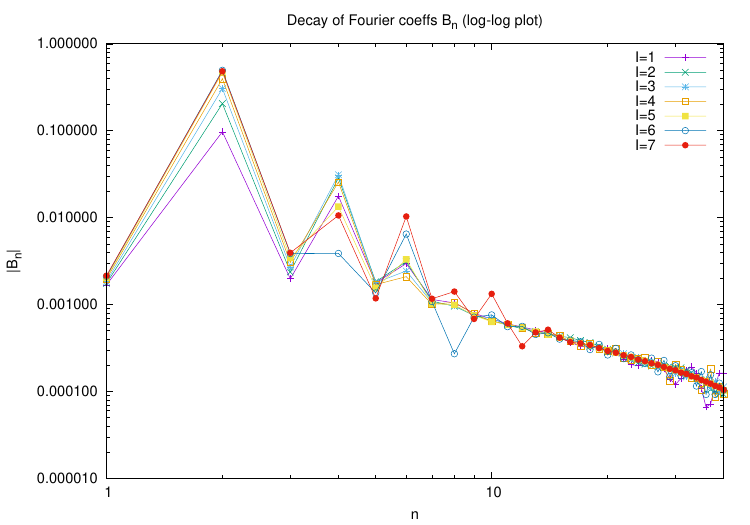}
	\includegraphics[width=.5\linewidth]{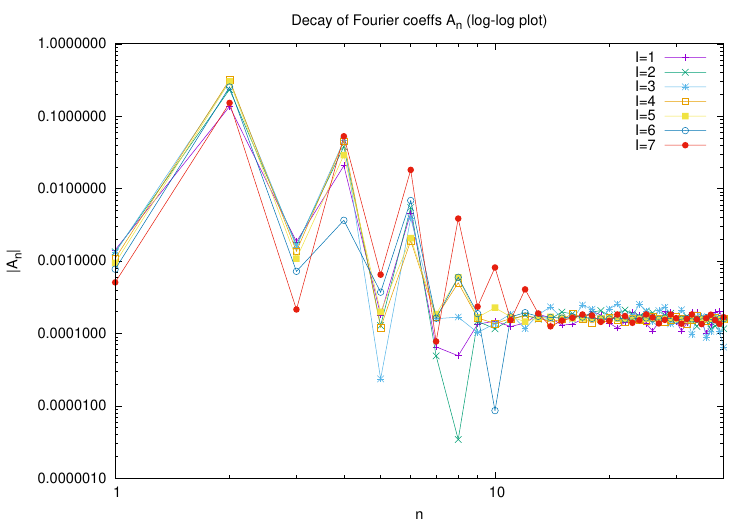}
	\includegraphics[width=.5\linewidth]{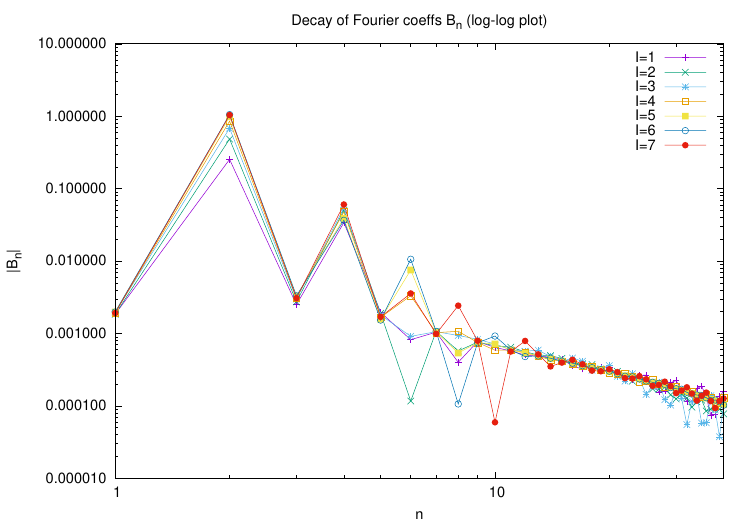}
	\caption{Decay of Fourier coefficients (log-log plot) for the scattering map $\sm_1$ (above) and $\sm_2$ (below).}
	\label{fig:decay}
\end{figure}

Next, let $I$ vary to obtain the Fourier coefficients $A_n(I)$, $B_n(I)$ for each torus $I=1,2, \dotsc, 7$. Figure~\ref{fig:decay} shows the decay of these Fourier coefficients for each torus.

\begin{remark}\label{rem:evenFourier}
	We know from Remark~\ref{rem:symmetry} that the scattering map is $\pi$-periodic.
	Thus, for each torus $I=const$, its image is a $\pi$-periodic curve $I'=\gamma(\phi')$, and the odd Fourier coefficients should all be zero: $A_{2k+1}(I) = B_{2k+1}(I) = 0$ for $k=0,1,2, \dotsc$. Of course, in Figure~\ref{fig:decay} they are not exactly zero because the numerical data is not exactly $\pi$-periodic. From now on, all odd Fourier coefficients are set to zero, to obtain a model that satisfies the theoretical $\pi$-periodicity of the scattering map.
\end{remark}

Let now the index $n$ of the Fourier coefficient be fixed, and consider the \AB{polynomial} approximation of $A_n(I)$ and $B_n(I)$ given in equation~\eqref{eq:Taylor}.
Alternatively, given a set of $L+1$ data points $(I_0, A_n(I_0)),\allowbreak \dotsc,\allowbreak  (I_L, A_n(I_L))$, we will express the polynomials~\eqref{eq:Taylor} in Newton's form
\begin{equation} \label{eq:Newton}
A_n(I) = \sum_{l=0}^{L} \tilde{a}_l^{(n)} N_l(I) \quad \text{and} \quad
B_n(I) = \sum_{l=0}^{L} \tilde{b}_l^{(n)} N_l(I),
\end{equation}
where $\tilde{a}_l^{(n)}, \tilde{b}_l^{(n)}$ are the \textit{divided differences}, and $N_l(I)$ are Newton's basis polynomials
\begin{equation*} \label{eq:NewtonBasis}
N_0(I)=1, \quad N_l(I) \coloneqq \prod_{i=0}^{l-1} (I-I_i) \ \text{ for }l=1,\dotsc,L.
\end{equation*}

\begin{table}\footnotesize
	\begin{subtable}{\textwidth}
		\pgfplotstableread[col sep=comma]{ddA_filtered.csv}\ddAtable
		
		\pgfplotstabletypeset[
		string type,
		columns/zero/.style={column name=\bm{$n$}},
		columns/one/.style={column name=\bm{ $\tilde{a}_1^{(n)}$}},
		columns/two/.style={column name=\bm{ $\tilde{a}_2^{(n)}$}},
		columns/three/.style={column name=\bm{ $\tilde{a}_3^{(n)}$}},
		columns/four/.style={column name=\bm{ $\tilde{a}_4^{(n)}$}},
		columns/five/.style={column name=\bm{ $\tilde{a}_5^{(n)}$}},
		columns/six/.style={column name=\bm{ $\tilde{a}_6^{(n)}$}},
		columns/seven/.style={column name=\bm{ $\tilde{a}_7^{(n)}$}},
		every row 0 column 1/.style={postproc cell content/.style={@cell content=\textbf{\textcolor{red}{##1}}}},
		every row 0 column 2/.style={postproc cell content/.style={@cell content=\textbf{\textcolor{blue}{##1}}}},
		every head row/.style={before row=\toprule,after row=\midrule},
		every last row/.style={after row={\bottomrule}},
		]\ddAtable
	\end{subtable}
	
	\begin{subtable}{\textwidth}
		\pgfplotstableread[col sep=comma]{ddB_filtered.csv}\ddBtable
		
		\pgfplotstabletypeset[
		string type,
		columns/zero/.style={column name=\bm{$n$}},
		columns/one/.style={column name=\bm{ $\tilde{b}_1^{(n)}$}},
		columns/two/.style={column name=\bm{ $\tilde{b}_2^{(n)}$}},
		columns/three/.style={column name=\bm{ $\tilde{b}_3^{(n)}$}},
		columns/four/.style={column name=\bm{ $\tilde{b}_4^{(n)}$}},
		columns/five/.style={column name=\bm{ $\tilde{b}_5^{(n)}$}},
		columns/six/.style={column name=\bm{ $\tilde{b}_6^{(n)}$}},
		columns/seven/.style={column name=\bm{ $\tilde{b}_7^{(n)}$}},
		every row 0 column 1/.style={postproc cell content/.style={@cell content=\textbf{\textcolor{red}{##1}}}},
		every head row/.style={before row=\toprule,after row=\midrule},
		every last row/.style={after row={\bottomrule}},
		]\ddBtable
		
		\caption{Scattering map $\sm_1$.}
	\end{subtable}	
	
	\vspace{0.5cm}
	
	\begin{subtable}{\textwidth}
		\pgfplotstableread[col sep=comma]{ddA_SM2_filtered.csv}\ddASMtwotable
		
		\pgfplotstabletypeset[
		string type,
		columns/zero/.style={column name=\bm{$n$}},
		columns/one/.style={column name=\bm{ $\tilde{a}_1^{(n)}$}},
		columns/two/.style={column name=\bm{ $\tilde{a}_2^{(n)}$}},
		columns/three/.style={column name=\bm{ $\tilde{a}_3^{(n)}$}},
		columns/four/.style={column name=\bm{ $\tilde{a}_4^{(n)}$}},
		columns/five/.style={column name=\bm{ $\tilde{a}_5^{(n)}$}},
		columns/six/.style={column name=\bm{ $\tilde{a}_6^{(n)}$}},
		columns/seven/.style={column name=\bm{ $\tilde{a}_7^{(n)}$}},
		every row 0 column 1/.style={postproc cell content/.style={@cell content=\textbf{\textcolor{red}{##1}}}},
		every row 0 column 2/.style={postproc cell content/.style={@cell content=\textbf{\textcolor{blue}{##1}}}},
		every head row/.style={before row=\toprule,after row=\midrule},
		every last row/.style={after row={\bottomrule}},
		]\ddASMtwotable
	\end{subtable}
	
	\begin{subtable}{\textwidth}
		\pgfplotstableread[col sep=comma]{ddB_SM2_filtered.csv}\ddBSMtwotable
		
		\pgfplotstabletypeset[
		string type,
		columns/zero/.style={column name=\bm{$n$}},
		columns/one/.style={column name=\bm{ $\tilde{b}_1^{(n)}$}},
		columns/two/.style={column name=\bm{ $\tilde{b}_2^{(n)}$}},
		columns/three/.style={column name=\bm{ $\tilde{b}_3^{(n)}$}},
		columns/four/.style={column name=\bm{ $\tilde{b}_4^{(n)}$}},
		columns/five/.style={column name=\bm{ $\tilde{b}_5^{(n)}$}},
		columns/six/.style={column name=\bm{ $\tilde{b}_6^{(n)}$}},
		columns/seven/.style={column name=\bm{ $\tilde{b}_7^{(n)}$}},
		every row 0 column 1/.style={postproc cell content/.style={@cell content=\textbf{\textcolor{red}{##1}}}},
		every row 0 column 2/.style={postproc cell content/.style={@cell content=\textbf{\textcolor{blue}{##1}}}},
		every head row/.style={before row=\toprule,after row=\midrule},
		every last row/.style={after row={\bottomrule}},
		]\ddBSMtwotable
		
		\caption{Scattering map $\sm_2$.}
	\end{subtable}	
	
	\caption[Divided Differences]{Divided differences $\tilde{a}_l^{(n)}$, $\tilde{b}_l^{(n)}$ of the first  Fourier coefficients $A_n(I)$, $B_n(I)$. All odd Fourier coefficients are zero and thus not listed (see Remark~\ref{rem:evenFourier}). All constant terms of the Newton expansion $\tilde{a}_0^{(n)}$, $\tilde{b}_0^{(n)}$ are zero, and thus not listed (see Remark~\ref{rem:FourierI0}). 
		Notice that the two coefficients $\tilde{a}_1^2$ and $\tilde{b}_1^2$ (in red) are much larger than the rest of coefficients for both scattering maps $\sigma_1$, $\sigma_2$. 
		For the scattering map $\sigma_1$, the coefficient $\tilde{a}_2^2$ (in blue) is much larger than the rest of coefficients in the last six columns, i.e. $\tilde{a}_j^n$ for $j\geq 2$. 
		For the scattering map $\sigma_2$, the two coefficients $\tilde{a}_2^2$ and $\tilde{b}_2^2$ are much larger than the rest with $j\geq 2$ .}
	\label{tab:dd}
\end{table}

Given the values $A_n(I)$ at $I=0,1,2, \dotsc, 7$ obtained in the previous step (resp. $B_n(I)$), we use polynomial interpolation (Newton's `divided differences' algorithm) to obtain the coefficients $\tilde{a}_l^{(n)}$ (resp. $\tilde{b}_l^{(n)}$).
The divided differences of the first $8$ Fourier coefficients are listed in Table~\ref{tab:dd}.


\begin{remark}
	Since we only have 8 data points for each Fourier coefficient, the maximum possible degree of the Newton expansion~\eqref{eq:Newton} is $L=7$.
\end{remark}

\begin{remark}\label{rem:FourierI0}
	The constant terms $\tilde{a}_0^{(n)}, \tilde{b}_0^{(n)}$ of the Newton expansion are all zero, since $\tilde{a}_0^{(n)} = A_n(0) = 0$ and  $\tilde{b}_0^{(n)} = B_n(0) = 0$ for all $n$.
\end{remark}



\subsection{\AB{Polynomial interpolation} of the Frequency \texorpdfstring{$\omega(I)$}{omega(I)}}

Finally, we approximate the frequency $\omega(I)$ in Equation~\eqref{eq:sm}.

\subsubsection*{Frequency $\omega(I)$ at $I=1,2,\dotsc,7$.}

For action levels $I=1,2,\dotsc,7$ we have available the scattering map data $(I,\phi)\to(I',\phi')$, so the frequency $\omega(I)$ can be obtained from Equation~\eqref{eq:sm:phi} as
\begin{equation} \label{eq:omega}
\omega(I) = \phi - \phi' - \pd{\gft}{I}(I,\phi').
\end{equation}
In an exact calculation, $\omega(I)$ should be independent of $\phi$. Numerically, using Equation~\eqref{eq:omega} to approximate the value of $\omega(I)$ would yield slightly different values for $\omega(I)$ depending on the data point $(I,\phi)$, and indeed depending on $\phi$. Thus, we will compute $\omega(I)$ as the \textit{average} of Equation~\eqref{eq:omega} over all $\phi\in[0,2\pi)$.



\PR{To give an idea of the numerical error in determining $\omega(I)$, Figure~\ref{fig:omega} compares the frequency $\omega(2,\phi')$ of torus $I=2$ to its average $\average{\omega}$ over all available $\phi'$ values. Their discrepancy is less than $0.002$ radians.}

\subsubsection*{Frequency $\omega(I)$ at an arbitrary $I$-value.}

Similarly to what we did in Section~\ref{Sec:FTApproximation} with the Fourier coefficients $A_n(I)$, $B_n(I)$, we use the Newton series representation of $\omega(I)$:
\begin{equation} \label{eq:NewtonOmega}
\omega(I) = \sum_{l=0}^{L} \tilde{c}_l N_l(I).
\end{equation}

\begin{figure}
	\begin{tikzpicture}[gnuplot]
\path (0.000,0.000) rectangle (12.500,8.750);
\gpcolor{color=gp lt color border}
\gpsetlinetype{gp lt border}
\gpsetdashtype{gp dt solid}
\gpsetlinewidth{1.00}
\draw[gp path] (1.320,0.985)--(1.500,0.985);
\draw[gp path] (11.947,0.985)--(11.767,0.985);
\node[gp node right] at (1.136,0.985) {$2$};
\draw[gp path] (1.320,1.917)--(1.500,1.917);
\draw[gp path] (11.947,1.917)--(11.767,1.917);
\node[gp node right] at (1.136,1.917) {$2.2$};
\draw[gp path] (1.320,2.849)--(1.500,2.849);
\draw[gp path] (11.947,2.849)--(11.767,2.849);
\node[gp node right] at (1.136,2.849) {$2.4$};
\draw[gp path] (1.320,3.781)--(1.500,3.781);
\draw[gp path] (11.947,3.781)--(11.767,3.781);
\node[gp node right] at (1.136,3.781) {$2.6$};
\draw[gp path] (1.320,4.713)--(1.500,4.713);
\draw[gp path] (11.947,4.713)--(11.767,4.713);
\node[gp node right] at (1.136,4.713) {$2.8$};
\draw[gp path] (1.320,5.645)--(1.500,5.645);
\draw[gp path] (11.947,5.645)--(11.767,5.645);
\node[gp node right] at (1.136,5.645) {$3$};
\draw[gp path] (1.320,6.577)--(1.500,6.577);
\draw[gp path] (11.947,6.577)--(11.767,6.577);
\node[gp node right] at (1.136,6.577) {$3.2$};
\draw[gp path] (1.320,7.509)--(1.500,7.509);
\draw[gp path] (11.947,7.509)--(11.767,7.509);
\node[gp node right] at (1.136,7.509) {$3.4$};
\draw[gp path] (1.320,8.441)--(1.500,8.441);
\draw[gp path] (11.947,8.441)--(11.767,8.441);
\node[gp node right] at (1.136,8.441) {$3.6$};
\draw[gp path] (1.320,0.985)--(1.320,1.165);
\draw[gp path] (1.320,8.441)--(1.320,8.261);
\node[gp node center] at (1.320,0.677) {$1$};
\draw[gp path] (3.091,0.985)--(3.091,1.165);
\draw[gp path] (3.091,8.441)--(3.091,8.261);
\node[gp node center] at (3.091,0.677) {$2$};
\draw[gp path] (4.862,0.985)--(4.862,1.165);
\draw[gp path] (4.862,8.441)--(4.862,8.261);
\node[gp node center] at (4.862,0.677) {$3$};
\draw[gp path] (6.634,0.985)--(6.634,1.165);
\draw[gp path] (6.634,8.441)--(6.634,8.261);
\node[gp node center] at (6.634,0.677) {$4$};
\draw[gp path] (8.405,0.985)--(8.405,1.165);
\draw[gp path] (8.405,8.441)--(8.405,8.261);
\node[gp node center] at (8.405,0.677) {$5$};
\draw[gp path] (10.176,0.985)--(10.176,1.165);
\draw[gp path] (10.176,8.441)--(10.176,8.261);
\node[gp node center] at (10.176,0.677) {$6$};
\draw[gp path] (11.947,0.985)--(11.947,1.165);
\draw[gp path] (11.947,8.441)--(11.947,8.261);
\node[gp node center] at (11.947,0.677) {$7$};
\draw[gp path] (1.320,8.441)--(1.320,0.985)--(11.947,0.985)--(11.947,8.441)--cycle;
\node[gp node center,rotate=-270] at (0.292,4.713) {$\omega(I)$};
\node[gp node center] at (6.633,0.215) {$I$};
\gpcolor{rgb color={0.580,0.000,0.827}}
\gpsetpointsize{12.00}
\gppoint{gp mark 1}{(1.320,1.114)}
\gppoint{gp mark 1}{(3.091,1.332)}
\gppoint{gp mark 1}{(4.862,1.560)}
\gppoint{gp mark 1}{(6.634,1.813)}
\gppoint{gp mark 1}{(8.405,2.099)}
\gppoint{gp mark 1}{(10.176,2.431)}
\gppoint{gp mark 1}{(11.947,2.844)}
\gpcolor{rgb color={0.000,0.620,0.451}}
\draw[gp path] (1.320,1.114)--(1.497,1.136)--(1.674,1.158)--(1.851,1.180)--(2.028,1.201)%
  --(2.206,1.223)--(2.383,1.245)--(2.560,1.266)--(2.737,1.288)--(2.914,1.310)--(3.091,1.332)%
  --(3.268,1.354)--(3.445,1.376)--(3.623,1.398)--(3.800,1.421)--(3.977,1.444)--(4.154,1.466)%
  --(4.331,1.489)--(4.508,1.513)--(4.685,1.536)--(4.862,1.560)--(5.039,1.584)--(5.217,1.608)%
  --(5.394,1.633)--(5.571,1.658)--(5.748,1.683)--(5.925,1.708)--(6.102,1.734)--(6.279,1.760)%
  --(6.456,1.786)--(6.634,1.813)--(6.811,1.840)--(6.988,1.867)--(7.165,1.895)--(7.342,1.923)%
  --(7.519,1.951)--(7.696,1.980)--(7.873,2.009)--(8.050,2.039)--(8.228,2.069)--(8.405,2.099)%
  --(8.582,2.130)--(8.759,2.161)--(8.936,2.193)--(9.113,2.225)--(9.290,2.258)--(9.467,2.291)%
  --(9.644,2.325)--(9.822,2.359)--(9.999,2.394)--(10.176,2.430)--(10.353,2.467)--(10.530,2.504)%
  --(10.707,2.542)--(10.884,2.580)--(11.061,2.620)--(11.239,2.660)--(11.416,2.702)--(11.593,2.744)%
  --(11.770,2.787)--(11.947,2.832);
\gpcolor{rgb color={0.337,0.706,0.914}}
\gppoint{gp mark 3}{(1.320,8.236)}
\gppoint{gp mark 3}{(3.091,8.020)}
\gppoint{gp mark 3}{(4.862,7.788)}
\gppoint{gp mark 3}{(6.634,7.537)}
\gppoint{gp mark 3}{(8.405,7.249)}
\gppoint{gp mark 3}{(10.176,6.921)}
\gppoint{gp mark 3}{(11.947,6.500)}
\gpcolor{rgb color={0.902,0.624,0.000}}
\draw[gp path] (1.320,8.236)--(1.497,8.217)--(1.674,8.197)--(1.851,8.176)--(2.028,8.154)%
  --(2.206,8.132)--(2.383,8.110)--(2.560,8.088)--(2.737,8.065)--(2.914,8.042)--(3.091,8.020)%
  --(3.268,7.997)--(3.445,7.974)--(3.623,7.951)--(3.800,7.928)--(3.977,7.905)--(4.154,7.881)%
  --(4.331,7.858)--(4.508,7.835)--(4.685,7.811)--(4.862,7.788)--(5.039,7.764)--(5.217,7.740)%
  --(5.394,7.716)--(5.571,7.691)--(5.748,7.666)--(5.925,7.641)--(6.102,7.615)--(6.279,7.589)%
  --(6.456,7.563)--(6.634,7.537)--(6.811,7.510)--(6.988,7.482)--(7.165,7.454)--(7.342,7.426)%
  --(7.519,7.398)--(7.696,7.369)--(7.873,7.339)--(8.050,7.310)--(8.228,7.280)--(8.405,7.249)%
  --(8.582,7.218)--(8.759,7.187)--(8.936,7.156)--(9.113,7.124)--(9.290,7.091)--(9.467,7.058)%
  --(9.644,7.025)--(9.822,6.991)--(9.999,6.957)--(10.176,6.922)--(10.353,6.886)--(10.530,6.849)%
  --(10.707,6.811)--(10.884,6.772)--(11.061,6.732)--(11.239,6.690)--(11.416,6.646)--(11.593,6.600)%
  --(11.770,6.552)--(11.947,6.501);
\gpcolor{color=gp lt color border}
\draw[gp path] (1.320,8.441)--(1.320,0.985)--(11.947,0.985)--(11.947,8.441)--cycle;
\gpdefrectangularnode{gp plot 1}{\pgfpoint{1.320cm}{0.985cm}}{\pgfpoint{11.947cm}{8.441cm}}
\end{tikzpicture}
	\caption{Frequency function $\omega(I)$ for $\sm_1$ (in green) and $\sm_2$ (yellow).}
	\label{fig:interp_omega}
\end{figure}

\begin{table}\footnotesize
	\begin{subtable}{\textwidth}
		\pgfplotstableread[col sep=comma]{ddOmega.csv}\ddOmegatable
		
		\pgfplotstabletypeset[
		string type,
		columns/zero/.style={column name=\bm{$\tilde{c}_0$}},
		columns/one/.style={column name=\bm{ $\tilde{c}_1$}},
		columns/two/.style={column name=\bm{ $\tilde{c}_2$}},
		columns/three/.style={column name=\bm{ $\tilde{c}_3$}},
		columns/four/.style={column name=\bm{ $\tilde{c}_4$}},
		columns/five/.style={column name=\bm{ $\tilde{c}_5$}},
		columns/six/.style={column name=\bm{ $\tilde{c}_6$}},
		every row 0 column 0/.style={postproc cell content/.style={@cell content=\textbf{\textcolor{red}{##1}}}},
		every row 0 column 1/.style={postproc cell content/.style={@cell content=\textbf{\textcolor{blue}{##1}}}},
		every head row/.style={before row=\toprule,after row=\midrule},
		every last row/.style={after row={\bottomrule}},
		]\ddOmegatable
		
		\caption{Scattering map $\sm_1$.}
	\end{subtable}	
	
	\vspace{0.5cm}
	
	\begin{subtable}{\textwidth}
		\pgfplotstableread[col sep=comma]{ddOmega_SM2.csv}\ddOmegaSMtwotable
		
		\pgfplotstabletypeset[
		string type,
		columns/zero/.style={column name=\bm{$\tilde{c}_0$}},
		columns/one/.style={column name=\bm{ $\tilde{c}_1$}},
		columns/two/.style={column name=\bm{ $\tilde{c}_2$}},
		columns/three/.style={column name=\bm{ $\tilde{c}_3$}},
		columns/four/.style={column name=\bm{ $\tilde{c}_4$}},
		columns/five/.style={column name=\bm{ $\tilde{c}_5$}},
		columns/six/.style={column name=\bm{ $\tilde{c}_6$}},
		every row 0 column 0/.style={postproc cell content/.style={@cell content=\textbf{\textcolor{red}{##1}}}},
		every row 0 column 1/.style={postproc cell content/.style={@cell content=\textbf{\textcolor{blue}{##1}}}},
		every head row/.style={before row=\toprule,after row=\midrule},
		every last row/.style={after row={\bottomrule}},
		]\ddOmegaSMtwotable
		\caption{Scattering map $\sm_2$.}
	\end{subtable}	
	
	\caption[Divided Differences]{Divided differences $\tilde{c}_l$ of function $\omega(I)$.
		Notice that the coefficient $\tilde{c}_0$ (in red) is much larger than the rest of coefficients of this table, and larger than those of Table~\ref{tab:dd}, for both scattering maps $\sigma_1$ and $\sigma_2$. Notice also that $\tilde{c}_1$ (in blue) is much larger than the rest of coefficients of the last five columns, i.e $\tilde{c}_j$ for $j\geq 2$, of this table, and larger than the double of those coefficients of the last 6 columns of Table~\ref{tab:dd} for both scattering maps $\sigma_1$ and $\sigma_2$.}
	\label{tab:ddOmega}
\end{table}

Given the values $\average{\omega}(I)$ at $I=1, 2, \dotsc, 7$, obtained in the previous step, we use polynomial interpolation to obtain the coefficients $\tilde{c}_l$. 
The frequency function $\omega(I)$ is plotted in Figure~\ref{fig:interp_omega}.

\begin{remark}
	Since we only have 7 data points for $\omega$, the maximum possible degree of the Newton expansion~\eqref{eq:NewtonOmega} is $L=6$.
\end{remark}

This completes the series representation of the scattering map, consisting of Equations~\eqref{eq:Fourier}-\eqref{eq:Taylor} and~\eqref{eq:NewtonOmega}.

\subsection{Applying the Scattering Map}

We will  use equations~\eqref{eq:sm} to apply the scattering map $(I',\phi')=\sm(I,\phi)$. \MR{However, it must be stressed that these equations do not give $I', \phi'$ \textit{explicitly} as functions of $I, \phi$. On the contrary, $I', \phi'$ are given \textit{implicitly}.}
\MB{Note that \eqref{eq:sm} give $I', \phi'$ \textit{implicitly} as functions of $I, \phi$.}
However, $\phi'$ can be obtained from Equation~\eqref{eq:sm:phi} as a fixed point of
\[ \phi' = f(\phi'; I, \phi) = \phi - \omega(I) - \pd{\gft}{I}(I,\phi'). \]

We simply use fixed point iteration, starting with the initial approximation $\phi'_0 = \phi - \omega(I)$. We require an absolute error smaller than $10^{-5}$ in the fixed point to stop the iteration. There is no point in requiring higher precision, because the error of our series representation in the angle variable is larger than $10^{-2}$; see Table~\ref{tab:FT_error} (bottom panel).

Once $\phi'$ is known, $I'$ is obtained directly from Equation~\eqref{eq:sm:I}.


\section{Approximation Error}
\label{sec:ApproximationError}

Now we have two different representations of the scattering map:

\begin{itemize}
	\item The \textbf{numerical scattering map} $\sigma(I,\phi)$ was obtained in Section~\ref{Sec:NumericalResults} using Birkhoff normal forms and numerical continuation of the invariant manifolds. It was calculated on a relatively coarse grid of points (Figure~\ref{fig:curves}).

	\item The \PB{\textbf{standard scattering map (SSM)}} $\tilde{\sigma}(I,\phi)$ consists in the Fourier-\AB{Polynomial interpolation}~\eqref{eq:Fourier}-\eqref{eq:Taylor} and~\eqref{eq:NewtonOmega}. The series approximation has been derived in Section~\ref{Sec:SeriesRepresentation} from the numerical map, so it is not as precise. However, it has the advantage that it can be evaluated at any desired point $(I,\phi)$.
\end{itemize}

To measure the quality of the series approximation, we do the following:

\begin{enumerate}
	\item Read the numerical scattering map  from file as a table:
	\[ (I',\phi') = \sigma(I,\phi).\]
	We have its values on a grid of points $(I,\phi)$.
	
	\item Evaluate the \PB{standard scattering map (SSM)}  \textit{on the same grid}:
	\[(\tilde{I'}, \tilde{\phi'}) = \tilde{\sigma}(I,\phi).\]
	
	\item Find the approximation error, defined as the \textit{maximum over all grid points} of
	\[ (\epsilon_I, \epsilon_\phi) = \left(\abs{\tilde{I'}-I'}, \abs{\tilde{\phi'}-\phi'}\right). \]
\end{enumerate}

\PB{
\begin{remark}
To better test the quality of the series approximation, we also compute the approximation error on a new set of data, independent of the original one. That is, we extend the original data set $I=1,2,\dotsc,7$, used to derive $\tilde{\sigma}$ with the new data set $I=0.5, 1.5, 2.5, \dotsc, 6.5$. See~Figure~\ref{fig:check_FourierTaylor}. Thus, only for testing purposes, we use 14 tori. As shown in Tables~\ref{tab:FT_error}, \ref{tab:FT_newerror}, the errors in the independent data set are comparable to those of the original data set.
\end{remark}
}

\begin{figure}
\begin{subfigure}{\textwidth}
\centering
    \scalebox{0.9}{\input{check_FourierTaylor_N2_M2}}
	\caption{$N=2, L=2$}
    \label{fig:check_FourierTaylor_N2_M2}
\end{subfigure}

\bigskip

\begin{subfigure}{\textwidth}
\centering
	\scalebox{0.9}{\input{check_FourierTaylor_N4_M5}}
	\caption{$N=4, L=5$}
	\label{fig:check_FourierTaylor_N4_M5}
\end{subfigure}
\caption{Image of the numerical scattering map (points) versus the standard scattering map (SSM) of degree $N,L$ (lines). \PB{For testing purposes, we check the approximation error both on the original data set $I=1,2,\dotsc,7$ (in purple) and on the independent data set $I=0.5, 1.5, 2.5, \dotsc, 6.5$ (in green).}}
\label{fig:check_FourierTaylor}
\end{figure}

Of course, the approximation error depends on the chosen degree $(N,L)$ of the Fourier-\AB{Polynomial interpolation}. For illustration, Figure~\ref{fig:check_FourierTaylor} compares the quality of a low order versus a high order approximation.

We will distinguish two different settings.
In the \textit{local} setting, one is interested in an accurate representation of the scattering map in a neighborhood of $I=0$, whereas in the \textit{global} setting, one is interested in an accurate representation in the whole domain of the global scattering map.

\subsection{Local Approximation Error}\label{sec:localapprox}

For definiteness, let us fix the local domain to be
\MB{
\begin{equation*} \label{eq:domainIloc} \domainloc = \{ (I, \phi)\colon\ I\in[1,3] \text{ and } \phi\in[0,2\pi) \}.\end{equation*}
}
%

\begin{table}[!h]
	
	\begin{filecontents*}{FT_error_local.csv}
N,zero,one,two,three
2,0.540000,0.134772,0.041698,0.037342
4,0.540000,0.145286,0.019337,0.011640
6,0.540000,0.151078,0.019536,0.013382
\end{filecontents*}
\csvreader[
	head to column names,
	tabular = lcccc,
	table head = \toprule
	& \multicolumn{4}{c}{L} \\
	\cmidrule(r){2-5}
	$N$ & 0 & 1 & 2 & 3 \\
	\midrule,
	table foot = \bottomrule,
	]{FT_error_local.csv}{}{%
		\N & \zero & \one & \two & \three	
	}%
    
		\begin{filecontents*}{T_error_local.csv}
N,zero,one,two
2,0.156427,0.052403,0.015123
4,0.156427,0.049706,0.016311
6,0.156427,0.049634,0.016270
\end{filecontents*}
\csvreader[
	head to column names,
	tabular = lccc,
	table head = \toprule
	& \multicolumn{3}{c}{L} \\
	\cmidrule(r){2-4}
	$N$ & 0 & 1 & 2 \\
	\midrule,
	table foot = \bottomrule,
	]{T_error_local.csv}{}{%
		\N & \zero & \one & \two
	}%
	
	\caption{Local approximation error  $\epsilon_I$ (top panel) and $\epsilon_\phi$ (bottom panel) as a function of $N,L$ for the first scattering map $\sm_1$.}
	\label{tab:FT_error_local}
\end{table}

We have computed the approximation error over the local domain $\domainloc$ as a function of the degree $(N,L)$ of the Fourier-\AB{Polynomial interpolation}; see Table~\ref{tab:FT_error_local}.

\begin{remark}
	Only the grid points $(I,\phi)$ belonging to the local domain (i.e tori $I=1,2,3$) are used in the computation of the local approximation error.
\end{remark}


The error of the Fourier-\AB{Polynomial interpolation} model decreases as $N$ and $L$ increase, but not monotonically. If we want an approximation error $\epsilon = \max\{\epsilon_I, \epsilon_\phi\}$ less than $0.05$, then it is enough to take $N=2$ and $L=2$. Notice that the improvement is mild beyond that point.
Indeed, Figure~\ref{fig:check_FourierTaylor_N2_M2} shows that $N=L=2$ gives a good approximation in the local domain $\domainloc$.

Thus, it is natural to choose $N=L=2$ to obtain an accurate model for the local scattering map.
In fact, as discussed before (Remark~\ref{rem:evenFourier}), we neglect the odd Fourier coefficients $A_1(I)$ and $B_1(I)$ due to the symmetry of the problem, and just keep the even ones $A_2(I)$ and $B_2(I)$.

Therefore, in the local setting, an \textbf{accurate model} for the scattering map is given by the Fourier-\AB{Polynomial interpolation}
\begin{equation*}
\gft(I,\phi') =
-\frac{B_2(I)}{2}\cos{2\phi'}
+\frac{A_2(I)}{2}\sin{2\phi'},
\end{equation*}
where
\begin{align*}
A_2(I) &= \tilde{a}_0^{(2)} + \tilde{a}_1^{(2)} I + \tilde{a}_2^{(2)} I(I-1) \\
B_2(I) &= \tilde{b}_0^{(2)} + \tilde{b}_1^{(2)} I + \tilde{b}_2^{(2)} I(I-1),
\end{align*}
\begin{equation*} \label{eq:TaylorQualitativeOmega}
\omega(I) = \tilde{c}_0 + \tilde{c}_1 (I-1) + \tilde{c}_2 (I-1)(I-2),
\end{equation*}
\AB{are polynomials of degree 2 in $I$, determined  by} only $9$ coefficients. These coefficients were given in Tables~\ref{tab:dd} and~\ref{tab:ddOmega}.

\AB{This part of the generating function can also be written as $\gft(I,\phi') = C_2(I)\cos\left(2\phi'-2\phi'_0\right)$, a cosine function of period $\pi$ at angle $\phi'$. The corresponding Standard Scattering Map (SSM) can thus be used as a universal local model around saddle-center libration points of the  RTBP problem, and its cosine expression explains the \MR{bimodal} shape of the KAM curves in Figure~\ref{fig:phase_port_SM}.}

%
%

\subsection{Global Approximation Error}\label{sec:globalapprox}

Suppose now that we are now interested in an accurate representation in the whole domain of the scattering map
\MB{
\begin{equation} \label{eq:domainI} \domain = \{ (I, \phi)\colon\ I\in[1,7] \text{ and } \phi\in[0,2\pi) \}.\end{equation}
}

%

\begin{table}
	
	\scalebox{0.8}{
		\begin{filecontents*}{FT_error.csv}
N,zero,one,two,three,four,five,six,seven
2,0.680000,0.963414,0.495653,0.261042,0.385496,0.131751,0.219163,0.094550
4,0.680000,0.975219,0.618115,0.276085,0.387552,0.087317,0.156552,0.021123
6,0.680000,0.996139,0.656461,0.266668,0.381078,0.107612,0.144998,0.013382
8,0.680000,1.003501,0.657922,0.287596,0.367730,0.110813,0.148937,0.013159
10,0.680000,1.007927,0.647908,0.287781,0.356098,0.105846,0.149545,0.012512
\end{filecontents*}
\csvreader[
		head to column names,
		tabular = lcccccccc,
		table head = \toprule
		& \multicolumn{7}{c}{L} \\
		\cmidrule(r){2-9}
		$N$ & 0 & 1 & 2 & 3 & 4 & 5 & 6 & 7 \\
		\midrule,
		table foot = \bottomrule,
		]{FT_error.csv}{}{%
			\N & \zero & \one & \two & \three & \four & \five & \six & \seven	
		}%
	}
	
	\scalebox{0.8}{
		\begin{filecontents*}{T_error.csv}
N,zero,one,two,three,four,five,six
2,0.449508,0.268271,0.157380,0.065119,0.148788,0.068206,0.121809
4,0.449508,0.270152,0.169874,0.074112,0.150692,0.085074,0.132709
6,0.449508,0.269731,0.169564,0.076822,0.153570,0.083003,0.132714
8,0.449508,0.269683,0.169068,0.079125,0.155540,0.083914,0.131502
10,0.449508,0.269748,0.169407,0.080213,0.157050,0.085329,0.131450
\end{filecontents*}
\csvreader[
		head to column names,
		tabular = lccccccc,
		table head = \toprule
		& \multicolumn{6}{c}{L} \\
		\cmidrule(r){2-8}
		$N$ & 0 & 1 & 2 & 3 & 4 & 5 & 6 \\
		\midrule,
		table foot = \bottomrule,
		]{T_error.csv}{}{%
			\N & \zero & \one & \two & \three & \four & \five & \six
		}%
	}	
	\caption{Global approximation error  $\epsilon_I$ (top panel) and $\epsilon_\phi$ (bottom panel) as a function of $N,L$ for the standard scattering map $\sm_1$. The error is evaluated in the original data set $I=1,2,\dotsc,7$.}
	\label{tab:FT_error}

\end{table}

\begin{table}
\PB{
    \scalebox{0.8}{
		\begin{filecontents*}{FT_newerror.csv}
N,zero,one,two,three,four,five,six,seven
2,0.680000,0.782085,0.367800,0.198887,0.233913,0.071975,0.104571,0.077636
4,0.680000,0.822172,0.465488,0.201108,0.227997,0.035848,0.062120,0.032394
6,0.680000,0.841491,0.497539,0.205914,0.224239,0.048805,0.052237,0.022763
8,0.680000,0.848342,0.499856,0.219921,0.216637,0.050972,0.053874,0.023362
10,0.680000,0.852520,0.491977,0.219181,0.209180,0.049182,0.054193,0.023105
\end{filecontents*}
\csvreader[
		head to column names,
		tabular = lcccccccc,
		table head = \toprule
		& \multicolumn{7}{c}{L} \\
		\cmidrule(r){2-9}
		$N$ & 0 & 1 & 2 & 3 & 4 & 5 & 6 & 7 \\
		\midrule,
		table foot = \bottomrule,
		]{FT_newerror.csv}{}{%
			\N & \zero & \one & \two & \three & \four & \five & \six & \seven	
		}%
	}
}

    \PB{
	\scalebox{0.8}{
		\begin{filecontents*}{T_newerror.csv}
N,zero,one,two,three,four,five,six
2,0.370667,0.210922,0.112123,0.063086,0.107801,0.040273,0.067627
4,0.370667,0.213108,0.123612,0.057792,0.107593,0.040863,0.062113
6,0.370667,0.212718,0.122866,0.059556,0.109339,0.039459,0.062510
8,0.370667,0.212655,0.122537,0.059490,0.110579,0.039398,0.061885
10,0.370667,0.212716,0.122826,0.059643,0.111556,0.039633,0.061842
\end{filecontents*}
\csvreader[
		head to column names,
		tabular = lccccccc,
		table head = \toprule
		& \multicolumn{6}{c}{L} \\
		\cmidrule(r){2-8}
		$N$ & 0 & 1 & 2 & 3 & 4 & 5 & 6 \\
		\midrule,
		table foot = \bottomrule,
		]{T_newerror.csv}{}{%
			\N & \zero & \one & \two & \three & \four & \five & \six
		}%
	}	
    }
	\caption{Global approximation error  $\epsilon_I$ (top panel) and $\epsilon_\phi$ (bottom panel) as a function of $N,L$ for the standard scattering map $\sm_1$. The error is evaluated in the independent data set $I=0.5, 1.5, 2.5,\dotsc,6.5$. Compare to the original data set in Table~\ref{tab:FT_error}.}
	\label{tab:FT_newerror}
\end{table}

The approximation error over the global domain is given in
Table~\ref{tab:FT_error}.
If we want an approximation error $\epsilon = \max\{\epsilon_I, \epsilon_\phi\}$ less than $0.1$, we need to increase the degree of Fourier-\AB{Polynomial interpolation} to $N=4$ and $L=5$.
In fact, Figure~\ref{fig:check_FourierTaylor_N4_M5} shows that $N=4,\ L=5$ gives a good approximation in the global domain $\domain$.

From now on, we will use $N=4,\ L=5$ as our model for the global scattering map $\sm_1$.

A similar analysis suggests that $N=4,\ L=6$ be used for the global scattering map $\sm_2$.

%
%
%
%

\section{Inner Map}
\label{sec:inner_map}
Recall that the \textit{inner flow} refers to the restriction of the RTBP flow to the normally hyperbolic invariant manifold $\Lambda_c$, while the \textit{inner map} refers to the restriction of the first return map $\innerMap$ to $\NHIMsec$. Abusing notation, the inner map will still be called $\innerMap$.

As explained in Section~\ref{Sec:RTBP}, the RTBP flow on the center manifold $W^{\mathrm{c}}(L_1)$ consists simply on a translation on the $2$-torus,
\begin{align*}
\dot J_p &= 0, &\quad \dot \phi_p &= \pd{H}{J_p} \eqqcolon \nu_p(J_p, J_v), \\
\dot J_v &= 0, &\quad \dot \phi_v &= \pd{H}{J_v} \eqqcolon \nu_v(J_p, J_v).
\end{align*}
The value of the planar and vertical frequencies $\nu_p$, $\nu_v$ of the torus are obtained differentiating the Hamiltonian in Birkhoff normal form.

\begin{figure}
	\begin{tikzpicture}[gnuplot]
\path (0.000,0.000) rectangle (12.500,8.750);
\gpcolor{color=gp lt color border}
\gpsetlinetype{gp lt border}
\gpsetdashtype{gp dt solid}
\gpsetlinewidth{1.00}
\draw[gp path] (1.688,0.985)--(1.868,0.985);
\draw[gp path] (11.947,0.985)--(11.767,0.985);
\node[gp node right] at (1.504,0.985) {$6.094$};
\draw[gp path] (1.688,2.228)--(1.868,2.228);
\draw[gp path] (11.947,2.228)--(11.767,2.228);
\node[gp node right] at (1.504,2.228) {$6.096$};
\draw[gp path] (1.688,3.470)--(1.868,3.470);
\draw[gp path] (11.947,3.470)--(11.767,3.470);
\node[gp node right] at (1.504,3.470) {$6.098$};
\draw[gp path] (1.688,4.713)--(1.868,4.713);
\draw[gp path] (11.947,4.713)--(11.767,4.713);
\node[gp node right] at (1.504,4.713) {$6.1$};
\draw[gp path] (1.688,5.956)--(1.868,5.956);
\draw[gp path] (11.947,5.956)--(11.767,5.956);
\node[gp node right] at (1.504,5.956) {$6.102$};
\draw[gp path] (1.688,7.198)--(1.868,7.198);
\draw[gp path] (11.947,7.198)--(11.767,7.198);
\node[gp node right] at (1.504,7.198) {$6.104$};
\draw[gp path] (1.688,8.441)--(1.868,8.441);
\draw[gp path] (11.947,8.441)--(11.767,8.441);
\node[gp node right] at (1.504,8.441) {$6.106$};
\draw[gp path] (1.688,0.985)--(1.688,1.165);
\draw[gp path] (1.688,8.441)--(1.688,8.261);
\node[gp node center] at (1.688,0.677) {$0$};
\draw[gp path] (3.154,0.985)--(3.154,1.165);
\draw[gp path] (3.154,8.441)--(3.154,8.261);
\node[gp node center] at (3.154,0.677) {$1$};
\draw[gp path] (4.619,0.985)--(4.619,1.165);
\draw[gp path] (4.619,8.441)--(4.619,8.261);
\node[gp node center] at (4.619,0.677) {$2$};
\draw[gp path] (6.085,0.985)--(6.085,1.165);
\draw[gp path] (6.085,8.441)--(6.085,8.261);
\node[gp node center] at (6.085,0.677) {$3$};
\draw[gp path] (7.550,0.985)--(7.550,1.165);
\draw[gp path] (7.550,8.441)--(7.550,8.261);
\node[gp node center] at (7.550,0.677) {$4$};
\draw[gp path] (9.016,0.985)--(9.016,1.165);
\draw[gp path] (9.016,8.441)--(9.016,8.261);
\node[gp node center] at (9.016,0.677) {$5$};
\draw[gp path] (10.481,0.985)--(10.481,1.165);
\draw[gp path] (10.481,8.441)--(10.481,8.261);
\node[gp node center] at (10.481,0.677) {$6$};
\draw[gp path] (11.947,0.985)--(11.947,1.165);
\draw[gp path] (11.947,8.441)--(11.947,8.261);
\node[gp node center] at (11.947,0.677) {$7$};
\draw[gp path] (1.688,8.441)--(1.688,0.985)--(11.947,0.985)--(11.947,8.441)--cycle;
\node[gp node center,rotate=-270] at (0.292,4.713) {$\nu$};
\node[gp node center] at (6.817,0.215) {$I$};
\gpcolor{rgb color={0.580,0.000,0.827}}
\draw[gp path] (1.688,7.667)--(3.154,6.677)--(4.619,5.702)--(6.085,4.743)--(7.550,3.797)%
  --(9.016,2.865)--(10.481,1.946)--(11.947,1.039);
\gpsetpointsize{4.00}
\gppoint{gp mark 1}{(1.688,7.667)}
\gppoint{gp mark 1}{(3.154,6.677)}
\gppoint{gp mark 1}{(4.619,5.702)}
\gppoint{gp mark 1}{(6.085,4.743)}
\gppoint{gp mark 1}{(7.550,3.797)}
\gppoint{gp mark 1}{(9.016,2.865)}
\gppoint{gp mark 1}{(10.481,1.946)}
\gppoint{gp mark 1}{(11.947,1.039)}
\gpcolor{color=gp lt color border}
\draw[gp path] (1.688,8.441)--(1.688,0.985)--(11.947,0.985)--(11.947,8.441)--cycle;
\gpdefrectangularnode{gp plot 1}{\pgfpoint{1.688cm}{0.985cm}}{\pgfpoint{11.947cm}{8.441cm}}
\end{tikzpicture}
	\caption{The inner shift $\nu(I)$.}
	\label{fig:FreqAll}
\end{figure}
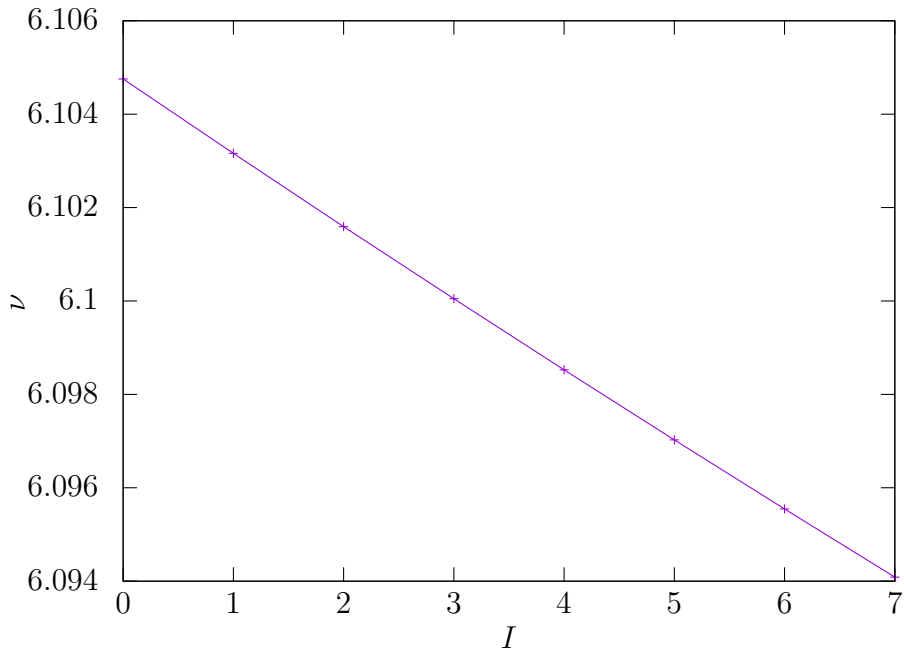

Upon restriction to the NHIM $\Lambda_c$, we get rid of the planar action $J_p$, which can be recovered if necessary using the energy condition, so the inner flow is
\begin{equation} \label{eq:InnerFlow}
\dot J_v = 0,
\quad \dot \phi_p = \nu_p(J_v),
\quad \dot \phi_v = \nu_v(J_v).
\end{equation}
Finally, the first return map of the inner flow to the section $\Sigma$ is
\[ J_v' = J_v,
\quad \phi_v' = \phi_v + \frac{2\pi\nu_v(J_v)}{\nu_p(J_v)}. \]
In terms of the scaled coordinates $I, \phi$, the inner map $\innerMap\colon\ \NHIMsecn\to \NHIMsecn$ is given by
\begin{equation} \label{eq:InnerMap}
I' = I, \quad \phi' = \phi + \nu(I),
\end{equation}
where we have introduced the new function
\[\nu(I) \coloneqq \frac{2\pi\nu_v(I/1000)}{\nu_p(I/1000)}.
\]
As seen in Figure~\ref{fig:FreqAll}, the inner shift $\nu(I)$ decreases almost linearly with $I$. In particular, this shows that the inner map $\innerMap$ is a twist map.

\section{Phase Space of the Scattering Map} \label{Sec:PhaseSpaceSM}

In Sections~\ref{Sec:SeriesRepresentation} and~\ref{sec:ApproximationError}, we have obtained a series representation of the scattering maps $\sm(I,\phi)$ with small approximation error. This representation is valid on an annulus $\domain$ inside the NHIM (see Equation~\eqref{eq:domainI}).
Thus we can now iterate the scattering map in $\domain$.

\begin{figure}
	%
	
	
	\includegraphics[width=0.8\linewidth]{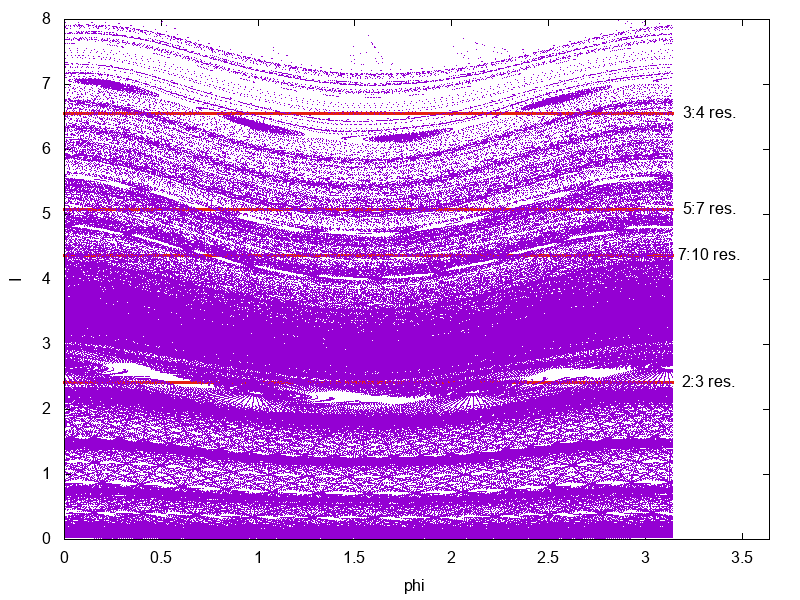}
	\includegraphics[width=0.8\linewidth]{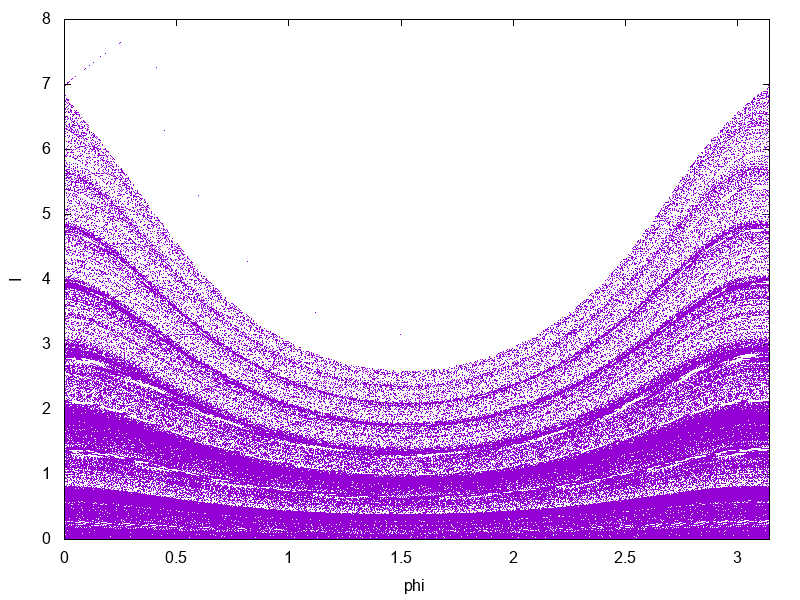}
	\caption{\PB{Phase portrait of the scattering maps $\sm_1$ (above) and $\sm_2$ (below). The theoretical location of the main resonances, found in Section~\ref{sec:KAM}, is superimposed in red.}}
	\label{fig:phase_port_SM}
\end{figure}

To explore the global phase space of the scattering map, we perform the following experiment: 
\PR{take $100$ initial conditions equi-distributed on the line $\left\{I\in[0,7],\ \phi=0\right\}$, and iterate each initial condition 1000 times by the scattering map just to draw smooth enough curves from continuous invariant sets.}
\PB{Take $300\times 300$ initial conditions evenly distributed in $\domain$, and iterate each initial condition 100 times by the scattering map.}
The resulting phase space portrait is shown in Figure~\ref{fig:phase_port_SM}.

Notice that the scattering map~\eqref{eq:sm}, rewritten as 
\begin{subequations} \label{eq:dampedSM}
	\begin{align}
	\phi' &= \phi - \omega(I) - \pd{\gft}{I}(I,\phi') \label{eq:dampedSM:phi}\\
	I' &= I + \pd{\gft}{\phi'}(I,\phi') \label{eq:dampedSM:I}
	\end{align}
\end{subequations}
can be seen as a perturbation of the integrable map $(I,\phi)\mapsto (I'=I,\phi'=\phi-\omega(I))$, as long as the derivatives of $\gft$ are small enough.
The frequency $-\omega(I)$ represents the \emph{phase shift} of the map  $(I,\phi)\mapsto (I'=I,\phi'=\phi-\omega(I))$. 

Since $\omega'(I)\neq 0$ (see Figure~\ref{fig:interp_omega}), the integrable map $(I,\phi)\mapsto (I,\phi-\omega(I))$ is a twist map. Below we check that the scattering map is also a twist map. (See `Twist Condition').

By KAM theory (for an area preserving map given by its generating function; see, for instance, \cite{ArnoldAvez68,Haro02}), when the derivatives of $\gft$ are small, we expect that many of the invariant tori of the integrable twist map persist.


\subsubsection*{Phase Shift}

\begin{figure}
	\scalebox{0.9}{\begin{tikzpicture}[gnuplot]
\path (0.000,0.000) rectangle (12.500,8.750);
\gpcolor{color=gp lt color border}
\gpsetlinetype{gp lt border}
\gpsetdashtype{gp dt solid}
\gpsetlinewidth{1.00}
\draw[gp path] (1.196,0.985)--(1.376,0.985);
\draw[gp path] (11.947,0.985)--(11.767,0.985);
\node[gp node right] at (1.012,0.985) {$-2.5$};
\draw[gp path] (1.196,2.228)--(1.376,2.228);
\draw[gp path] (11.947,2.228)--(11.767,2.228);
\node[gp node right] at (1.012,2.228) {$-2.4$};
\draw[gp path] (1.196,3.470)--(1.376,3.470);
\draw[gp path] (11.947,3.470)--(11.767,3.470);
\node[gp node right] at (1.012,3.470) {$-2.3$};
\draw[gp path] (1.196,4.713)--(1.376,4.713);
\draw[gp path] (11.947,4.713)--(11.767,4.713);
\node[gp node right] at (1.012,4.713) {$-2.2$};
\draw[gp path] (1.196,5.956)--(1.376,5.956);
\draw[gp path] (11.947,5.956)--(11.767,5.956);
\node[gp node right] at (1.012,5.956) {$-2.1$};
\draw[gp path] (1.196,7.198)--(1.376,7.198);
\draw[gp path] (11.947,7.198)--(11.767,7.198);
\node[gp node right] at (1.012,7.198) {$-2$};
\draw[gp path] (1.196,8.441)--(1.376,8.441);
\draw[gp path] (11.947,8.441)--(11.767,8.441);
\node[gp node right] at (1.012,8.441) {$-1.9$};
\draw[gp path] (1.196,0.985)--(1.196,1.165);
\draw[gp path] (1.196,8.441)--(1.196,8.261);
\node[gp node center] at (1.196,0.677) {$1$};
\draw[gp path] (2.988,0.985)--(2.988,1.165);
\draw[gp path] (2.988,8.441)--(2.988,8.261);
\node[gp node center] at (2.988,0.677) {$2$};
\draw[gp path] (4.780,0.985)--(4.780,1.165);
\draw[gp path] (4.780,8.441)--(4.780,8.261);
\node[gp node center] at (4.780,0.677) {$3$};
\draw[gp path] (6.572,0.985)--(6.572,1.165);
\draw[gp path] (6.572,8.441)--(6.572,8.261);
\node[gp node center] at (6.572,0.677) {$4$};
\draw[gp path] (8.363,0.985)--(8.363,1.165);
\draw[gp path] (8.363,8.441)--(8.363,8.261);
\node[gp node center] at (8.363,0.677) {$5$};
\draw[gp path] (10.155,0.985)--(10.155,1.165);
\draw[gp path] (10.155,8.441)--(10.155,8.261);
\node[gp node center] at (10.155,0.677) {$6$};
\draw[gp path] (11.947,0.985)--(11.947,1.165);
\draw[gp path] (11.947,8.441)--(11.947,8.261);
\node[gp node center] at (11.947,0.677) {$7$};
\draw[gp path] (1.196,8.441)--(1.196,0.985)--(11.947,0.985)--(11.947,8.441)--cycle;
\node[gp node center] at (6.571,0.215) {$I$};
\node[gp node right] at (10.479,8.107) {$-\omega(I)$};
\gpcolor{rgb color={0.580,0.000,0.827}}
\gpsetdashtype{gp dt 2}
\draw[gp path] (10.663,8.107)--(11.579,8.107);
\draw[gp path] (1.196,6.853)--(1.375,6.795)--(1.554,6.737)--(1.734,6.679)--(1.913,6.621)%
  --(2.092,6.564)--(2.271,6.506)--(2.450,6.448)--(2.629,6.390)--(2.809,6.332)--(2.988,6.274)%
  --(3.167,6.215)--(3.346,6.156)--(3.525,6.096)--(3.705,6.036)--(3.884,5.976)--(4.063,5.915)%
  --(4.242,5.853)--(4.421,5.791)--(4.600,5.728)--(4.780,5.665)--(4.959,5.601)--(5.138,5.536)%
  --(5.317,5.470)--(5.496,5.404)--(5.676,5.337)--(5.855,5.269)--(6.034,5.201)--(6.213,5.132)%
  --(6.392,5.061)--(6.572,4.990)--(6.751,4.918)--(6.930,4.846)--(7.109,4.772)--(7.288,4.697)%
  --(7.467,4.622)--(7.647,4.545)--(7.826,4.467)--(8.005,4.389)--(8.184,4.309)--(8.363,4.228)%
  --(8.543,4.146)--(8.722,4.062)--(8.901,3.978)--(9.080,3.892)--(9.259,3.804)--(9.438,3.716)%
  --(9.618,3.625)--(9.797,3.533)--(9.976,3.440)--(10.155,3.344)--(10.334,3.247)--(10.514,3.148)%
  --(10.693,3.047)--(10.872,2.944)--(11.051,2.838)--(11.230,2.731)--(11.409,2.620)--(11.589,2.507)%
  --(11.768,2.392)--(11.947,2.273);
\gpcolor{rgb color={0.000,0.620,0.451}}
\gpsetdashtype{gp dt solid}
\draw[gp path] (1.196,6.139)--(1.375,6.079)--(1.554,6.021)--(1.734,5.964)--(1.913,5.909)%
  --(2.092,5.855)--(2.271,5.802)--(2.450,5.750)--(2.629,5.699)--(2.809,5.648)--(2.988,5.597)%
  --(3.167,5.546)--(3.346,5.495)--(3.525,5.444)--(3.705,5.393)--(3.884,5.341)--(4.063,5.289)%
  --(4.242,5.236)--(4.421,5.182)--(4.600,5.123)--(4.780,5.063)--(4.959,5.003)--(5.138,4.942)%
  --(5.317,4.880)--(5.496,4.818)--(5.676,4.755)--(5.855,4.692)--(6.034,4.628)--(6.213,4.564)%
  --(6.392,4.500)--(6.572,4.436)--(6.751,4.371)--(6.930,4.307)--(7.109,4.243)--(7.288,4.180)%
  --(7.467,4.117)--(7.647,4.055)--(7.826,3.994)--(8.005,3.933)--(8.184,3.874)--(8.363,3.816)%
  --(8.543,3.759)--(8.722,3.705)--(8.901,3.652)--(9.080,3.590)--(9.259,3.520)--(9.438,3.448)%
  --(9.618,3.374)--(9.797,3.298)--(9.976,3.219)--(10.155,3.095)--(10.334,2.931)--(10.514,2.757)%
  --(10.693,2.574)--(10.872,2.381)--(11.051,2.180)--(11.230,1.967)--(11.409,1.744)--(11.589,1.510)%
  --(11.768,1.264)--(11.947,1.006);
\gpcolor{color=gp lt color border}
\node[gp node right] at (10.479,7.799) {$-\omega(I) \pm \max_{\phi\prime} \abs{\pd{\gft}{I}(I,\phi\prime)}$};
\gpcolor{rgb color={0.000,0.620,0.451}}
\draw[gp path] (10.663,7.799)--(11.579,7.799);
\draw[gp path] (1.196,7.568)--(1.375,7.511)--(1.554,7.453)--(1.734,7.394)--(1.913,7.334)%
  --(2.092,7.272)--(2.271,7.209)--(2.450,7.146)--(2.629,7.082)--(2.809,7.016)--(2.988,6.950)%
  --(3.167,6.884)--(3.346,6.816)--(3.525,6.748)--(3.705,6.680)--(3.884,6.610)--(4.063,6.540)%
  --(4.242,6.470)--(4.421,6.399)--(4.600,6.333)--(4.780,6.266)--(4.959,6.198)--(5.138,6.130)%
  --(5.317,6.061)--(5.496,5.990)--(5.676,5.919)--(5.855,5.847)--(6.034,5.774)--(6.213,5.699)%
  --(6.392,5.623)--(6.572,5.545)--(6.751,5.465)--(6.930,5.384)--(7.109,5.300)--(7.288,5.214)%
  --(7.467,5.126)--(7.647,5.035)--(7.826,4.941)--(8.005,4.844)--(8.184,4.744)--(8.363,4.640)%
  --(8.543,4.532)--(8.722,4.420)--(8.901,4.304)--(9.080,4.194)--(9.259,4.089)--(9.438,3.983)%
  --(9.618,3.877)--(9.797,3.769)--(9.976,3.661)--(10.155,3.593)--(10.334,3.564)--(10.514,3.540)%
  --(10.693,3.520)--(10.872,3.506)--(11.051,3.497)--(11.230,3.494)--(11.409,3.497)--(11.589,3.505)%
  --(11.768,3.519)--(11.947,3.540);
\gpcolor{color=gp lt color border}
\draw[gp path] (1.196,8.441)--(1.196,0.985)--(11.947,0.985)--(11.947,8.441)--cycle;
\gpdefrectangularnode{gp plot 1}{\pgfpoint{1.196cm}{0.985cm}}{\pgfpoint{11.947cm}{8.441cm}}
\end{tikzpicture}
}
	\scalebox{0.9}{\begin{tikzpicture}[gnuplot]
\path (0.000,0.000) rectangle (12.500,8.750);
\gpcolor{color=gp lt color border}
\gpsetlinetype{gp lt border}
\gpsetdashtype{gp dt solid}
\gpsetlinewidth{1.00}
\draw[gp path] (1.196,0.985)--(1.376,0.985);
\draw[gp path] (11.947,0.985)--(11.767,0.985);
\node[gp node right] at (1.012,0.985) {$-3.7$};
\draw[gp path] (1.196,2.228)--(1.376,2.228);
\draw[gp path] (11.947,2.228)--(11.767,2.228);
\node[gp node right] at (1.012,2.228) {$-3.6$};
\draw[gp path] (1.196,3.470)--(1.376,3.470);
\draw[gp path] (11.947,3.470)--(11.767,3.470);
\node[gp node right] at (1.012,3.470) {$-3.5$};
\draw[gp path] (1.196,4.713)--(1.376,4.713);
\draw[gp path] (11.947,4.713)--(11.767,4.713);
\node[gp node right] at (1.012,4.713) {$-3.4$};
\draw[gp path] (1.196,5.956)--(1.376,5.956);
\draw[gp path] (11.947,5.956)--(11.767,5.956);
\node[gp node right] at (1.012,5.956) {$-3.3$};
\draw[gp path] (1.196,7.198)--(1.376,7.198);
\draw[gp path] (11.947,7.198)--(11.767,7.198);
\node[gp node right] at (1.012,7.198) {$-3.2$};
\draw[gp path] (1.196,8.441)--(1.376,8.441);
\draw[gp path] (11.947,8.441)--(11.767,8.441);
\node[gp node right] at (1.012,8.441) {$-3.1$};
\draw[gp path] (1.196,0.985)--(1.196,1.165);
\draw[gp path] (1.196,8.441)--(1.196,8.261);
\node[gp node center] at (1.196,0.677) {$1$};
\draw[gp path] (2.988,0.985)--(2.988,1.165);
\draw[gp path] (2.988,8.441)--(2.988,8.261);
\node[gp node center] at (2.988,0.677) {$2$};
\draw[gp path] (4.780,0.985)--(4.780,1.165);
\draw[gp path] (4.780,8.441)--(4.780,8.261);
\node[gp node center] at (4.780,0.677) {$3$};
\draw[gp path] (6.572,0.985)--(6.572,1.165);
\draw[gp path] (6.572,8.441)--(6.572,8.261);
\node[gp node center] at (6.572,0.677) {$4$};
\draw[gp path] (8.363,0.985)--(8.363,1.165);
\draw[gp path] (8.363,8.441)--(8.363,8.261);
\node[gp node center] at (8.363,0.677) {$5$};
\draw[gp path] (10.155,0.985)--(10.155,1.165);
\draw[gp path] (10.155,8.441)--(10.155,8.261);
\node[gp node center] at (10.155,0.677) {$6$};
\draw[gp path] (11.947,0.985)--(11.947,1.165);
\draw[gp path] (11.947,8.441)--(11.947,8.261);
\node[gp node center] at (11.947,0.677) {$7$};
\draw[gp path] (1.196,8.441)--(1.196,0.985)--(11.947,0.985)--(11.947,8.441)--cycle;
\node[gp node center] at (6.571,0.215) {$I$};
\node[gp node right] at (10.479,1.935) {$-\omega(I)$};
\gpcolor{rgb color={0.580,0.000,0.827}}
\gpsetdashtype{gp dt 2}
\draw[gp path] (10.663,1.935)--(11.579,1.935);
\draw[gp path] (1.196,2.775)--(1.375,2.826)--(1.554,2.879)--(1.734,2.935)--(1.913,2.992)%
  --(2.092,3.051)--(2.271,3.110)--(2.450,3.170)--(2.629,3.230)--(2.809,3.291)--(2.988,3.352)%
  --(3.167,3.413)--(3.346,3.474)--(3.525,3.535)--(3.705,3.596)--(3.884,3.658)--(4.063,3.720)%
  --(4.242,3.782)--(4.421,3.844)--(4.600,3.907)--(4.780,3.970)--(4.959,4.033)--(5.138,4.097)%
  --(5.317,4.162)--(5.496,4.228)--(5.676,4.294)--(5.855,4.361)--(6.034,4.429)--(6.213,4.498)%
  --(6.392,4.568)--(6.572,4.639)--(6.751,4.711)--(6.930,4.784)--(7.109,4.859)--(7.288,4.934)%
  --(7.467,5.010)--(7.647,5.087)--(7.826,5.165)--(8.005,5.244)--(8.184,5.325)--(8.363,5.406)%
  --(8.543,5.488)--(8.722,5.571)--(8.901,5.655)--(9.080,5.741)--(9.259,5.827)--(9.438,5.915)%
  --(9.618,6.003)--(9.797,6.094)--(9.976,6.186)--(10.155,6.280)--(10.334,6.375)--(10.514,6.473)%
  --(10.693,6.574)--(10.872,6.678)--(11.051,6.786)--(11.230,6.898)--(11.409,7.015)--(11.589,7.137)%
  --(11.768,7.265)--(11.947,7.401);
\gpcolor{rgb color={0.000,0.620,0.451}}
\gpsetdashtype{gp dt solid}
\draw[gp path] (1.196,1.167)--(1.375,1.242)--(1.554,1.320)--(1.734,1.401)--(1.913,1.484)%
  --(2.092,1.567)--(2.271,1.652)--(2.450,1.736)--(2.629,1.821)--(2.809,1.906)--(2.988,1.990)%
  --(3.167,2.074)--(3.346,2.157)--(3.525,2.240)--(3.705,2.323)--(3.884,2.405)--(4.063,2.487)%
  --(4.242,2.569)--(4.421,2.652)--(4.600,2.729)--(4.780,2.805)--(4.959,2.883)--(5.138,2.961)%
  --(5.317,3.040)--(5.496,3.121)--(5.676,3.205)--(5.855,3.290)--(6.034,3.378)--(6.213,3.468)%
  --(6.392,3.560)--(6.572,3.656)--(6.751,3.754)--(6.930,3.855)--(7.109,3.959)--(7.288,4.065)%
  --(7.467,4.175)--(7.647,4.286)--(7.826,4.400)--(8.005,4.516)--(8.184,4.633)--(8.363,4.752)%
  --(8.543,4.871)--(8.722,4.991)--(8.901,5.112)--(9.080,5.231)--(9.259,5.349)--(9.438,5.466)%
  --(9.618,5.580)--(9.797,5.689)--(9.976,5.794)--(10.155,5.897)--(10.334,5.998)--(10.514,6.097)%
  --(10.693,6.187)--(10.872,6.266)--(11.051,6.334)--(11.230,6.391)--(11.409,6.432)--(11.589,6.457)%
  --(11.768,6.465)--(11.947,6.453);
\gpcolor{color=gp lt color border}
\node[gp node right] at (10.479,1.627) {$-\omega(I) \pm \max_{\phi\prime} \abs{\pd{\gft}{I}(I,\phi\prime)}$};
\gpcolor{rgb color={0.000,0.620,0.451}}
\draw[gp path] (10.663,1.627)--(11.579,1.627);
\draw[gp path] (1.196,4.382)--(1.375,4.410)--(1.554,4.439)--(1.734,4.469)--(1.913,4.501)%
  --(2.092,4.534)--(2.271,4.568)--(2.450,4.603)--(2.629,4.639)--(2.809,4.675)--(2.988,4.713)%
  --(3.167,4.751)--(3.346,4.790)--(3.525,4.830)--(3.705,4.870)--(3.884,4.911)--(4.063,4.952)%
  --(4.242,4.994)--(4.421,5.036)--(4.600,5.084)--(4.780,5.134)--(4.959,5.184)--(5.138,5.234)%
  --(5.317,5.284)--(5.496,5.334)--(5.676,5.384)--(5.855,5.432)--(6.034,5.481)--(6.213,5.529)%
  --(6.392,5.576)--(6.572,5.623)--(6.751,5.669)--(6.930,5.713)--(7.109,5.758)--(7.288,5.802)%
  --(7.467,5.845)--(7.647,5.887)--(7.826,5.930)--(8.005,5.973)--(8.184,6.016)--(8.363,6.059)%
  --(8.543,6.104)--(8.722,6.151)--(8.901,6.199)--(9.080,6.250)--(9.259,6.305)--(9.438,6.363)%
  --(9.618,6.426)--(9.797,6.499)--(9.976,6.578)--(10.155,6.662)--(10.334,6.753)--(10.514,6.850)%
  --(10.693,6.961)--(10.872,7.091)--(11.051,7.238)--(11.230,7.405)--(11.409,7.597)--(11.589,7.816)%
  --(11.768,8.066)--(11.947,8.349);
\gpcolor{color=gp lt color border}
\node[gp node right] at (10.479,1.319) {$-\pi$};
\gpcolor{rgb color={0.902,0.624,0.000}}
\draw[gp path] (10.663,1.319)--(11.579,1.319);
\draw[gp path] (1.196,7.924)--(1.305,7.924)--(1.413,7.924)--(1.522,7.924)--(1.630,7.924)%
  --(1.739,7.924)--(1.848,7.924)--(1.956,7.924)--(2.065,7.924)--(2.173,7.924)--(2.282,7.924)%
  --(2.391,7.924)--(2.499,7.924)--(2.608,7.924)--(2.716,7.924)--(2.825,7.924)--(2.934,7.924)%
  --(3.042,7.924)--(3.151,7.924)--(3.259,7.924)--(3.368,7.924)--(3.477,7.924)--(3.585,7.924)%
  --(3.694,7.924)--(3.802,7.924)--(3.911,7.924)--(4.019,7.924)--(4.128,7.924)--(4.237,7.924)%
  --(4.345,7.924)--(4.454,7.924)--(4.562,7.924)--(4.671,7.924)--(4.780,7.924)--(4.888,7.924)%
  --(4.997,7.924)--(5.105,7.924)--(5.214,7.924)--(5.323,7.924)--(5.431,7.924)--(5.540,7.924)%
  --(5.648,7.924)--(5.757,7.924)--(5.866,7.924)--(5.974,7.924)--(6.083,7.924)--(6.191,7.924)%
  --(6.300,7.924)--(6.409,7.924)--(6.517,7.924)--(6.626,7.924)--(6.734,7.924)--(6.843,7.924)%
  --(6.952,7.924)--(7.060,7.924)--(7.169,7.924)--(7.277,7.924)--(7.386,7.924)--(7.495,7.924)%
  --(7.603,7.924)--(7.712,7.924)--(7.820,7.924)--(7.929,7.924)--(8.038,7.924)--(8.146,7.924)%
  --(8.255,7.924)--(8.363,7.924)--(8.472,7.924)--(8.581,7.924)--(8.689,7.924)--(8.798,7.924)%
  --(8.906,7.924)--(9.015,7.924)--(9.124,7.924)--(9.232,7.924)--(9.341,7.924)--(9.449,7.924)%
  --(9.558,7.924)--(9.666,7.924)--(9.775,7.924)--(9.884,7.924)--(9.992,7.924)--(10.101,7.924)%
  --(10.209,7.924)--(10.318,7.924)--(10.427,7.924)--(10.535,7.924)--(10.644,7.924)--(10.752,7.924)%
  --(10.861,7.924)--(10.970,7.924)--(11.078,7.924)--(11.187,7.924)--(11.295,7.924)--(11.404,7.924)%
  --(11.513,7.924)--(11.621,7.924)--(11.730,7.924)--(11.838,7.924)--(11.947,7.924);
\gpcolor{color=gp lt color border}
\draw[gp path] (1.196,8.441)--(1.196,0.985)--(11.947,0.985)--(11.947,8.441)--cycle;
\gpdefrectangularnode{gp plot 1}{\pgfpoint{1.196cm}{0.985cm}}{\pgfpoint{11.947cm}{8.441cm}}
\end{tikzpicture}
}
	\caption{Enclosure of the phase shift for the scattering map $\sm_1$ (above) and $\sm_2$ (below). The phase shift $-\omega(I) - \pd{\gft}{I}(I,\phi')$ is enclosed inside the green lines. }
	\label{fig:phase_shift}
\end{figure}
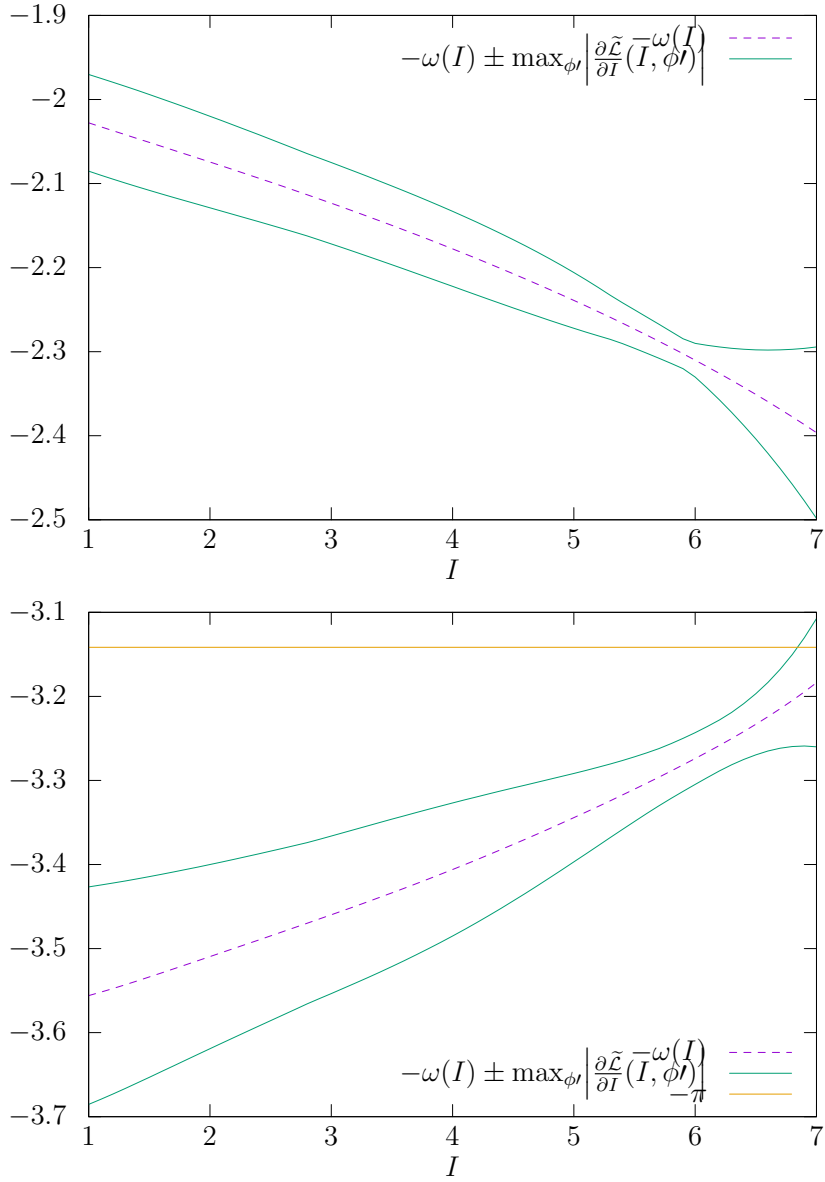

The phase shift of the scattering map~\eqref{eq:dampedSM} is $\phi'-\phi = - \omega(I) - \pd{\gft}{I}(I,\phi')$.
The divided differences of $\omega(I)$ are listed in Table~\ref{tab:ddOmega}.
From these data we see that $|\tilde{c}_0|\gg |\tilde{c}_1|\gg |\tilde{c}_j|$ for $j=2,\dots,6$. Comparing with Table~\ref{tab:dd}, we notice that $|\tilde{c}_0|\gg |\tilde{a}^n_j|, |\tilde{b}^n_j|$, $j=1,\dots,7$, $n=2,4,6,8$. This implies that $\displaystyle \frac{\partial \gft}{\partial I}$ is much smaller than $\omega(I)$, which is non-zero, at least for $I$ small enough. This argument could be used to assert that the phase shift is non-zero in the local domain $\domainloc$.

To deal with the global domain $\domain$, we determine the range of $I$ values where the phase shift is non-zero. The phase shift is bounded by
\[ -\omega(I) - \max_{\phi'} \abs{\pd{\gft}{I}(I,\phi')} \leq 
-\omega(I) - \pd{\gft}{I}(I,\phi') \leq 
-\omega(I) + \max_{\phi'} \abs{\pd{\gft}{I}(I,\phi')}.\]
We have computed these bounds explicitly, using the series expansions of $\omega(I)$ and $\pd{\gft}{I}(I,\phi')$. The result is shown in Figure~\ref{fig:phase_shift}. Note that the phase shift for $\sm_1$ is non-zero for all $I\in[1,7]$, while the phase shift for $\sm_2$ (modulo $\pi$) is non-zero except possibly for a small range of $I$ values close to $I=7$.

\begin{figure}
	\includegraphics[width=\linewidth]{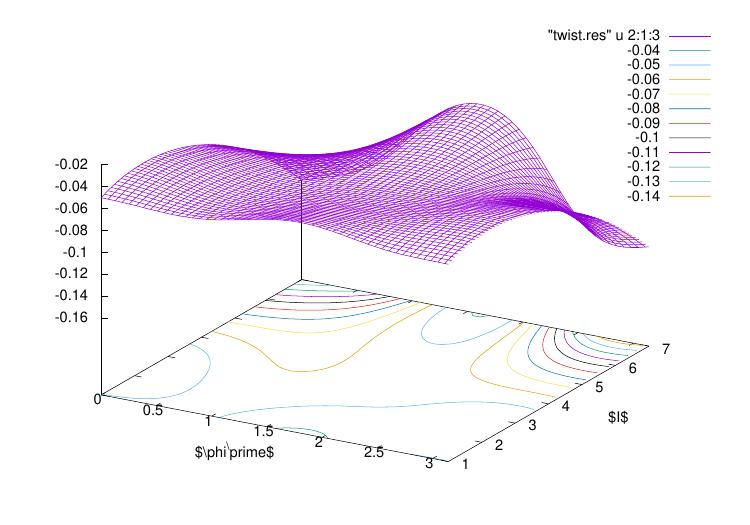}
	\includegraphics[width=\linewidth]{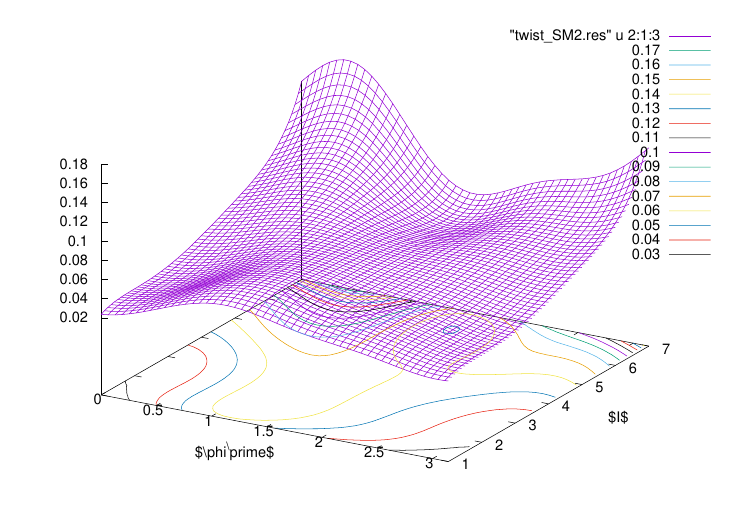}
	
	\caption{Twist for the scattering map $\sm_1$ (above) and $\sm_2$ (below).}
	\label{fig:twist}
\end{figure}

\subsubsection*{Twist Condition}

%

Let us compute the twist of the scattering map~\eqref{eq:dampedSM}:
\[ \pd{\phi'}{I}(I,\phi') = - \omega'(I) - \pd[2]{\gft}{I} (I,\phi') - \md{\gft}{2}{\phi'}{}{I}{} (I,\phi') \pd{\phi'}{I}. \]
Thus
\begin{equation} \label{eq:twist}
\pd{\phi'}{I}(I,\phi') = - \frac{\omega'(I) + \pd[2]{\gft}{I}}{1 + \md{\gft}{2}{\phi'}{}{I}{} },
\end{equation}
provided that the denominator is non-zero.

We have computed the twist~\eqref{eq:twist} explicitly on the global domain $\domain$, using the series expansions of $\omega(I)$ and $\gft(I,\phi')$. The result is shown in Figure~\ref{fig:twist}. Note that the twist for both $\sm_1$ and $\sm_2$ is non-zero in $\domain$.

\section{KAM tori and resonant zones for the Scattering Maps}\label{sec:KAM}
As seen in Figure~\ref{fig:phase_port_SM}, the phase portrait of the scattering maps $
\sm_1, \sm_2$ is filled up with invariant curves and
some resonant zones between them. 
Recall from Section~\ref{Sec:PhaseSpaceSM} that both $\sm_1$ and $\sm_2$ are twist maps.
To compute these invariant curves, we now introduce the exponential form in the sine-cosine Fourier expansion~\eqref{eq:dgft/dphi} of $ \pd{\gft}{\phi'}$
\[ \pd{\gft}{\phi'} =
\sum_{n=1}^N A_n\cos{n\phi'}+\sum_{n=1}^N B_n\sin{n\phi'}= \sum_{n=-N,n\neq 0}^N C_n \ee^{in\phi'},
\]
where $C=\left(C_{-n},\dots,C_{-1},C_{1},\dots,C_{N}\right)\in\mathbb{C}^{2N}$ satisfies
\[
C_n=\frac{1}{2}(A_n-i B_n), \quad C_{-n}=\frac{1}{2}(A_n+i B_n)=\overline{C_n}, \text{ for } n>0.
\]

Notice that for $\gft=0$, or equivalently $C=0$, any torus $I=I_0$ is invariant since then $I'=I_0$ in the expression~\eqref{eq:dampedSM:I} of a scattering map, with an inner dynamics $\phi'=\phi-\omega_0$ given by~\eqref{eq:dampedSM:phi}, where $\omega_0:=\omega(I_0)$.

For $\left|\gft\right|$ small enough, or equivalently $|C|$ small enough, a lot of these invariant curves survive.
An invariant curve $I=I_0+h(\phi)$ of a scattering map~\eqref{eq:dampedSM} satisfies $I'=I_0+h(\phi')$, that is
\[
h(\phi')=h\left(\phi'+\omega(I)+\frac{\partial\gft}{\partial I}(I,\phi')\right)
+ \pd{\gft}{\phi'}(I,\phi'), \text{ where } I=I_0+h(\phi).
\]
Expanding in $C$ and $h$ we get
\[
h(\phi')=h(\phi'+\omega_0)+\pd{\gft}{\phi'}(I,\phi') +O(h\,C).
\]
Writing $ h(\phi')=\sum_{n=-N,n\neq 0}^N h_n \ee^{in\phi'}$ we get
\[
\sum h_n \ee^{in\phi'}=\sum\ee^{in\omega_0} h_n \ee^{in\phi'} + \sum C_n \ee^{in\phi'} + O(h\,C),
\]
which, equating Fourier coefficients, gives
\[
h_n =-\frac{C_n}{\ee^{in\omega_0}-1}+O(C^2)  \text{ for } 0<|n|\leq N.
\]
For this approximate formula one needs that $\dfrac{n\omega_0}{2\pi}\not\in\mathbb{Z}$ for $|n|\leq N$ (non-resonant condition),
and one sees that, up to order $O(C^2)$, the coefficients $h_n=O(C)$ are uniquely determined by $C$. KAM theorem consists 
in proving the convergence of these expansions for \emph{diophantine frequencies} $\omega_0$, using that $\sm_1$ and $\sm_2$ are twist maps.



\subsubsection*{Resonant Zones for $\sm_1$}

Resonant zones for the symplectic map~\eqref{eq:sm}, or, equivalently,~\eqref{eq:dampedSM},
where the KAM theorem does not provide invariant curves for small $\mathcal{L}$, appear around the values $I$
such that
$\displaystyle\frac{\omega(I)}{\pi}$ is a rational number. In the global setting,
\[
\omega(I)=\Omega'(I)\approx \sum_{l=0}^{L}c_l N_l(I),
\]
where the coefficients $c_l$ are given in table~\ref{tab:ddOmega} for both scattering maps $\sigma_1$ and $\sigma_2$. 

In particular for $\sigma_1$ (a totally analogous study for $\sigma_2$ can be carried out) $\displaystyle\frac{\omega(0)}{\pi}=0.630128\dots$, which is not too far for $2/3$ whose continued fraction is [1,2], which means that
\[
\frac{2}{3}=[\mathbf{1},\mathbf{2}]=\frac{1}{\mathbf{1}+\displaystyle\frac{1}{\mathbf{2}}}.
\]
Therefore for $I$ such that $\displaystyle\frac{\omega(I)}{\pi}=\frac{2}{3}$, which happens to be $I\approx 2.4175$, there should appear a resonance,
indeed the largest one, since the width of the `eyes' of a resonance is related to the denominator, in this case $3$.

%
%
%
%

Other close rationals to 2/3 are given by close modified continued fractions. For instance
\[
[\mathbf{1},\mathbf{2},\mathbf{1}]=\frac{1}{\mathbf{1}+\displaystyle\frac{1}{\mathbf{2} + \displaystyle\frac{1}{\mathbf{1}} }}=\frac{3}{4}.
\]

We can compute some of the largest ones, ordered by their denominators:

\centerline{\renewcommand{\arraystretch}{1.5}
	\begin{tabular}{l|cl}
		\multicolumn{1}{c}{$I$}&$\displaystyle\frac{\omega(I)}{\pi}$&continued fraction\\
		\hline
		2.4175&$\frac{2}{3}$&[1,2]\\
		6.5550&$\frac{3}{4}$&[1,2,1]\\
		5.0752&$\frac{5}{7}$&[1,2,2]\\
		4.3631&$\frac{7}{10}$&[1,2,3]\\
		3.9523&$\frac{9}{13}$&[1,2,4]
\end{tabular}}

The continued fraction of $\displaystyle\frac{2}{3}$ can be also written as [1,1,1], and smaller resonant values can be obtained for smaller continued fractions
like $\displaystyle [1,1,1,3]=\frac{7}{11} = 0.636363\dots$, etc

For the scattering map $\sm_1$, the two main resonances are clearly visible near $I\approx 2.4175$ and $I\approx 6.5550$. 
See Figure~\ref{fig:phase_port_SM} (top panel).

\section{\MB{Arnold diffusion}}
\label{sec:HeuristicArnoldDiffusion}

\MB{In this section we give details on the verification of the Numerical Result~\ref{thm:main_numerical}. Specifically, we verify numerically that the conditions of Theorem \ref{th:main} are satisfied.} 

\MB{Let \[\domain=\{(I,\phi)\,|\, I\in [1,7],\, \phi\in [0,2\pi)\} \] be the annulus where we want to show diffusion. The inner map
$f$ is  $\innerMap$  restricted to $\domain$, and the system of scattering maps is either one of $\{\sigma_1\}$, $\{\sigma_2\}$, or 
$\{\sigma_1,\sigma_2\}$.
Recall that $f$ is a twist map. In Section~\ref{Sec:PhaseSpaceSM} we verified numerically that $\sigma_1$ and $\sigma_2$ are twist maps as well.
}

As seen in Table~\ref{tab:dd}, the harmonics of ${\gft}_i$ for both scattering maps $\sigma_1$ and $\sigma_2$, particularly those of degree two, are not zero, so that the inner map and any of these scattering maps can not have common invariant curves, and by Theorem~\ref{th:main}  both the double dynamical systems $\left\{\innerMap, \sm_1\right\}$ and $\left\{\innerMap, \sm_2\right\}$ formed by the inner map and one of the scattering maps have diffusing  pseudo-orbits along $\domain$.

Even more, for the two scattering maps described in  Table~\ref{tab:dd}, one sees that the difference between the respective coefficients $C_n$ is greater than $0.04$ for the coefficients $C_{\pm 1}$ due to $\tilde{ a}_1$. This, together with the fact that the magnitude of the denominator $\ee^{\pm i\omega_0}-1$ in the formula of $h_{\pm 1}$ is much smaller than 1, prevents the two scattering maps from having common invariant curves, which, on the other hand, is clearly observed in the juxtaposition of the curves found numerically for the two scattering maps.
This implies that the double  dynamical system $\left\{\sm_1, \sm_2\right\}$ also has diffusing orbits along the NHIM.

We can take advantage of these dynamics to construct fast diffusing pseudo-orbits obtained from the triple dynamical system $\left\{\innerMap, \sm_1, \sm_2\right\}$ formed by the inner map and the two scattering maps.
We will give explicit constructions  of diffusing pseudo-orbits, including fast ones, in Section~\ref{Sec:DriftOrbits}.

Once such pseudo-orbits are obtained, the Shadowing Lemma \cite[Theorem 3.7]{GideaLlaveSeara20-CPAM} gives true orbits that shadow the pseudo-orbits, thus achieving Arnold diffusion.

Note that the above construction of diffusing pseudo-orbits assumes that the inner dynamics is derived from the  Birkhoff normal form approximation, which is given by an integrable Hamiltonian.
However, the original Hamiltonian is not integrable, and the NHIM $\Lambda^\Sigma_c$  for the original Hamiltonian  is not foliated by circles invariant under  the true inner map.
We will now argue the existence of diffusing orbits for the original system.

Recall that the inner map for the Birkhoff normal form is an integrable twist map, and  the global error in the numerical integration of orbits with initial condition $(J,\phi_p,\phi_v)$  is less than $10^{-12}/5000$; see Section \ref{Sec:NumericalResults}.
(We recall $J=J_v$ is the vertical amplitude of the motion.)
Since each level set of $J$ is preserved by the inner dynamics $\innerMap$ for the Birkhoff normal form,
and since $I = 1000J$, it
follows that each essential invariant circle $\Gamma$ for the inner map $\innerMap$ for the original
Hamiltonian is less than $\rho_1=10^{-12}<10^{-9}/1500$ away from a level set of $I$.

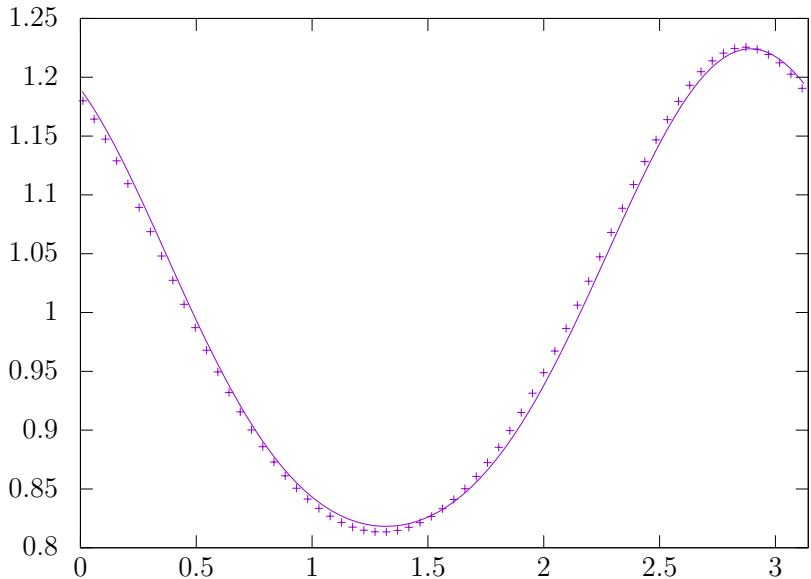
\begin{figure}
	\scalebox{0.9}{\begin{tikzpicture}[gnuplot]
\path (0.000,0.000) rectangle (12.500,8.750);
\gpcolor{color=gp lt color border}
\gpsetlinetype{gp lt border}
\gpsetdashtype{gp dt solid}
\gpsetlinewidth{1.00}
\draw[gp path] (1.196,0.616)--(1.376,0.616);
\draw[gp path] (11.947,0.616)--(11.767,0.616);
\node[gp node right] at (1.012,0.616) {$0.8$};
\draw[gp path] (1.196,1.485)--(1.376,1.485);
\draw[gp path] (11.947,1.485)--(11.767,1.485);
\node[gp node right] at (1.012,1.485) {$0.85$};
\draw[gp path] (1.196,2.355)--(1.376,2.355);
\draw[gp path] (11.947,2.355)--(11.767,2.355);
\node[gp node right] at (1.012,2.355) {$0.9$};
\draw[gp path] (1.196,3.224)--(1.376,3.224);
\draw[gp path] (11.947,3.224)--(11.767,3.224);
\node[gp node right] at (1.012,3.224) {$0.95$};
\draw[gp path] (1.196,4.094)--(1.376,4.094);
\draw[gp path] (11.947,4.094)--(11.767,4.094);
\node[gp node right] at (1.012,4.094) {$1$};
\draw[gp path] (1.196,4.963)--(1.376,4.963);
\draw[gp path] (11.947,4.963)--(11.767,4.963);
\node[gp node right] at (1.012,4.963) {$1.05$};
\draw[gp path] (1.196,5.833)--(1.376,5.833);
\draw[gp path] (11.947,5.833)--(11.767,5.833);
\node[gp node right] at (1.012,5.833) {$1.1$};
\draw[gp path] (1.196,6.702)--(1.376,6.702);
\draw[gp path] (11.947,6.702)--(11.767,6.702);
\node[gp node right] at (1.012,6.702) {$1.15$};
\draw[gp path] (1.196,7.572)--(1.376,7.572);
\draw[gp path] (11.947,7.572)--(11.767,7.572);
\node[gp node right] at (1.012,7.572) {$1.2$};
\draw[gp path] (1.196,8.441)--(1.376,8.441);
\draw[gp path] (11.947,8.441)--(11.767,8.441);
\node[gp node right] at (1.012,8.441) {$1.25$};
\draw[gp path] (1.196,0.616)--(1.196,0.796);
\draw[gp path] (1.196,8.441)--(1.196,8.261);
\node[gp node center] at (1.196,0.308) {$0$};
\draw[gp path] (2.907,0.616)--(2.907,0.796);
\draw[gp path] (2.907,8.441)--(2.907,8.261);
\node[gp node center] at (2.907,0.308) {$0.5$};
\draw[gp path] (4.618,0.616)--(4.618,0.796);
\draw[gp path] (4.618,8.441)--(4.618,8.261);
\node[gp node center] at (4.618,0.308) {$1$};
\draw[gp path] (6.329,0.616)--(6.329,0.796);
\draw[gp path] (6.329,8.441)--(6.329,8.261);
\node[gp node center] at (6.329,0.308) {$1.5$};
\draw[gp path] (8.040,0.616)--(8.040,0.796);
\draw[gp path] (8.040,8.441)--(8.040,8.261);
\node[gp node center] at (8.040,0.308) {$2$};
\draw[gp path] (9.751,0.616)--(9.751,0.796);
\draw[gp path] (9.751,8.441)--(9.751,8.261);
\node[gp node center] at (9.751,0.308) {$2.5$};
\draw[gp path] (11.462,0.616)--(11.462,0.796);
\draw[gp path] (11.462,8.441)--(11.462,8.261);
\node[gp node center] at (11.462,0.308) {$3$};
\draw[gp path] (1.196,8.441)--(1.196,0.616)--(11.947,0.616)--(11.947,8.441)--cycle;
\gpcolor{rgb color={0.580,0.000,0.827}}
\gpsetpointsize{4.00}
\gppoint{gp mark 1}{(1.233,7.224)}
\gppoint{gp mark 1}{(1.399,6.953)}
\gppoint{gp mark 1}{(1.565,6.657)}
\gppoint{gp mark 1}{(1.731,6.336)}
\gppoint{gp mark 1}{(1.897,5.998)}
\gppoint{gp mark 1}{(2.063,5.647)}
\gppoint{gp mark 1}{(2.229,5.289)}
\gppoint{gp mark 1}{(2.395,4.929)}
\gppoint{gp mark 1}{(2.561,4.570)}
\gppoint{gp mark 1}{(2.727,4.215)}
\gppoint{gp mark 1}{(2.893,3.870)}
\gppoint{gp mark 1}{(3.059,3.536)}
\gppoint{gp mark 1}{(3.225,3.216)}
\gppoint{gp mark 1}{(3.391,2.912)}
\gppoint{gp mark 1}{(3.557,2.625)}
\gppoint{gp mark 1}{(3.723,2.358)}
\gppoint{gp mark 1}{(3.889,2.110)}
\gppoint{gp mark 1}{(4.055,1.884)}
\gppoint{gp mark 1}{(4.221,1.679)}
\gppoint{gp mark 1}{(4.387,1.496)}
\gppoint{gp mark 1}{(4.553,1.336)}
\gppoint{gp mark 1}{(4.719,1.198)}
\gppoint{gp mark 1}{(4.885,1.083)}
\gppoint{gp mark 1}{(5.051,0.991)}
\gppoint{gp mark 1}{(5.217,0.921)}
\gppoint{gp mark 1}{(5.383,0.875)}
\gppoint{gp mark 1}{(5.548,0.852)}
\gppoint{gp mark 1}{(5.715,0.851)}
\gppoint{gp mark 1}{(5.881,0.873)}
\gppoint{gp mark 1}{(6.047,0.919)}
\gppoint{gp mark 1}{(6.212,0.987)}
\gppoint{gp mark 1}{(6.379,1.079)}
\gppoint{gp mark 1}{(6.545,1.193)}
\gppoint{gp mark 1}{(6.711,1.330)}
\gppoint{gp mark 1}{(6.876,1.489)}
\gppoint{gp mark 1}{(7.043,1.671)}
\gppoint{gp mark 1}{(7.209,1.875)}
\gppoint{gp mark 1}{(7.375,2.101)}
\gppoint{gp mark 1}{(7.540,2.348)}
\gppoint{gp mark 1}{(7.707,2.615)}
\gppoint{gp mark 1}{(7.873,2.901)}
\gppoint{gp mark 1}{(8.039,3.204)}
\gppoint{gp mark 1}{(8.204,3.524)}
\gppoint{gp mark 1}{(8.371,3.858)}
\gppoint{gp mark 1}{(8.537,4.203)}
\gppoint{gp mark 1}{(8.703,4.557)}
\gppoint{gp mark 1}{(8.868,4.917)}
\gppoint{gp mark 1}{(9.035,5.277)}
\gppoint{gp mark 1}{(9.201,5.635)}
\gppoint{gp mark 1}{(9.366,5.986)}
\gppoint{gp mark 1}{(9.532,6.325)}
\gppoint{gp mark 1}{(9.699,6.646)}
\gppoint{gp mark 1}{(9.865,6.945)}
\gppoint{gp mark 1}{(10.030,7.215)}
\gppoint{gp mark 1}{(10.196,7.453)}
\gppoint{gp mark 1}{(10.363,7.654)}
\gppoint{gp mark 1}{(10.529,7.814)}
\gppoint{gp mark 1}{(10.694,7.929)}
\gppoint{gp mark 1}{(10.860,7.996)}
\gppoint{gp mark 1}{(11.027,8.015)}
\gppoint{gp mark 1}{(11.193,7.986)}
\gppoint{gp mark 1}{(11.358,7.909)}
\gppoint{gp mark 1}{(11.524,7.785)}
\gppoint{gp mark 1}{(11.691,7.617)}
\gppoint{gp mark 1}{(11.857,7.408)}
\draw[gp path] (1.224,7.358)--(1.317,7.223)--(1.412,7.075)--(1.507,6.916)--(1.603,6.746)%
  --(1.700,6.566)--(1.798,6.376)--(1.897,6.177)--(1.997,5.969)--(2.098,5.755)--(2.199,5.534)%
  --(2.303,5.308)--(2.407,5.077)--(2.512,4.843)--(2.619,4.606)--(2.727,4.369)--(2.836,4.131)%
  --(2.947,3.895)--(3.059,3.662)--(3.173,3.433)--(3.288,3.208)--(3.404,2.990)--(3.522,2.779)%
  --(3.641,2.576)--(3.761,2.383)--(3.883,2.200)--(4.006,2.028)--(4.130,1.867)--(4.255,1.719)%
  --(4.382,1.583)--(4.509,1.460)--(4.637,1.350)--(4.766,1.253)--(4.896,1.169)--(5.026,1.098)%
  --(5.157,1.040)--(5.288,0.995)--(5.420,0.962)--(5.551,0.941)--(5.683,0.933)--(5.814,0.937)%
  --(5.946,0.952)--(6.076,0.980)--(6.207,1.020)--(6.337,1.073)--(6.466,1.137)--(6.594,1.214)%
  --(6.721,1.302)--(6.848,1.403)--(6.973,1.517)--(7.097,1.642)--(7.220,1.778)--(7.342,1.926)%
  --(7.462,2.085)--(7.581,2.254)--(7.699,2.433)--(7.815,2.621)--(7.930,2.817)--(8.043,3.021)%
  --(8.155,3.231)--(8.266,3.446)--(8.375,3.666)--(8.483,3.889)--(8.590,4.114)--(8.695,4.341)%
  --(8.799,4.569)--(8.902,4.795)--(9.004,5.020)--(9.104,5.242)--(9.204,5.461)--(9.302,5.675)%
  --(9.400,5.884)--(9.497,6.086)--(9.593,6.282)--(9.688,6.470)--(9.782,6.649)--(9.876,6.820)%
  --(9.969,6.980)--(10.062,7.131)--(10.154,7.271)--(10.246,7.399)--(10.337,7.516)--(10.428,7.621)%
  --(10.519,7.714)--(10.609,7.793)--(10.700,7.860)--(10.790,7.913)--(10.880,7.953)--(10.970,7.979)%
  --(11.060,7.992)--(11.150,7.990)--(11.241,7.975)--(11.331,7.945)--(11.422,7.902)--(11.513,7.845)%
  --(11.604,7.774)--(11.696,7.690)--(11.788,7.592)--(11.881,7.482);
\gpcolor{color=gp lt color border}
\draw[gp path] (1.196,8.441)--(1.196,0.616)--(11.947,0.616)--(11.947,8.441)--cycle;
\gpdefrectangularnode{gp plot 1}{\pgfpoint{1.196cm}{0.616cm}}{\pgfpoint{11.947cm}{8.441cm}}
\end{tikzpicture}
}
	\caption{Image of the action level $I=1$ by the numerical scattering map (points) versus the scattering map series $N=4, M=5$ (lines).}
	\label{fig:check_FourierTaylor_I_1}
\end{figure}

We also know from Section \ref{subsubsec:GlobalizingSM} that the scattering maps $\sigma_1$, $\sigma_2$ are globally defined on $\domain$. Moreover, for each $I$ in the interval $[1,7]$, the oscillation, 
\[\sup_{\phi} \frac{\partial \gft_1}{\partial \phi'}(I,\phi')-\inf_{\phi} \frac{\partial \gft_1}{\partial \phi'}(I,\phi')\]
corresponding to $\sigma_1$, is bigger than $\rho_2=0.2$.

To see this, note that the smallest oscillations in Figure~\ref{fig:check_FourierTaylor_N4_M5} happen for the action level $I=1$. Zooming in that action level (Figure~\ref{fig:check_FourierTaylor_I_1}), it is clear that the oscillation is bigger than $\rho_2=0.2$.
Alternatively, approximate the image of $I=1$ by the truncated Fourier-\AB{Polynomial interpolation} consisting of the dominant coefficients
\[ I' = I + \tilde{a}_1^{(2)}\cos 2\phi' +
\tilde{b}_1^{(2)}\sin 2\phi', \]
where $\tilde{a}_1^{(2)} = 0.178180$ and $\tilde{b}_1^{(2)} = -0.097275$ (See Table~\ref{tab:dd}).
This function has oscillations of size twice its amplitude $\sqrt{\left(\tilde{a}_1^{(2)}\right)^2 + \left(\tilde{b}_1^{(2)}\right)^2} = 0.203003$.

Similarly, the oscillation of $\frac{\partial \gft_2}{\partial \phi'}$ corresponding to $\sigma_2$ is bigger than $0.3$.

These facts imply that no essential invariant circle $\Gamma$ for $\innerMap$ is  invariant under $\sigma_1$ or $\sigma_2$.
Therefore, Theorem \ref{th:main} applies and there are orbits of the IFS $\{\innerMap,\sigma_1\}$, as well as orbits  of the IFS $\{\innerMap,\sigma_2\}$, that go from the lower boundary of the annulus $\domain$ to its upper boundary.

Again, the Shadowing Lemma \cite[Theorem 3.7]{GideaLlaveSeara20-CPAM} gives true orbits that shadow the obtained pseudo-orbits, thus achieving Arnold diffusion.

\section{Time Estimates for Inner and Transition Map}\label{sec:TimeEstimates}

One of our main goals is to estimate the drift time spent by drift orbits, constructed in Section~\ref{Sec:DriftOrbits}. As an intermediate step, we measure the time spent on one iterate of the inner map (`inner time', or $\tin$), and the time spent on one iterate of the transition map (`outer time', or $\tout$).

In terms of the RTBP inner flow~\eqref{eq:InnerFlow}, one application of the inner map corresponds to integrating an initial condition $(J_v, \phi_p=0, \phi_v)\in \NHIMsec$ during the amount of time that it takes to return to the section $\Sigma$. Thus, each iterate of the inner map takes time
\[ \tin = \frac{2\pi}{\nu_p}. \]
Numerically, we find that $2.0764 < \nu_p(I) < 2.0781$, and therefore the inner time is bounded by
\[ 3.0235 < \tin < 3.0261. \]

In terms of the RTBP flow, one application of the transition map  corresponds to a (segment of) a homoclinic trajectory from $y_- \in W^u_\textrm{loc}(\NHIM_c)$ to $y_+ W^s_\textrm{loc}(\NHIM_c)$. As noted in Subsection~\ref{subsec:HomoclinicOrbits}, the flight time of homoclinic segments is bounded by
\[ 5.936738 \leq \tout \leq 6.000688. \]
Therefore in our setting we find that outer times $\tout$ are approximately twice as long as inner times $\tin$.

\section{Drift Orbits}\label{Sec:DriftOrbits}

\MR{In the previous section we provided some heuristic arguments for the existence of trajectories that undergo Arnold diffusion in the sense that they drift in terms of the  scaled  vertical amplitude $I$ (see \eqref{eq:scaled}).}
\MB{In Section~\ref{sec:HeuristicArnoldDiffusion} we provided numerical evidence that the conditions of Theorem~\ref{th:main} are fulfilled, hence there exist trajectories that exhibit Arnold diffusion, characterized by a drift in the scaled vertical amplitude 
$I$ (see \eqref{eq:scaled}).
Our arguments  do not provide a way to find such orbits, much less those orbits that drift `fast'.}

In this section, we propose different algorithms to produce drift orbits, i.e., orbits of the iterated function system (IFS) consisting of the inner and outer map, whose action variable $I$ increases from $I=1$ to $I>7$. We produce two different type of orbits:

\begin{itemize}
	\item Orbits of the iterated function system $\{\innerMap, \sm_1\}$, $\{\innerMap, \sm_2\}$, or $\{\innerMap, \sm_1, \sm_2\}$. These orbits are a realization of the existence \MB{Theorem~\ref{th:main}}. We look for short orbits, i.e. we try to minimize the number of iterates. However, it is important to realize that \textbf{these orbits do not directly translate to pseudo-orbits for the R3BP flow} \MR{. The reason is that the scattering map is not dynamically defined; it is geometrically defined instead (i.e.} (since one iterate of the scattering map does not correspond to a segment of homoclinic trajectory).
	
	\item Orbits of the iterated function system $\{\innerMap, \tm_1\}$, $\{\innerMap, \tm_2\}$, or $\{\innerMap, \tm_1, \tm_2\}$. These are also a realization of the existence \MB{Theorem~\ref{th:main}}, and they \textbf{directly translate to pseudo-orbits of the R3BP flow}. Each iterate of the transition map corresponds to a segment of homoclinic trajectory. Using the time estimates of Section~\ref{sec:TimeEstimates}, we can estimate the total drift time of the pseudo-orbit as:
	\[ t = n_0 \tin + n_1 \tout + n_2 \tout, \]
	where $n_0$, $n_1$ and $n_2$ denote the number of iterates of $\innerMap$, $\tm_1$ and $\tm_2$ respectively.
	Keeping an eye on Astrodynamics applications, we want to minimize the total drift time.
\end{itemize}

\subsection{Double Dynamical System \texorpdfstring{$\{\innerMap, \sm_i\}$}{(innerMap, sm i)}. Greedy Algorithm.}

In Section~\ref{sec:KAM} we already established that $\sm_i$ ($i=1,2$) has many invariant curves. All iterates of $\sm_i$ either belong to an invariant curve, or are confined between two invariant curves. Thus it is not possible to cross from $I=1$ to $I>7$ using just one scattering map $\sm_i$.

However, according to \MB{Theorem~\ref{th:main}}, one can combine the inner and outer map to produce drift orbits.
Now we will explicitly construct such drift orbits for the double dynamical system $\{\innerMap, \sm_i\}$.

Let us partition the domain $\domain$ of the scattering map into three sets:
\[ \domain = \negDomain \cup \neuDomain \cup \posDomain. \]
$\posDomain$ denotes the subdomain where $\sm$ gains action, $\negDomain$ where it looses action, and $\neuDomain$ where it neither gains nor looses action:
\begin{align*}
\posDomain &= \left\{ (I,\phi') \mid I'-I>0 \right\}
= \left\{ (I,\phi') \mid \pd{\gft}{\phi'}(I,\phi')>0 \right\} \\
\negDomain &= \left\{ (I,\phi') \mid I'-I<0 \right\}
= \left\{ (I,\phi') \mid \pd{\gft}{\phi'}(I,\phi')<0 \right\} \\
\neuDomain &= \left\{ (I,\phi') \mid I'-I=0 \right\}
= \left\{ (I,\phi') \mid \pd{\gft}{\phi'}(I,\phi')=0 \right\}.
\end{align*}
These sets are readily identified in Figure~\ref{fig:dL_dphip}. For the first scattering map, $\neuDomain$ roughly consists of two vertical lines at $\phi'\approx 0.5$ and $\phi'\approx 2$, and $\posDomain$ roughly consists of the vertical strip \MB{$(I,\phi') \in [1,7]\times(2, 0.5)$}.
For the second scattering map, $\neuDomain$ roughly consists of two vertical lines at $\phi'\approx 1.25$ and $\phi'\approx 3$, and $\posDomain$ roughly consists of the vertical strip \MB{$(I,\phi') \in [1,7]\times(3, 1.25)$}.

\begin{figure}
	\scalebox{0.9}{\input{dL_dphip}}
	\scalebox{0.9}{\input{dL_dphip_SM2}}
	\caption{The function $\pd{\gft}{\phi'}(I,\phi')$ for $\sm_1$ (above) and $\sm_2$ (below). Note that the function $\phi' \to \pd{\gft}{\phi'}(I_0,\phi')$ attains its maximum at
		$\phi'\approx 2.8$ (above) and $\phi'\approx 0.6$ (below).}
	\label{fig:dL_dphip}
\end{figure}

A simple strategy to produce drift orbits is to \textit{always apply the scattering map if it increases the action} (even if the action gain $I'-I$ is small). Otherwise, apply the inner map.

This `greedy' algorithm is guaranteed to produce a drift orbit independently of the initial condition, due to the following simple observations: In our model's domain $\domain$ of validity,

\begin{itemize}
	\item The inner map $(I', \phi') = \innerMap(I, \phi)$, given in Equation~\eqref{eq:InnerMap}, is a twist map with frequency $\nu(I) \approx 6.1$ on the universal cover (see Figure~\ref{fig:FreqAll}), or $\nu(I)\approx -0.2$ on the base space (where angles are identified modulo $\pi$). Thus the angle $\phi$ decreases approximately by $0.2$ radians at every iterate of the inner map.
	
	\item Hence, for any given point $(I, \phi) \in \domain$, its forward orbit by the inner map eventually enters $\posDomain$.
	
	\item If $(I,\phi)$ already belongs to $\posDomain$, we apply the scattering map, increasing the action. Else, we apply the inner map until the orbit enters $\posDomain$, and then apply the scattering map.
\end{itemize}

\begin{figure}
	\scalebox{0.85}{\includegraphics[width=\linewidth]{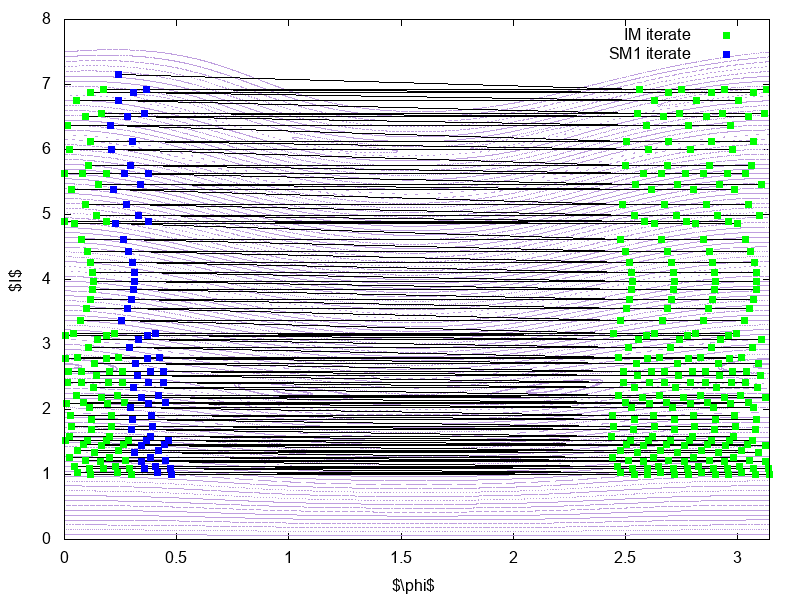}}
	\scalebox{0.85}{\includegraphics[width=\linewidth]{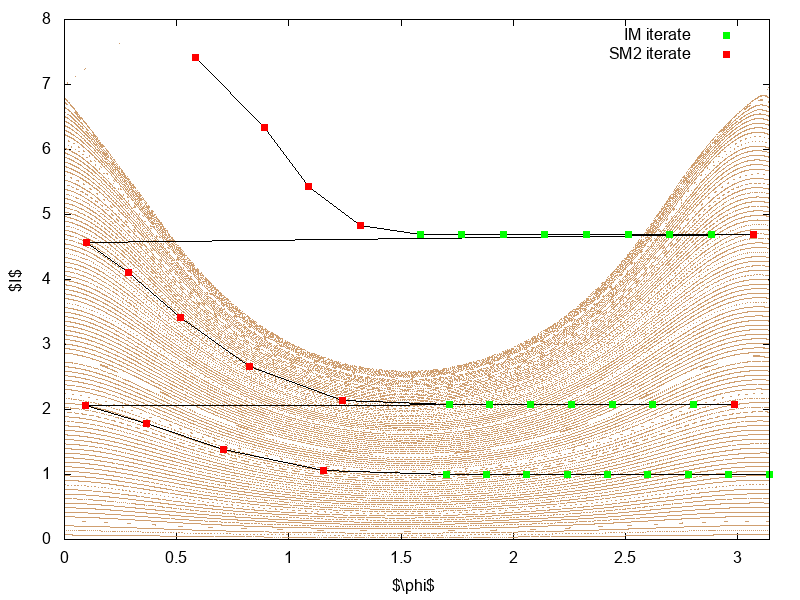}}
	\caption{Drift orbit of $\{\innerMap, \sm_1\}$ (top panel) and $\{\innerMap, \sm_2\}$ (bottom panel) using the greedy algorithm. Green points correspond to iterates of the inner map; blue (resp. red) points correspond to iterates of the first (resp. second) scattering map. Iterates have been joined by line segments to make the orbit more visible. For reference, the orbit is shown against a background consisting of the phase space of the scattering map.}
	\label{fig:diffusion_greedy}
\end{figure}

For example, Figure~\ref{fig:diffusion_greedy} shows the drift orbit produced starting from the initial condition $(I=1, \phi=0)$.
Notice that the drift orbit produced by $\{\innerMap, \sm_1\}$ is much longer than the one produced by $\{\innerMap, \sm_2\}$. However, as explained before, these orbits do not directly translate to pseudo-orbits of the RTBP flow, and we don't have control over their drift time.

\subsection{Triple Dynamical System \texorpdfstring{$\{\innerMap, \tm_1, \tm_2\}$}{(innerMap, tm1, tm2)}. Shortest-Time Algorithm.}\label{sec:triple_dynamical_system}

The algorithms described in previous sections are relatively simple to implement, but they yield sub-optimal pseudo-orbits in terms of their drift time. Now we focus on finding the optimal drift time. This is specially challenging when combining three dynamical systems (inner map, transition map 1, and transition map 2) to construct the pseudo-orbit. Obviously we can't consider all the possible combinations of $\{\innerMap, \tm_1, \tm_2\}$, since this number grows exponentially with respect to the length of the orbit.
The main idea is to leverage the classic Dijkstra algorithm~\cite{Dijkstra59} for finding shortest paths in a graph.

First we partition the domain $\domain$ into a uniform grid of $m \times n$ two-cells (rectangles) of equal size by dividing \MB{$I\in[1,7]$} into $m$ intervals and $\phi\in[0,\pi)$ into $n$ intervals. (In practice, we will use $m=n=30$, so the grid consists of $900$ small cells).

We introduce a \textit{directed} graph $G=<V,E>$ whose vertices $V$ represent the different cells.
An edge $e\in E$ from $u\in V$ to $v\in V$ means that the center point $(I,\phi)$ of cell $u$ is mapped into cell $v$ either by $\innerMap$, $\tm_1$, or $\tm_2$. In each case, the edge records the `distance' between cells, defined as the integration time corresponding to applying $\innerMap, \tm_1$, resp. $\tm_2$.

More precisely, an edge from $u$ to $v$ is a pair $e=(map, distance)$, where
\begin{itemize}
	\item $e=(\innerMap, \tin)$ if $(I,\phi)$ is mapped into $v$ by the inner map;
	\item $e=(\tm_1, \tout)$ if $(I,\phi)$ is mapped into $v$ by the first transition map;
	\item $e=(\tm_2, \tout)$ if $(I,\phi)$ is mapped into $v$ by the second transition map;
	\item $e=(\emptyset, \infty)$ if $(I,\phi)$ is not mapped into $v$ by neither map.
\end{itemize}

%

On rare occasions, $\tm_1$ and/or $\tm_2$ can map the center point $(I,\phi)$ to the same cell as the inner map. If this happens, the inner map is preferred since it takes shorter time, so we set the edge to $(\innerMap, \tin)$.

\begin{remark}
	The image $(I',\phi')$ of $(I,\phi)$ by the transition map may be outside $\domain$. However, this can only happen when $I'>7$ (by construction, $I'<0$ can never happen). In this case, we associate $u$ with the closest cell $v$ to the point $(I', \phi')$, namely the cell containing $(7, \phi')$.
\end{remark}


Given a source cell $s$ and a destination cell $t$, Dijkstra's algorithm applied to $G$ provides the \textbf{shortest (directed) path in the graph} from $s$ to $t$ in terms of the distance defined above.

Notice that this path does not exactly correspond to an orbit of the IFS, since we have only considered iterates of center points to construct $G$ (and the orbit does not necessarily pass through center points, but rather through arbitrary cell points). However, this path clearly \emph{informs} the choice of map $\{\innerMap, \tm_1, \tm_2\}$ that we should apply when the orbit passes through a given cell.

For example, suppose that the current iterate is inside cell $u$, and the shortest path from $u$ to $t$ starts with, say,
\[ u \xrightarrow{(\tm_1, \tout)} v \longrightarrow \cdots \longrightarrow t. \]
Then, the best choice given the available information is to apply the first transition map to the current iterate.

Our algorithm to construct optimal orbits (shortest drift time) is given next.

\begin{algorithm}
	\caption{Shortest-Time Algorithm}
	\label{shortest}
	\begin{algorithmic}[1] 
		\Procedure{OrbitShortestTime}{$x,y$} \Comment{Shortest-time orbit from point $x\in \domain$ to (a neighborhood of) point $y\in \domain$}
		\State $t\gets cell(y)$ \Comment{Destination cell}
		\State $orbit \gets x$ \Comment{Initialize orbit with $x$}
		\While{$x \notin neighborhood(y)$} \Comment{End when close enough to $y$}
		\State $u\gets cell(x)$ \Comment{Update current cell}
		\State $path\gets Dijkstra(u,t)$ \Comment{Shortest path from $u$ to $t$}
		\If{$path$ starts with $\innerMap$}
		\State $x \gets \innerMap(x)$
		\ElsIf{$path$ starts with $\tm_1$}
		\State $x \gets \tm_1(x)$
		\Else 	\Comment{$path$ starts with $\tm_2$}
		\State $x \gets \tm_2(x)$
		\EndIf
		\State $orbit \gets concat(orbit, x)$ 	\Comment{Add iterate $x$ to orbit}
		\EndWhile\label{shortestendwhile}
		\State \textbf{return} $orbit$
		\EndProcedure
	\end{algorithmic}
\end{algorithm}

\begin{figure}
	\scalebox{0.9}{\includegraphics[width=\linewidth]{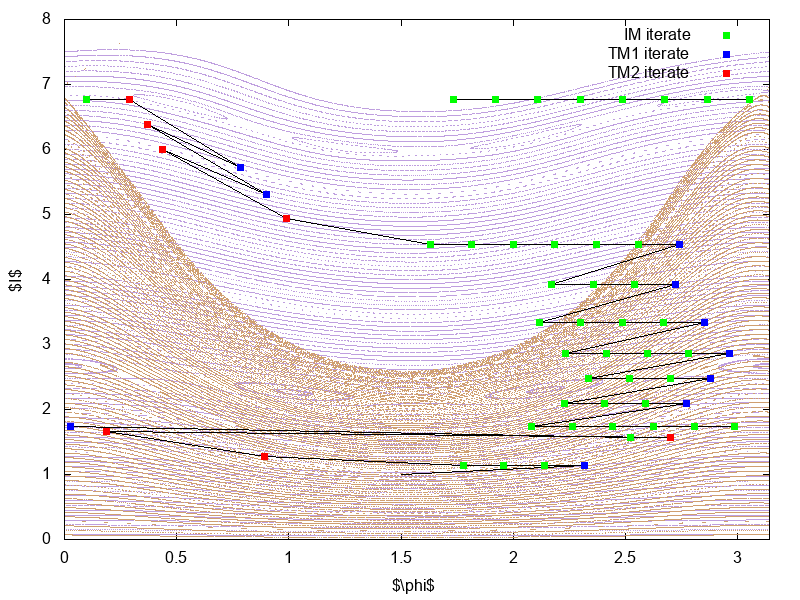}}
	\caption{Drift orbit of the triple system $\{\innerMap, \tm_1, \tm_2\}$ using the shortest path algorithm. Notice that some iterates of the transition map actually \textit{reduce} the action, in order to increase it more efficiently overall.}
	\label{fig:diffusion_shortest_path_TM}
\end{figure}

Figure~\ref{fig:diffusion_shortest_path_TM} shows the shortest-time orbit from $x=(I,\phi)=(1,1.5)$ to a neighborhood of $y=(I,\phi)=(7,1.5)$. The corresponding pseudo-orbit for the RTBP flow takes time $34\tin+17\tout \approx 204$ RTBP time units, i.e. about 32 years (optimal drift time).
Compare this to the orbits obtained in previous sections.

Notice that the optimal orbit uses all three dynamics ($\innerMap$, $\tm_1$ and $\tm_2$) for maximum flexibility.


Notice that some iterates actually \textit{decrease} the action. The key point is that, sometimes, one needs to take an iterate that decreases action in order to quickly move to a region where it later increases sharply. This way the pseudo-orbit's time is globally optimized.

\AR{Notice that in the last iterations of the scattering maps, it is when more vertical amplitude $I$ is gained (or lost), and that the scattering maps are only used 1 or 2 times near $I=7$, so it is not necessary to achieve much precision near the last invariant torus corresponding to $I=7$.}


Finally, from the orbit of the IFS $\{\innerMap, \tm_1, \tm_2\}$ we can construct a pseudo-orbit by concatenating segments of trajectories of the RTBP flow. Each iterate of the inner map $\innerMap$ corresponds to its flow suspension, which is integrated using the Birkhoff normal form. Each iterate of the transition map $\tm_1$ or $\tm_2$ corresponds to a finite piece of homoclinic trajectory, which is computed by continuation of those previously found in Section~\ref{Sec:NumericalResults}.

\begin{figure}
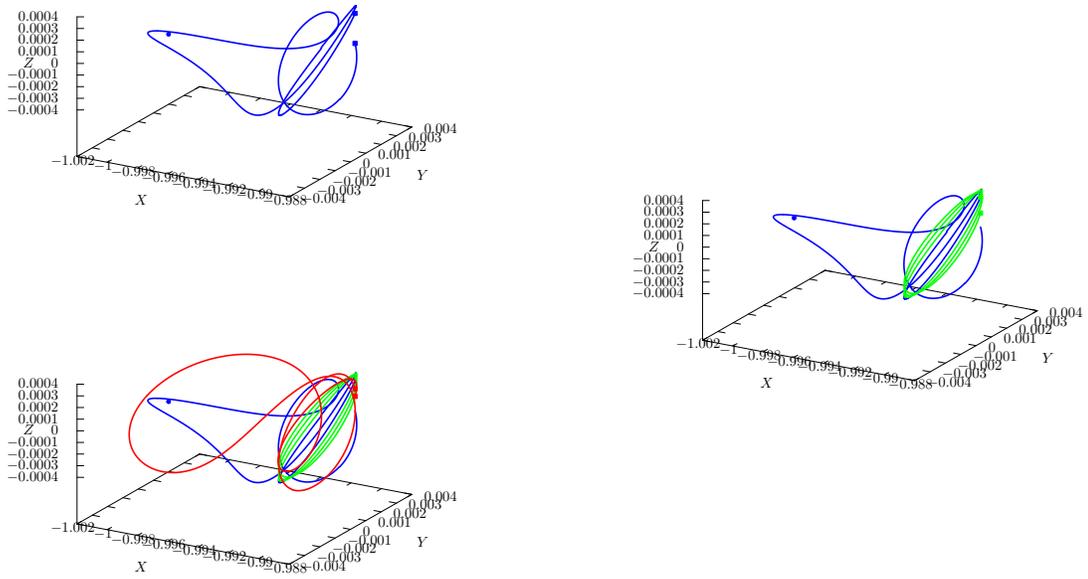

	\begin{minipage}{0.6\textwidth}
		\small
		\scalebox{.5}{\input{diffusion_orbit_1}}
		
		\hspace{0.5mm}
		
		\scalebox{.5}{\input{diffusion_orbit_3}}
	\end{minipage}
	\begin{minipage}{0.6\textwidth}
		\small
		\scalebox{.5}{\input{diffusion_orbit_2}}
	\end{minipage}
	\caption{Top left: Trajectory segment corresponding to the first application of $\tm_1$ in Figure~\ref{fig:diffusion_shortest_path_TM} (endpoints are marked with squares). The Earth is represented by a blue circle (not to scale). Middle right: Segments corresponding to $\tm_1 \circ \innerMap^3$. Bottom left: Segments corresponding to $\tm_1 \circ \innerMap^3 \circ \tm_2$. }
	\label{fig:diffusion_orbit}
\end{figure}

Figure~\ref{fig:diffusion_orbit} illustrates the construction of the pseudo-orbit corresponding to the orbit in Figure~\ref{fig:diffusion_shortest_path_TM}: First, $\tm_1$ is applied once; this corresponds to the blue homoclinic segment. Then $\innerMap$ is applied three times; this corresponds to the green segment. Next, $\tm_2$ is applied once; this corresponds to the homoclinic red segment. Notice that all segments start and end on the Poincar\'e section $\Sigma$ (endpoints are marked with squares). This construction would continue until the whole pseudo-orbit is obtained (not displayed).

\PB{
\begin{remark}\label{rem:thrust}
    We have obtained a pseudo-orbit, not a true trajectory of the RTBP. The endpoints of consecutive segments do not exactly match in positions or velocities, but the discontinuities are small (as is apparent in Figure~\ref{fig:diffusion_orbit}).
    In fact, the velocity discontinuities $\Delta v = \norm{(\Delta \dot X, \Delta \dot Y, \Delta \dot Z)}$, are smaller than $3.5\times 10^{-4}$.  Since the complete pseudo-orbit consists of 27 consecutive segments, the total required thrust is less than $26 \times 3.5\times 10^{-4} \approx 0.0091$ non-dimensional RTBP units. When converted into metric units, this amounts to $271$ m/s.
\end{remark}
}

Of course, to this pseudo-orbit \AB{one can apply shadowing results~\cite[Theorem 3.7]{GideaLlaveSeara20-CPAM}} to \MB{establish the existence of a}  true RTBP trajectory that shadows it.
For applications, however, obtaining the pseudo-orbit is often the crucial step, since it is ultimately refined in a much more realistic model than the RTBP, e.g. using JPL's Ephemeris.

%


\section*{Acknowledgement}{AD and PR supported by Spanish grant PID2021-123968NB-I00 (MI-CIU/AEI/10.13039/501100011033/FEDER/UE). Research of M.G. was partially supported by NSF grant DMS-2307718 and DMS-2154725.} 
\AB{
The authors are very grateful to the anonymous reviewers for their comments and suggestions that have contributed to improving this article in a very significant way.
}

\bibliographystyle{elsarticle-num} 
\bibliography{mybiblio}



%
%
%
\end{document}